\documentclass[11pt]{article}

\usepackage[margin=0.9in]{geometry}
\usepackage[english]{babel}
\usepackage{makecell}
\usepackage{amsthm,amsmath,amsfonts,amssymb,mathtools,diagbox,listings,verbatim}
\usepackage{mathrsfs}
\usepackage{bm}
\usepackage{bbm}
\usepackage{graphicx}
\usepackage{booktabs}
\usepackage{enumitem}
\usepackage{xcolor}
\usepackage{aliascnt}
\usepackage[numbers,sort&compress]{natbib}
\usepackage{extarrows}
\usepackage{tikz}
\usepackage{float}
\usepackage{pdfpages}
\usepackage{microtype}
\usepackage{titlesec}
\titlespacing*{\section}{0pt}{1.8ex plus .3ex minus .2ex}{0.9ex plus .2ex}
\titlespacing*{\subsection}{0pt}{1.45ex plus .25ex minus .15ex}{0.7ex plus .15ex}
\titlespacing*{\subsubsection}{0pt}{1.1ex plus .2ex minus .1ex}{0.5ex plus .1ex}
\setlist{topsep=4pt,itemsep=1pt,parsep=1pt,partopsep=0pt}
\usepackage[colorlinks,citecolor=blue,urlcolor=blue,linkcolor=blue]{hyperref}

\usetikzlibrary{arrows.meta,calc}
\tikzset{
  sepfig/.style={
    scale=0.88,
    every node/.style={font=\scriptsize},
    line cap=round,
    line join=round,
    >=Stealth
  },
  sepedge/.style={thick},
  separrow/.style={thick,-{Stealth[length=1.8mm,width=1.4mm]}},
  seplabel/.style={fill=white,inner sep=0.7pt}
}

\numberwithin{equation}{section}

\theoremstyle{plain}
\newtheorem{theorem}{Theorem}[section]

\newaliascnt{lemma}{theorem}
\newtheorem{lemma}[lemma]{Lemma}
\aliascntresetthe{lemma}

\newaliascnt{proposition}{theorem}
\newtheorem{proposition}[proposition]{Proposition}
\aliascntresetthe{proposition}

\newaliascnt{corollary}{theorem}

\aliascntresetthe{corollary}

\theoremstyle{definition}
\newaliascnt{definition}{theorem}

\aliascntresetthe{definition}

\newaliascnt{example}{theorem}

\aliascntresetthe{example}

\newaliascnt{remark}{theorem}
\newtheorem{remark}[remark]{Remark}
\aliascntresetthe{remark}

\newcommand{\EE}{\mathbb{E}}

\newcommand{\RR}{\mathbb{R}}
\newcommand{\ZZ}{\mathbb{Z}}

\newcommand{\one}{\mathbbm{1}}

\newcommand{\dd}{\,\mathrm{d}}

\DeclareMathOperator{\conv}{conv}
\DeclareMathOperator{\Var}{Var}

\providecommand{\E}{\mathbb E}
\providecommand{\Prob}{\mathbb P}
\providecommand{\Cov}{\operatorname{Cov}}

\providecommand{\M}{\mathcal M}

\providecommand{\fall}[2]{(#1)_{#2}}

\providecommand{\bara}{\overline{\alpha}}
\providecommand{\dir}[1]{\overrightarrow{#1}}

\providecommand{\putlab}[2]{\node[lab] at #1 {#2};}

\tikzset{
  hexbg/.style={draw=gray!55,line width=.35pt},
  edge/.style={-{Latex[length=1.8mm,width=1.3mm]},line width=.75pt},
  edgeB/.style={-{Latex[length=1.8mm,width=1.3mm]},line width=.75pt},
  pot/.style={draw=black!80,line width=.55pt},
  potB/.style={draw=black!80,line width=.55pt},
  lab/.style={font=\scriptsize,inner sep=.7pt,fill=white,fill opacity=.94,text opacity=1,text=black},
  tinycase/.style={font=\scriptsize,inner sep=0pt,fill=white,fill opacity=.86,text opacity=1}
}

\newcommand{\hexpoints}{%
  \coordinate (A) at (-0.55,0.95);
  \coordinate (B) at (0.55,0.95);
  \coordinate (C) at (1.10,0.00);
  \coordinate (D) at (0.55,-0.95);
  \coordinate (E) at (-0.55,-0.95);
  \coordinate (F) at (-1.10,0.00);
}
\newcommand{\hexframe}{\draw[hexbg] (A)--(B)--(C)--(D)--(E)--(F)--cycle;}
\newcommand{\hexcoords}{\hexpoints\hexframe}

\newcommand{\smallhexbare}[1]{%
\begin{tikzpicture}[scale=.68,baseline=-0.55ex]
  \hexpoints
  #1
\end{tikzpicture}}
\newcommand{\widehexbare}[1]{%
\begin{tikzpicture}[scale=.86,baseline=-0.55ex]
  \hexpoints
  #1
\end{tikzpicture}}

\newcommand{\FigCaseIPotential}{\widehexbare{%
  \draw[edge] (A)--(B);      \node[lab] at (0,1.12) {$e$};
  \draw[edge] (C)--(D);      \node[lab] at (1.14,-0.46) {$f$};
  \draw[edge] (E)--(F);      \node[lab] at (-1.14,-0.46) {$g$};
  \draw[potB] (B)--(C);      \node[lab] at (1.14,0.46) {$x_1$};
  \draw[potB] (D)--(E);      \node[lab] at (0,-1.12) {$x_2$};
  \draw[potB] (F)--(A);      \node[lab] at (-1.14,0.46) {$x_3$};
  \draw[potB] (A)--(D);      \node[lab] at (0.27,0.46) {$x_6$};
  \draw[potB] (F)--(C);      \node[lab] at (0.55,0) {$x_5$};
  \draw[potB] (B)--(E);      \node[lab] at (0.27,-0.46) {$x_4$};
}}
\newcommand{\FigCaseIone}{\smallhexbare{%
  \draw[edge] (A)--(B);      \node[lab] at (0,1.12) {$e$};
  \draw[edge] (C)--(D);      \node[lab] at (1.14,-0.46) {$f$};
  \draw[edge] (E)--(F);      \node[lab] at (-1.14,-0.46) {$g$};
  \draw[potB] (F)--(A);      \node[lab] at (-1.14,0.46) {$x_3$};
  \draw[potB] (B)--(E);      \node[lab] at (0.27,0.46) {$x_6$};
}}
\newcommand{\FigCaseItwo}{\smallhexbare{%
  \draw[edge] (A)--(B);      \node[lab] at (0,1.12) {$e$};
  \draw[edge] (C)--(D);      \node[lab] at (1.14,-0.46) {$f$};
  \draw[edge] (E)--(F);      \node[lab] at (-1.14,-0.46) {$g$};
  \draw[potB] (B)--(C);      \node[lab] at (1.14,0.46) {$x_1$};
  \draw[potB] (D)--(E);      \node[lab] at (0,-1.12) {$x_2$};
  \draw[potB] (F)--(A);      \node[lab] at (-1.14,0.46) {$x_3$};
}}
\newcommand{\FigCaseIthree}{\smallhexbare{%
  \draw[edge] (A)--(B);      \node[lab] at (0,1.12) {$e$};
  \draw[edge] (C)--(D);      \node[lab] at (1.14,-0.46) {$f$};
  \draw[edge] (E)--(F);      \node[lab] at (-1.14,-0.46) {$g$};
  \draw[potB] (B)--(C);      \node[lab] at (1.14,0.46) {$x_1$};
  \draw[potB] (F)--(A);      \node[lab] at (-1.14,0.46) {$x_3$};
  \draw[potB] (A)--(D);      \node[lab] at (0.27,-0.46) {$x_4$};
  \draw[potB] (B)--(E);      \node[lab] at (0.27,0.46) {$x_6$};
}}
\newcommand{\FigCaseIfour}{\smallhexbare{%
  \draw[edge] (A)--(B);      \node[lab] at (0,1.12) {$e$};
  \draw[edge] (C)--(D);      \node[lab] at (1.14,-0.46) {$f$};
  \draw[edge] (E)--(F);      \node[lab] at (-1.14,-0.46) {$g$};
  \draw[potB] (B)--(C);      \node[lab] at (1.14,0.46) {$x_1$};
  \draw[potB] (D)--(E);      \node[lab] at (0,-1.12) {$x_2$};
  \draw[potB] (F)--(A);      \node[lab] at (-1.14,0.46) {$x_3$};
  \draw[potB] (F)--(C);      \node[lab] at (0.55,0) {$x_5$};
}}
\newcommand{\FigCaseIfive}{\smallhexbare{%
  \draw[edge] (A)--(B);      \node[lab] at (0,1.12) {$e$};
  \draw[edge] (C)--(D);      \node[lab] at (1.14,-0.46) {$f$};
  \draw[edge] (E)--(F);      \node[lab] at (-1.14,-0.46) {$g$};
  \draw[potB] (B)--(C);      \node[lab] at (1.14,0.46) {$x_1$};
  \draw[potB] (D)--(E);      \node[lab] at (0,-1.12) {$x_2$};
  \draw[potB] (F)--(A);      \node[lab] at (-1.14,0.46) {$x_3$};
  \draw[potB] (A)--(D);      \node[lab] at (0.27,-0.46) {$x_4$};
  \draw[potB] (B)--(E);      \node[lab] at (0.27,0.46) {$x_6$};
}}
\newcommand{\FigCaseIsix}{\smallhexbare{%
  \draw[edge] (A)--(B);      \node[lab] at (0,1.12) {$e$};
  \draw[edge] (C)--(D);      \node[lab] at (1.14,-0.46) {$f$};
  \draw[edge] (E)--(F);      \node[lab] at (-1.14,-0.46) {$g$};
  \draw[potB] (B)--(C);      \node[lab] at (1.14,0.46) {$x_1$};
  \draw[potB] (D)--(E);      \node[lab] at (0,-1.12) {$x_2$};
  \draw[potB] (F)--(A);      \node[lab] at (-1.14,0.46) {$x_3$};
  \draw[potB] (A)--(D);      \node[lab] at (0.27,-0.46) {$x_4$};
  \draw[potB] (F)--(C);      \node[lab] at (0.55,0) {$x_5$};
  \draw[potB] (B)--(E);      \node[lab] at (0.27,0.46) {$x_6$};
}}

\newcommand{\IIPentPic}[2]{%
\begin{tikzpicture}[scale=#1,baseline=-0.55ex]
  \coordinate (T) at (0,1.00);
  \coordinate (L) at (-0.95,0.31);
  \coordinate (R) at (0.95,0.31);
  \coordinate (BR) at (0.59,-0.81);
  \coordinate (BL) at (-0.59,-0.81);
  #2
\end{tikzpicture}}
\newcommand{\FigIIaBase}{\IIPentPic{0.96}{%
  \draw[edge] (L)--(T);\putlab{(-0.58,0.72)}{$e$}
  \draw[edge] (R)--(T);\putlab{(0.58,0.72)}{$f$}
  \draw[edge] (BL)--(BR);\putlab{(0,-1.00)}{$g$}
  \draw[potB] (L)--(BR);\putlab{(-0.72,-0.04)}{$x_1$}
  \draw[potB] (T)--(BL);\putlab{(-0.08,0.14)}{$x_2$}
  \draw[potB] (R)--(BR);\putlab{(0.92,-0.15)}{$x_3$}
}}
\newcommand{\FigIIaOne}{\IIPentPic{0.82}{%
  \draw[edge] (L)--(T);\putlab{(-0.58,0.72)}{$e$}
  \draw[edge] (R)--(T);\putlab{(0.58,0.72)}{$f$}
  \draw[edge] (BL)--(BR);\putlab{(0,-1.00)}{$g$}
  \draw[potB] (L)--(BR);\putlab{(-0.72,-0.04)}{$x_1$}
  \draw[potB] (T)--(BL);\putlab{(-0.08,0.14)}{$x_2$}
}}
\newcommand{\FigIIaTwo}{\IIPentPic{0.82}{%
  \draw[edge] (L)--(T);\putlab{(-0.58,0.72)}{$e$}
  \draw[edge] (R)--(T);\putlab{(0.58,0.72)}{$f$}
  \draw[edge] (BL)--(BR);\putlab{(0,-1.00)}{$g$}
  \draw[potB] (L)--(BR);\putlab{(-0.72,-0.04)}{$x_1$}
  \draw[potB] (T)--(BL);\putlab{(-0.08,0.14)}{$x_2$}
  \draw[potB] (R)--(BR);\putlab{(0.92,-0.15)}{$x_3$}
}}
\newcommand{\FigIIbBase}{%
\begin{tikzpicture}[scale=1.05,baseline=-0.55ex]
  \coordinate (T) at (0,1.05);
  \coordinate (UL) at (-1.00,0.32);
  \coordinate (UR) at (1.00,0.32);
  \coordinate (BR) at (0.62,-0.86);
  \coordinate (BL) at (-0.62,-0.86);
  \draw[edge] (UL)--(T);\putlab{(-0.60,0.74)}{$e$}
  \draw[edge] (T)--(UR);\putlab{(0.62,0.74)}{$f$}
  \draw[edge] (BL)--(BR);\putlab{(0,-1.10)}{$g$}
\end{tikzpicture}}

\newcommand{\VCasePic}[2]{%
\begin{tikzpicture}[scale=#1,baseline=-0.55ex,line cap=round]
  #2
\end{tikzpicture}}
\newcommand{\VCaseOne}{\VCasePic{.70}{%
  \draw[potB] (-1.00,0.87)--(-0.50,1.73);\putlab{(-1.10,1.28)}{$e$}
  \draw[potB] (0.50,1.73)--(1.00,0.87);\putlab{(1.10,1.28)}{$f$}
  \draw[potB] (-0.50,0)--(0.50,0);\putlab{(0.00,0.22)}{$g$}\putlab{(0.00,-0.22)}{$k$}
  \draw[potB] (-1.00,-0.87)--(-0.50,-1.73);\putlab{(-1.10,-1.28)}{$l$}
  \draw[potB] (1.00,-0.87)--(0.50,-1.73);\putlab{(1.10,-1.28)}{$m$}
}}
\newcommand{\VCaseTwo}{\VCasePic{.70}{%
  \draw[potB] (-1.00,0.87)--(-0.50,1.73);\putlab{(-1.10,1.28)}{$e$}
  \draw[potB] (0.50,1.73)--(1.00,0.87);\putlab{(1.10,1.28)}{$f$}
  \draw[potB] (-0.50,0)--(0.50,0);\putlab{(0.00,0.22)}{$g$}
  \draw[potB] (-0.50,0)--(-1.00,-0.87);\putlab{(-1.05,-0.42)}{$k$}
  \draw[potB] (-0.50,-1.73)--(0.50,-1.73);\putlab{(0.00,-1.98)}{$l$}
  \draw[potB] (0.50,0)--(1.00,-0.87);\putlab{(1.05,-0.42)}{$m$}
}}
\newcommand{\VCaseThree}{\VCasePic{.70}{
  \draw[potB] (-0.50,1.73)--(0.50,1.73);\putlab{(0.00,1.98)}{$e$}
  \draw[potB] (1.00,0.87)--(0.50,0);\putlab{(1.10,0.42)}{$f$}
  \draw[potB] (-0.50,0)--(-1.00,0.87);\putlab{(-1.05,0.42)}{$g$}
  \draw[potB] (-0.50,0)--(-1.00,-0.87);\putlab{(-1.05,-0.42)}{$k$}
  \draw[potB] (-0.50,-1.73)--(0.50,-1.73);\putlab{(0.00,-1.98)}{$l$}
  \draw[potB] (0.50,0)--(1.00,-0.87);\putlab{(1.05,-0.42)}{$m$}
}}

\newcommand{\Conf}[2]{\begin{minipage}{.19\linewidth}\centering #1\\[1.8mm]{\scriptsize (#2)}\end{minipage}}

\newcommand{\ProofFullPic}[1]{%
\begin{tikzpicture}[scale=#1,baseline=-0.55ex]
  \hexcoords
  \draw[potB] (F)--(A);\putlab{(-0.98,0.56)}{$e$}
  \draw[potB] (B)--(C);\putlab{(0.98,0.56)}{$f$}
  \draw[potB] (E)--(D);\putlab{(0,-1.18)}{$s$}
  \draw[potB] (A)--(B);\putlab{(0,1.18)}{$x_1$}
  \draw[potB] (C)--(D);\putlab{(1.08,-0.44)}{$x_2$}
  \draw[potB] (F)--(E);\putlab{(-1.08,-0.44)}{$x_3$}
  \draw[potB] (F)--(C);\putlab{(0.52,0)}{$x_4$}
  \draw[potB] (B)--(E);\putlab{(0.32,0.56)}{$x_5$}
  \draw[potB] (A)--(D);\putlab{(-0.32,0.56)}{$x_6$}
\end{tikzpicture}}
\newcommand{\ProofRootFig}{\ProofFullPic{.92}}
\newcommand{\ProofMini}[2]{\begin{minipage}{.22\linewidth}\centering #1\\[2.4mm]{\scriptsize #2}\end{minipage}}
\newcommand{\ProofSubPic}[2]{%
\begin{tikzpicture}[scale=#1,baseline=-0.55ex]
  \hexpoints
  #2
\end{tikzpicture}}
\newcommand{\ProofAlphaOne}{\ProofSubPic{.74}{%
  \draw[edge] (F)--(A);\putlab{(-1.02,0.48)}{$e$}
  \draw[edge] (B)--(C);\putlab{(1.02,0.48)}{$f$}
  \draw[edge] (D)--(E);\putlab{(0,-1.14)}{$s$}
}}
\newcommand{\ProofAlphaTwo}{\ProofSubPic{.74}{%
  \draw[edge] (B)--(A);\putlab{(0,1.16)}{$x_1$}
  \draw[edge] (F)--(C);\putlab{(0.52,0)}{$x_4$}
  \draw[edge] (D)--(E);\putlab{(0,-1.14)}{$s$}
}}
\newcommand{\ProofBetaOne}{\ProofSubPic{.74}{%
  \draw[edge] (B)--(A);\putlab{(0,1.16)}{$x_1$}
  \draw[edge] (F)--(E);\putlab{(-1.04,-0.42)}{$x_3$}
  \draw[edge,line width=.9pt] (D)--(C);\putlab{(1.08,-0.40)}{$x_2$}
}}
\newcommand{\ProofBetaTwo}{\ProofSubPic{.74}{%
  \draw[edge] (F)--(C);\putlab{(0.52,0)}{$x_4$}
  \draw[edge] (B)--(E);\putlab{(0.22,0.46)}{$x_5$}
  \draw[edge] (D)--(A);\putlab{(-0.22,0.46)}{$x_6$}
}}
\newcommand{\ProofBetaThree}{\ProofSubPic{.74}{%
  \draw[edge] (F)--(A);\putlab{(-1.02,0.48)}{$e$}
  \draw[edge] (B)--(E);\putlab{(0.22,0.46)}{$x_5$}
  \draw[edge,line width=.9pt] (D)--(C);\putlab{(1.06,-0.40)}{$x_2$}
}}
\newcommand{\ProofBetaFour}{\ProofSubPic{.74}{%
  \draw[edge] (F)--(E);\putlab{(-1.04,-0.42)}{$x_3$}
  \draw[edge] (D)--(A);\putlab{(-0.18,0.42)}{$x_6$}
  \draw[edge] (B)--(C);\putlab{(1.02,0.48)}{$f$}
}}

\newcommand{\facoords}{\hexpoints}
\newcommand{\fabg}{}
\newcommand{\faBase}{%
  \draw[edge] (F)--(A);\putlab{(-1.08,0.45)}{$e$}
  \draw[edge] (B)--(C);\putlab{(1.08,0.45)}{$f$}
  \draw[edge] (D)--(E);\putlab{(0,-1.15)}{$s$}
}
\newcommand{\faTop}{\draw[potB] (A)--(B);}
\newcommand{\faMid}{\draw[potB] (F)--(C);}
\newcommand{\faLeft}{\draw[potB] (F)--(E);}
\newcommand{\faRight}{\draw[potB] (C)--(D);}
\newcommand{\faDiagA}{\draw[potB] (A)--(D);}
\newcommand{\faDiagB}{\draw[potB] (B)--(E);}
\newcommand{\fafig}[1]{\begin{tikzpicture}[scale=.72,baseline=-0.55ex]\facoords\fabg#1\end{tikzpicture}}
\newcommand{\FAone}{\fafig{\faBase\faTop\faMid}}
\newcommand{\FAtwo}{\fafig{\faBase\faLeft\faDiagA}}
\newcommand{\FAthree}{\fafig{\faBase\faTop\faLeft\faRight}}
\newcommand{\FAfour}{\fafig{\faBase\faLeft\faRight\faDiagA\faTop}}
\newcommand{\FAfive}{\fafig{\faBase\faTop\faMid\faLeft\faRight}}
\newcommand{\FAsix}{\fafig{\faBase\faTop\faLeft\faDiagA\faMid}}
\newcommand{\FAseven}{\fafig{\faBase\faLeft\faRight\faDiagA\faDiagB}}
\newcommand{\FAeight}{\fafig{\faBase\faTop\faLeft\faRight\faDiagA\faDiagB}}
\newcommand{\FAnine}{\fafig{\faBase\faTop\faMid\faLeft\faRight\faDiagA}}
\newcommand{\FAten}{\fafig{\faBase\faTop\faMid\faLeft\faRight\faDiagA\faDiagB}}

\newcommand{\fbBase}{%
  \draw[edge] (A)--(B);\putlab{(0,1.15)}{$e$}
  \draw[edge] (C)--(D);\putlab{(1.28,-0.42)}{$f$}
  \draw[edge] (E)--(F);\putlab{(-1.28,-0.42)}{$g$}
  \draw[line width=.75pt] (D)--(E);\putlab{(0,-1.18)}{$s$}
}
\newcommand{\fbUL}{\draw[potB] (F)--(A);}
\newcommand{\fbUR}{\draw[potB] (B)--(C);}
\newcommand{\fbMid}{\draw[potB] (F)--(C);}
\newcommand{\fbAD}{\draw[potB] (A)--(D);}
\newcommand{\fbBE}{\draw[potB] (B)--(E);}
\newcommand{\fbfig}[1]{\begin{tikzpicture}[scale=.68,baseline=-0.55ex]\hexpoints #1\end{tikzpicture}}
\newcommand{\FBone}{\fbfig{\fbBase\fbAD\fbUR}}
\newcommand{\FBtwo}{\fbfig{\fbBase\fbMid}}
\newcommand{\FBthree}{\fbfig{\fbBase\fbUL\fbUR}}
\newcommand{\FBfour}{\fbfig{\fbBase\fbUL\fbUR\fbBE}}
\newcommand{\FBfive}{\fbfig{\fbBase\fbUL\fbUR\fbMid}}
\newcommand{\FBsix}{\fbfig{\fbBase\fbUL\fbMid\fbBE}}
\newcommand{\FBseven}{\fbfig{\fbBase\fbUL\fbUR\fbAD\fbBE}}
\newcommand{\FBeight}{\fbfig{\fbBase\fbUL\fbUR\fbMid\fbBE}}
\newcommand{\FBnine}{\fbfig{%
  \fbBase\fbUL\fbUR\fbMid
  \draw[potB,line width=.9pt] (A)--(D);
  \draw[potB,line width=.9pt] (B)--(E);
}}

\newcommand{\fbZeroBase}{%
  \draw[edge] (A)--(B);\putlab{(0,1.15)}{$e$}
  \draw[edge] (C)--(D);\putlab{(1.28,-0.42)}{$f$}
  \draw[edge] (E)--(F);\putlab{(-1.28,-0.42)}{$g$}
}
\newcommand{\fbzfig}[1]{\begin{tikzpicture}[scale=.68,baseline=-0.55ex]\hexpoints #1\end{tikzpicture}}
\newcommand{\FBZone}{\fbzfig{\fbZeroBase\draw[potB] (F)--(A);\putlab{(-1.03,0.54)}{$x_1$}\draw[potB] (B)--(E);\putlab{(0.24,0.50)}{$x_4$}}}
\newcommand{\FBZtwo}{\fbzfig{\fbZeroBase\draw[potB] (F)--(C);\putlab{(0.48,0)}{$x_3$}\draw[potB] (A)--(D);
\putlab{(-0.24,0.50)}{$x_5$}\draw[potB] (B)--(E);\putlab{(0.24,0.50)}{$x_4$}}}
\newcommand{\FBZthree}{\fbzfig{\fbZeroBase\draw[potB] (F)--(A);\putlab{(-1.03,0.54)}{$x_1$}\draw[potB] (B)--(C);
\putlab{(1.00,0.56)}{$x_2$}\draw[potB] (A)--(D);\putlab{(-0.24,0.50)}{$x_5$}\draw[potB] (B)--(E);\putlab{(0.24,0.50)}{$x_4$}}}
\newcommand{\FBZfour}{\fbzfig{\fbZeroBase\draw[potB] (F)--(A);\putlab{(-1.03,0.54)}{$x_1$}\draw[potB] (F)--(C);
\putlab{(0.48,0)}{$x_3$}\draw[potB] (A)--(D);\putlab{(-0.24,0.50)}{$x_5$}\draw[potB] (B)--(E);\putlab{(0.24,0.50)}{$x_4$}}}
\newcommand{\FBZfive}{\fbzfig{\fbZeroBase\draw[potB] (F)--(A);\putlab{(-1.03,0.54)}{$x_1$}\draw[potB] (B)--(C);
\putlab{(1.00,0.56)}{$x_2$}\draw[potB] (F)--(C);\putlab{(0.48,0)}{$x_3$}\draw[potB] (A)--(D);\putlab{(-0.24,0.50)}{$x_5$}\draw[potB] (B)--(E);\putlab{(0.24,0.50)}{$x_4$}}}

\newcommand{\RefPic}[1]{%
\begin{tikzpicture}[
  scale=.76,
  baseline=-0.55ex,
  line cap=round,
  line join=round
]
#1
\end{tikzpicture}}

\newcommand{\RefConf}[2]{%
\begin{minipage}[c]{.29\linewidth}
  \centering
  #1\par\vspace{1.4mm}
  {\scriptsize #2}
\end{minipage}%
}

\newcommand{\RefHexPoints}{%
  \coordinate (T)  at (0,1.10);
  \coordinate (UR) at (.953,.55);
  \coordinate (LR) at (.953,-.55);
  \coordinate (B)  at (0,-1.10);
  \coordinate (LL) at (-.953,-.55);
  \coordinate (UL) at (-.953,.55);
}

\newcommand{\RefHexInterior}{%
  \coordinate (P) at (-.476,.275);
  \coordinate (Q) at (.476,-.275);

  \draw[potB] (T)--(P);
  \node[lab] at (-.31,.68) {$x_1$};

  \draw[potB] (P)--(Q);
  \node[lab] at (0,.15) {$x_2$};

  \draw[potB] (Q)--(B);
  \node[lab] at (.31,-.68) {$x_3$};
}

\newcommand{\RefAdjOne}{\RefPic{%
  \RefHexPoints

  \draw[potB]
    (T)--(UR)--(LR)--(B)--(LL)--(UL)--cycle;

  \node[lab] at (-.61,.88) {$e$};
  \node[lab] at (.61,.88) {$f$};

  \RefHexInterior
}}

\newcommand{\RefAdjTwo}{\RefPic{%
  \coordinate (A) at (-.82,.82);
  \coordinate (B) at (.82,.82);
  \coordinate (C) at (.82,-.82);
  \coordinate (D) at (-.82,-.82);
  \coordinate (M) at (0,.82);

  \draw[potB] (A)--(B)--(C)--(D)--cycle;

  \node[lab] at (-1.00,0) {$x_1$};
  \node[lab] at (0,-1.00) {$x_2$};
  \node[lab] at (1.00,0) {$x_3$};

  \draw[potB] (D)--(M);
  \node[lab] at (-.37,.08) {$e$};

  \draw[potB] (M)--(C);
  \node[lab] at (.37,.08) {$f$};
}}
\newcommand{\RefSquarePoints}{%
  \coordinate (A) at (-.82,.82);
  \coordinate (B) at (.82,.82);
  \coordinate (C) at (.82,-.82);
  \coordinate (D) at (-.82,-.82);

  \coordinate (P) at (.2197,-.2197);
}

\newcommand{\RefSquareInterior}{%
  \draw[potB] (D)--(P);
  \node[lab] at (-.40,-.57) {$x_1$};

  \draw[potB] (P)--(B);
  \node[lab] at (.56,.34) {$x_2$};

  \draw[potB] (1.18,0)--(1.98,0);
  \node[lab] at (1.58,-.21) {$x_3$};
}

\newcommand{\RefAdjThree}{\RefPic{%
  \RefSquarePoints

  \draw[potB] (A)--(B)--(C)--(D)--cycle;

  \node[lab] at (-1.00,0) {$e$};
  \node[lab] at (0,-1.00) {$f$};

  \RefSquareInterior
}}

\newcommand{\RefDisOne}{\RefPic{%
  \RefHexPoints

  \draw[potB]
    (T)--(UR)--(LR)--(B)--(LL)--(UL)--cycle;

  \node[lab] at (-.61,.88) {$e$};
  \node[lab] at (.61,-.88) {$f$};

  \RefHexInterior
}}

\newcommand{\RefDisTwo}{\RefPic{%
  \RefHexPoints

  \draw[potB]
    (T)--(UR)--(LR)--(B)--(LL)--(UL)--cycle;

  \node[lab] at (-.61,.88) {$e$};
  \node[lab] at (1.15,0) {$f$};

  \RefHexInterior
}}

\newcommand{\RefDisThree}{\RefPic{%
  \RefSquarePoints

  \draw[potB] (A)--(B)--(C)--(D)--cycle;

  \node[lab] at (-1.00,0) {$e$};
  \node[lab] at (1.00,0) {$f$};

  \RefSquareInterior
}}

\title{Berry--Esseen Bounds for Faces in Unimodular Triangulations of Random Polytopes}
\author{Zhonggen Su \qquad Huilian Zheng\\[0.6em]
\small School of Mathematical Sciences, Zhejiang University\\
\small \texttt{suzhonggen@zju.edu.cn}; \texttt{12535018@zju.edu.cn}}
\date{}

\begin{document}

\maketitle

\begin{abstract}
Consider the number of faces arising from unimodular triangulations of random
polytopes associated with Erd\H{o}s--R\'{e}nyi graphs.
Focusing on two-dimensional faces of regular unimodular triangulations of random
symmetric edge polytopes, we derive precise asymptotic expansions for the
expectation and variance and then establish a quantitative central limit theorem
with an explicit Berry--Esseen bound.
The proof combines an algebraically derived characterization of faces in terms
of local graph configurations with a detailed combinatorial--geometric study of
the relevant configurations.
The normal approximation is established using the discrete Malliavin--Stein
method.
We also determine the asymptotic orders of the expectation and variance and
establish a quantitative central limit theorem for the number of faces in each
fixed higher dimension.
These results identify a common local mechanism underlying the fluctuations of
face numbers and extend the distributional theory of random symmetric edge
polytopes \href{https://arxiv.org/abs/2603.10650}{[Donzelmann et al.,
arXiv:2603.10650]} beyond edge statistics.
\end{abstract}

\noindent\textbf{2020 Mathematics Subject Classification.} Primary 60F05; secondary 52B20, 05C80, 60D05.\par
\medskip
\noindent\textbf{Keywords.} Berry--Esseen bound; Discrete Malliavin--Stein method; Erd\H{o}s--R\'{e}nyi random graph; Symmetric edge polytope; Unimodular triangulation.
\medskip

\section{Introduction and main results}\label{sec:intro}

\subsection{Background and motivation}

Random polytopes occupy a central place in stochastic geometry.
In the classical setting, they are obtained as convex hulls of points sampled
from continuous distributions in fixed-dimensional ambient spaces.
Following the pioneering work of R\'{e}nyi and Sulanke
\cite{RenyiSulanke1963,RenyiSulanke1964}, a substantial literature has developed
on the asymptotic behavior of geometric and combinatorial functionals of random
polytopes, such as intrinsic volumes and face numbers.
Central limit theorems are now available for several fundamental models,
including random polytopes in smooth convex bodies, Gaussian polytopes, and
Poisson polytopes
\cite{Vu2006,BaranyVu2007,BaranyReitzner2010,ThaeleTurchiWespi2018}.
In many of these results, the asymptotic behavior is governed by stabilization
phenomena near the boundary of the underlying convex body; see Lachi\`eze-Rey et
al. \cite{LachiezeReySchulteYukich2019} for a general normal-approximation
framework.

Random lattice polytopes arising from discrete structures form a distinct class
of models.
Rather than sampling points from a continuous distribution, one starts with a
random graph, complex, or matroid and encodes its combinatorial data in a
lattice point configuration.
In many natural examples, the ambient dimension grows with the size of the
underlying discrete structure, and the geometry is governed by combinatorial
relations rather than by the boundary of a fixed convex body.
Distributional questions in this setting therefore call for a combination of
probabilistic enumeration, discrete geometry, and combinatorial commutative
algebra.

Symmetric edge polytopes provide a natural and particularly tractable model of
this type.
To avoid confusion, we use the terms \emph{node} and \emph{arc} for the vertices
and edges of the underlying graph, and reserve \emph{vertex} and \emph{edge} for
the corresponding geometric objects.
Let $G=(V,E)$ be a finite simple graph, and let $(e_i)_{i\in V}$ denote the
standard basis of $\mathbb R^V$.
The symmetric edge polytope associated with $G$ is defined as
\begin{equation}\label{eq:symmetric-edge-polytope}
  P_G:=\conv\{\pm(e_i-e_j):\{i,j\}\in E\}\subseteq\mathbb R^V.
\end{equation}
Each arc $\{i,j\}\in E$ gives rise to the antipodal vertices $e_i-e_j$ and
$e_j-e_i$, which we identify with the directed arcs $i\to j$ and $j\to i$,
respectively.
This identification links the facial structure of $P_G$ to directed
configurations in $G$.
In particular, short directed cycles give rise to local obstructions to faces of
the polytope and of its regular unimodular triangulations.

This connection has led to extensive work on toric ideals, Gr\"obner bases,
Ehrhart theory, and combinatorial invariants such as $f$-vectors and
$\gamma$-vectors; see
\cite{OhsugiHibi1999,OhsugiHibi2011,HigashitaniJochemkoMichalek2019,
DaliDelucchiMichalek2022,DaliJuhnkeKubitzkeKohneVenturello2023,
CodenottiRiccardiVenturello2025}.
Symmetric edge polytopes therefore provide a natural setting in which graph
theory, lattice-polytope geometry, and combinatorial commutative algebra
interact.
This correspondence is particularly useful in the random setting, where the
relevant directed configurations become local events in the underlying random
graph.

To obtain a random polytope, let $G_{n,p}$ be the Erd\H{o}s--R\'{e}nyi random
graph on node set $[n]$, in which every arc is present independently with
probability $p=p(n)$, and set
\[
  P_{n,p}:=P_{G_{n,p}}.
\]
Donzelmann et al. \cite{DonzelmannJuhnkeRednossThaele2026} initiated the
distributional study of this model, and they obtained sharp expectation and
variance asymptotics and quantitative central limit theorems for the number of
edges of both $P_{n,p}$ and its unimodular triangulations.
Their results raise a natural and challenging question: is there a common
fluctuation mechanism underlying higher-dimensional face statistics of random
symmetric edge polytopes?
The passage from edges to higher-dimensional faces is not merely a routine
extension.
Indeed, for two-dimensional faces, both the relevant local events and their
dependence structure are already substantially more intricate.

To make this precise, recall that all unimodular triangulations of a fixed
lattice polytope have the same $f$-vector.
Hence the face numbers considered below do not depend on the particular
unimodular triangulation.
We may therefore fix a regular unimodular triangulation $\Delta_{n,p}$ of
$P_{n,p}$, on the vertices of $P_{n,p}$ and the origin.
The triangulation is chosen so that its facial structure can be described in
terms of local obstructions arising from directed cycles in $G_{n,p}$.
For $k\geq 0$, faces of dimension $k$ are called $k$-faces, and we write
$f_k(\Delta_{n,p})$ for the number of $k$-faces of $\Delta_{n,p}$.
A $k$-face is called a non-origin $k$-face if it does not contain the origin.
Such a face is generated by $k+1$ directed arcs, and whether these arcs generate
a face can be characterized by local conditions arising from directed cycles.

In this paper, we shall first derive refined asymptotic formulas with explicit
leading coefficients for both the expectation and variance of
$f_2(\Delta_{n,p})$.
To this end, we need to develop two novel ideas.
One is to equip each arc with an auxiliary continuous random weight to encode
the relative order of arcs, so that we reduce the dominant covariance
contribution to a one-dimensional integral by conditioning on a common arc.
The other is the surprising discovery that the resulting integrand, which
appears to depend on the integration variable, is in fact constant, thanks to an
underlying combinatorial invariance with respect to the ordering of the arcs,
allowing us to avoid extensive combinatorial enumeration.
We shall also extend the preceding analysis to the case of $k$-faces where
$k\ge3$, but it becomes increasingly cumbersome to evaluate the leading
coefficients explicitly.
Hence, we are content with providing the asymptotic orders of the expectation
and variance of $f_k(\Delta_{n,p})$.
Furthermore, we shall use the well-known discrete Malliavin--Stein method to
establish a quantitative central limit theorem with an explicit convergence rate
for $f_k(\Delta_{n,p})$ for every fixed $k\ge2$.

\subsection{Main results}

Throughout the paper, $p=p(n)\in(0,1)$, and we set $q:=1-p$.
For $r\in\mathbb{N}$, let
\[
  (n)_r:=n(n-1)\cdots(n-r+1)
\]
denote the falling factorial.
For positive sequences $a_n$ and $b_n$, we write $a_n\lesssim b_n$ if
$a_n=O(b_n)$, $a_n\ll b_n$ if $a_n=o(b_n)$, and $a_n\asymp b_n$ if both
$a_n\lesssim b_n$ and $b_n\lesssim a_n$.
The implicit constants may depend on fixed parameters, but not on $n$.

We shall use the following polynomials:
\[
  \begin{aligned}
    H_1(p) &:=1-\frac{3}{2}p^2-p^3+3p^4-\frac{3}{2}p^5+\frac{1}{6}p^6,\\
    H_2(p) &:=1-p^2+\frac{1}{3}p^3,\\
    H_3(p) &:=p^4-6p^3+p^2+4p+2.
  \end{aligned}
\]

We begin with the precise result for $2$-faces.
For convenience, write $F_{n,p}:=f_2(\Delta_{n,p}).$

\begin{theorem}\label{thm:main}
Assume that $p\gg n^{-2}$ and $\limsup\limits_{n\to\infty}p<1.$ Then, as
$n\to\infty$, the following assertions hold.

\begin{enumerate}
\item[(i)]
The expectation and variance satisfy
\[
  \mathbb{E}F_{n,p}
  = \frac{1}{6}(n)_6p^3H_1(p) + (n)_5p^3H_2(p) + o\bigl(n^5p^{5/2}\bigr),
\]
and
\[
  \operatorname{Var}(F_{n,p})
  = \frac{1}{8}(n)_{10}p^5q^5H_3(p)^2 + o\bigl(n^{10}p^5\bigr).
\]

\item[(ii)]
If $N$ is a standard Gaussian random variable, then
\[
  d_K\left( \frac{F_{n,p}-\mathbb{E}F_{n,p}} {\sqrt{\operatorname{Var}(F_{n,p})}}, N
  \right)
  \lesssim \frac{1}{n\sqrt{p}}.
\]
In particular,
\[
  \frac{F_{n,p}-\mathbb{E}F_{n,p}} {\sqrt{\operatorname{Var}(F_{n,p})}} \xrightarrow{d}
  N(0,1).
\]
\end{enumerate}
\end{theorem}

We further extend the preceding results to the number of faces in every fixed
dimension $k\ge 3$.
Write
\[
  F^{[k]}_{n,p}:=f_k(\Delta_{n,p}).
\]
In the following statement, the subscript $k$ indicates that the implicit
constants may depend on $k$.

\begin{theorem}\label{thm:k-face-clt}
Let $k\ge3$ be fixed.
Under the assumptions of Theorem~\ref{thm:main}, as $n\to\infty$, the following
assertions hold.

\begin{enumerate}
\item[(i)]
The expectation and variance satisfy
\[
  \mathbb E F^{[k]}_{n,p}
  \asymp_k n^{2k+2}p^{k+1},
  \qquad \operatorname{Var}(F^{[k]}_{n,p})
  \asymp_k n^{4k+2}p^{2k+1}.
\]

\item[(ii)]
If $N$ is a standard Gaussian random variable, then
\[
  d_K\left( \frac{F^{[k]}_{n,p}-\E F^{[k]}_{n,p}} {\sqrt{\Var(F^{[k]}_{n,p})}},N\right)
  \lesssim_k\frac{1}{n\sqrt p\,q^{2k+1/2}}.
\]
In particular,
\[
  \frac{F^{[k]}_{n,p}-\E F^{[k]}_{n,p}} {\sqrt{\Var(F^{[k]}_{n,p})}} \xrightarrow{d}
  N(0,1).
\]
\end{enumerate}
\end{theorem}

\begin{remark}[Variance and parameter regime]
Since $H_3(p)=2+p(1-p)(4+5p-p^2)>2$ for $0<p<1$, and $q$ is bounded away from
zero, the variance asymptotic in Theorem~\ref{thm:main} is nondegenerate.
In particular, $\operatorname{Var}(F_{n,p})\asymp n^{10}p^5$.
The condition $p\gg n^{-2}$ is equivalent to $n\sqrt p\to\infty$.
Thus, under the stated assumptions, the bounds in Theorems~\ref{thm:main} and~\ref{thm:k-face-clt}
tend to zero throughout this regime.
\end{remark}

The remainder of the paper is organized as follows.
Section~\ref{sec:prelim} briefly reviews the geometric--algebraic construction
of the triangulation, describes its non-origin $2$-faces in terms of directed
arcs, and states the discrete Malliavin--Stein bound used for normal
approximation.
Section~\ref{sec:expectation-variance-gradients} establishes the expectation and
variance asymptotics for $F_{n,p}$ and derives the required estimates for its
first- and second-order discrete gradients.
Section~\ref{sec:clt} combines these ingredients to prove Theorem~\ref{thm:main}.
Finally, Section~\ref{sec:higher-dimensional} proves Theorem~\ref{thm:k-face-clt} for
higher-dimensional faces.
In contrast to the two-dimensional case, we determine the asymptotic orders of
the expectation and variance rather than refined asymptotic formulas.
The normal approximation is then obtained within the same Malliavin--Stein
framework.

\section{Preliminaries}\label{sec:prelim}

This section presents the geometric, algebraic, and probabilistic framework for
the analysis that follows.
We first review lattice polytopes and unimodular triangulations.
We then describe the order-induced regular unimodular triangulation of a
symmetric edge polytope and determine when three directed arcs generate a
$2$-face.
Finally, we state the discrete Malliavin--Stein bound for normal approximation.

\subsection{Lattice polytopes and unimodular triangulations}
\label{subsec:lattice-polytopes}

In this subsection, we briefly recall standard notions on lattice polytopes and
triangulations.
For more details, we refer the reader to Beck and Robins \cite{BeckRobins2015}
and Ziegler \cite{Ziegler1995}.

A lattice polytope $P\subset\RR^d$ is the convex hull of finitely many points in
$\ZZ^d$, that is,
\[
  P=\conv(x_1,\ldots,x_r), \qquad\text{for }x_1,\ldots,x_r\in\ZZ^d.
\]
We write $\dim P$ for the dimension of its affine hull.

A supporting hyperplane of $P$ is a hyperplane such that $P$ is contained in
exactly one of its two closed halfspaces.
A face of $P$ is the intersection of $P$ with a supporting hyperplane.
Throughout this paper, we do not regard $P$ itself as a face.
The dimension of a face is defined as the dimension of its affine hull.
Faces of dimensions zero and one are called vertices and edges, respectively.
Faces of dimension $k$ are called $k$-faces, and $f_k$ denotes the number of
$k$-dimensional faces of a polytope or simplicial complex.

A triangulation $\mathcal{T}$ of a lattice polytope $P$ is a simplicial complex
whose underlying space is $P$ and whose vertices are lattice points of $P$.
It is unimodular if every maximal simplex of $\mathcal{T}$ is unimodular.
Since a unimodular simplex contains no lattice points other than its vertices,
every lattice point of $P$ is a vertex of any unimodular triangulation of $P$.
We write $f_k(\mathcal{T})$ for the number of $k$-dimensional faces of
$\mathcal{T}$.

A triangulation is regular if it is induced by a height function on its
vertices: the vertices are lifted to one higher dimension, and the lower faces
of their convex hull are projected back to $P$.
Regular triangulations can be described by initial ideals of toric ideals.
In particular, a squarefree initial ideal determines a regular unimodular
triangulation; see Sturmfels~\cite{Sturmfels1996}.

Although the $f$-vector of a triangulation generally depends on the
triangulation, it is independent of the choice of unimodular triangulation.
Indeed, the $h$-polynomial of a unimodular triangulation agrees with the Ehrhart
$h^*$-polynomial of $P$, and the $h$-vector determines the $f$-vector; see
Stanley~\cite{Stanley1980}, Sturmfels~\cite{Sturmfels1996} and Beck and
Robins~\cite{BeckRobins2015}.
Consequently, $f_k(\mathcal{T})$ depends only on $P$ whenever $\mathcal{T}$ is
unimodular.
We nevertheless work with a particular order-induced regular unimodular
triangulation, whose faces can be described in terms of directed arcs and
cycles.

\subsection{Symmetric edge polytopes and their triangulations}
\label{subsec:symmetric-edge-polytopes}

Let $G=(V,E)$ be a finite simple graph, and let $P_G$ be the symmetric edge
polytope defined in (\ref{eq:symmetric-edge-polytope}).
We retain the identification of the vertex $e_i-e_j$ with the directed arc
$i\to j$.

We now describe the order-induced regular unimodular triangulation of $P_G$.
Higashitani et al. \cite[Proposition~3.8]{HigashitaniJochemkoMichalek2019}
described a Gr\"obner basis for the toric ideal $I_{\mathcal A_G}$ associated
with the lattice configuration
\[
  \mathcal A_G:=\{0\}\cup\{\pm(e_i-e_j):\{i,j\}\in E\},
\]
obtained by adjoining the origin to the vertex configuration of $P_G$.
For each arc $e\in E$, let $x_e$ and $y_e$ denote the variables associated with
its two orientations, and let $z$ denote the variable associated with the
origin.
Fix a total order $e_1\prec\cdots\prec e_m$ on the underlying arcs of $G$.
We use the degree reverse lexicographic order induced by the paired variable
order
\[
  z<x_{e_1}<y_{e_1}<\cdots<x_{e_m}<y_{e_m}.
\]

The resulting Gr\"obner basis consists of binomials arising from directed
cycles, together with $x_ey_e-z^2$ for $e\in E$.
Its initial ideal is squarefree, so the associated initial complex is a regular
unimodular triangulation of $P_G$ on the vertices of $P_G$ together with the
origin.

The squarefree monomials in the initial ideal encode the nonfaces of this
triangulation.
Restricting to degrees two and three gives the following local obstructions:
\begin{enumerate}[label=\textup{(\roman*)}]
    \item the monomial $x_ey_e$ forbids a pair of antipodal vertices;

    \item a directed $3$-cycle gives a forbidden pair of directed arcs;

    \item a directed $4$-cycle gives a forbidden pair unless the pair contains
    the $\prec$-minimal arc of the cycle;

    \item a directed $5$-cycle gives a forbidden triple of directed arcs;

    \item a directed $6$-cycle gives a forbidden triple unless the triple
    contains the $\prec$-minimal arc of the cycle.
\end{enumerate}

Consequently, a triple of directed arcs generates a $2$-face if and only if the
product of the variables associated with them is not divisible by any of the
initial monomials listed above.
The following proposition expresses these conditions in terms of directed arcs
and cycles.

\begin{proposition}\label{prop:two-face-criterion}
Let $G$ be a graph equipped with a total order $\prec$ on its arcs, and let
$\mathcal{T}_{\prec}$ be the corresponding order-induced regular unimodular
triangulation of $P_G$, with the origin adjoined as a vertex.
For a set $\alpha=\{\dir e,\dir f,\dir g\}$ of three directed arcs, the
corresponding vertices span a $2$-face of $\mathcal T_\prec$ if and only if all
of the following conditions hold:
\begin{enumerate}[label=(\roman*)]
\item No two arcs of $\alpha$ have the same underlying undirected arc;

\item No two arcs in $\alpha$ lie on a directed $3$-cycle;

\item For every pair $\beta\subseteq\alpha$ whose arcs lie on a directed
$4$-cycle, the $\prec$-minimal arc of this cycle belongs to $\beta$;

\item The three arcs in $\alpha$ do not lie on a directed $5$-cycle;

\item If the three arcs in $\alpha$ lie on a directed $6$-cycle, then the
$\prec$-minimal arc of this cycle belongs to $\alpha$.
\end{enumerate}
\end{proposition}

Applied to $G_{n,p}$, Proposition~\ref{prop:two-face-criterion} gives an indicator
representation for the non-origin $2$-faces, while the faces containing the
origin can be recovered from lower-dimensional face counts.
This representation provides the geometric--algebraic input for the expectation
and variance calculations.
We next recall the probabilistic tool used to control its fluctuations.

\begingroup \emergencystretch=1em
\subsection{The discrete Malliavin--\allowbreak Stein method}\label{subsec:malliavin-stein}

The Malliavin--Stein method combines Stein's method \cite{Stein1972,Stein1986}
with Malliavin calculus to obtain quantitative normal-approximation bounds; see
Nualart and Peccati \cite{NualartPeccati2005} and Nourdin and Peccati
\cite{NourdinPeccati2009,NourdinPeccati2012}.
For functionals of independent Rademacher variables, a discrete version was
developed by Nourdin et al. \cite{NourdinPeccatiReinert2010} and further refined
by Krokowski et al.
\cite{KrokowskiReichenbachsThaele2016,KrokowskiReichenbachsThaele2017}.
We use the simplified second-order inequality of Eichelsbacher et al.
\cite{EichelsbacherRednossThaeleZheng2023}.
It is particularly effective for local graph functionals because it requires
only moments of first- and second-order discrete gradients.

\par \endgroup

Let $m\ge 1$, and let $X_1,\ldots,X_m$ be independent Rademacher random
variables such that
\[
  \mathbb P(X_i=+1)=p_i, \qquad \mathbb P(X_i=-1)=q_i:=1-p_i .
\]
Let $F=F(X_1,\ldots,X_m)$ be a real-valued measurable functional.
For $i\in[m]$, define
\[
  F_i^+ := F(X_1,\ldots,X_{i-1},+1,X_{i+1},\ldots,X_m),
\]
and
\[
  F_i^- := F(X_1,\ldots,X_{i-1},-1,X_{i+1},\ldots,X_m).
\]
The first-order discrete gradient is
\[
  D_iF:=\sqrt{p_iq_i}\,(F_i^+-F_i^-).
\]
Thus $D_iF$ measures, up to the natural normalization, how much $F$ changes when
the $i$th coordinate is flipped.

For $i,j\in[m]$, define
\[
  F_{j,i}^{+,+}:=(F_j^+)_i^+, \qquad F_{j,i}^{-,+}:=(F_j^-)_i^+,
\]
\[
  F_{j,i}^{+,-}:=(F_j^+)_i^-, \qquad F_{j,i}^{-,-}:=(F_j^-)_i^-.
\]
The second-order discrete gradient is
\[
  D_iD_jF := \sqrt{p_iq_ip_jq_j} \left( F_{j,i}^{+,+} - F_{j,i}^{-,+} - F_{j,i}^{+,-} +
  F_{j,i}^{-,-} \right),
\]
with the convention $D_iD_iF=0$.

For real-valued random variables $X$ and $Y$, the Kolmogorov distance is
\[
  d_K(X,Y) := \sup_{z\in\mathbb R} \left| \mathbb P(X\le z)-\mathbb P(Y\le z) \right|.
\]

The following finite-dimensional second-order bound is the normal-approximation
tool used in the proof of Theorem~\ref{thm:main}.

\begin{proposition}\label{prop:MS}
Let $F=F(X_1,\ldots,X_m)$ satisfy $\E F=0$ and $\E F^2=1$.
Define
\begin{align*}
  B_1:=&\sum_{j,k,\ell=1}^m\sqrt{\E[(D_jF)^2(D_kF)^2]}\,\sqrt{\E[(D_\ell D_jF)^2(D_\ell
  D_kF)^2]},\\
  B_2:=&\sum_{j,k,\ell=1}^m\frac1{p_\ell q_\ell}\E[(D_\ell D_jF)^2(D_\ell D_kF)^2],\\
  B_3:=&\sum_{k=1}^m\frac1{p_kq_k}\E[(D_kF)^4],\\
  B_4:=&\sum_{k,\ell=1}^m\frac1{p_kq_k}\sqrt{\E[(D_kF)^4]}\,\sqrt{\E[(D_\ell
  D_kF)^4]},\\
  B_5:=&\sum_{k,\ell=1}^m\frac1{p_kq_kp_\ell q_\ell}\E[(D_\ell D_kF)^4].
\end{align*}
Then, for a standard Gaussian random variable $N$,
\[
  d_K(F,N)
  \le \frac{\sqrt{15}}2\sqrt{B_1} + \frac{\sqrt3}2\sqrt{B_2} + 2\sqrt{B_3} +
  2\sqrt6\sqrt{B_4} + 2\sqrt3\sqrt{B_5}.
\]
\end{proposition}

\begin{remark}[Application to the face count]\label{rem:MS-history}
Proposition~\ref{prop:MS} is the finite-dimensional version of the simplified second-order
Gaussian Poincar\'e inequality of Eichelsbacher et al.
\cite[Theorem~4.1]{EichelsbacherRednossThaeleZheng2023}.
Compared with the inequality of Krokowski et al.
\cite[Theorem~4.1]{KrokowskiReichenbachsThaele2017}, it does not contain the
terms corresponding to $A_3$ and $A_4$ therein.
All remaining terms are expressed through the moments of the first- and
second-order discrete gradients.

Note that for the local graph functional $F_{n,p}$, changing one or two arcs
affects only the directed-cycle configurations containing those arcs.
Hence, the first- and second-order discrete gradients can be bounded by
analyzing local enumeration.
No explicit chaotic decomposition or multiplication formula for discrete
multiple stochastic integrals is needed.
\end{remark}

\section{Expectation, Variance and Discrete Gradients of \texorpdfstring{$F_{n,p}$}{F(n,p)}}
\label{sec:expectation-variance-gradients}
In this section, we first derive asymptotics for the expectation and variance of
$F_{n,p}$, and then obtain suitable bounds for its first- and second-order
discrete gradients, together with estimates on their fourth moments.
\subsection{Expectation}\label{sec:expectation}
We begin with the expectation of $F_{n,p}$.
A full exact computation is rather involved, so we content ourselves with
computing the leading term explicitly and bound the remainder to the required
order. Although individual face indicators may depend on the
chosen total order, the total face count does not.
Thus, when computing expectations, we may average
uniformly over all total orders on $\binom{[n]}{2}$.
Under this averaging, the relative order of any fixed
collection of distinct arcs is uniform.
\begin{theorem}\label{thm:expectation}
If $p\gg n^{-2}$ and $\limsup\limits_{n\to\infty}p<1$, then as $n\to\infty$,
\begin{align}\label{eq:EFmain-full-en}
  \E F_{n,p}=\frac16\fall n6 p^3H_1(p)+\fall n5p^3H_2(p)+o\bigl(n^5p^{5/2}\bigr),
\end{align}
where
\begin{align}
  H_1(p)&=1-\frac32p^2-p^3+3p^4-\frac32p^5+\frac16p^6,\label{eq:H1-polynomial}\\
  H_2(p)&=1-p^2+\frac13p^3.\label{eq:H2-polynomial}
\end{align}
\end{theorem}

Let $E_{n,p}$ be the number of arcs of $G_{n,p}$, $K_{n,p}$ the number of edges
of the unimodular triangulation $\Delta_{n,p}$ of $P_{n,p}$, and $H_{n,p}$ the
number of non-origin $2$-faces of $\Delta_{n,p}$.
Denote by $\mathcal{M}$ the set of all triples of possible directed arcs on
$[n]$.
It follows that $H_{n,p}=\sum_{\alpha\in\mathcal{M}} I_{\alpha},$ where
$I_{\alpha}$ is the indicator function of the event that the triple $\alpha$
generates a $2$-face in the unimodular triangulation $\Delta_{n,p}$ of
$P_{n,p}$.
The following lemma records the relationship among these quantities.
\begin{lemma}\label{lem:two-face-decomposition}
  $F_{n,p}=H_{n,p}+K_{n,p}-2E_{n,p}.$
\end{lemma}

\begin{proof}
We divide the 2-faces into two classes according to whether they contain the
origin.
The non-origin 2-faces are counted by $H_{n,p}$, while the others are formed by
the origin and an edge of $\Delta_{n,p}$ not incident to the origin, so the
number of them is $K_{n,p}-2E_{n,p}$.
Adding the contributions of the two classes gives the desired identity.
\end{proof}
Note that \cite[Theorem~4.1]{DonzelmannJuhnkeRednossThaele2026} gives
\begin{equation}\label{eq:KminusE-full-en}
  \begin{aligned}
    \E K_{n,p}-2\E E_{n,p} &=6\binom n4p^2(2-p^2)+6\binom n3p^2 \\
    &\quad+2\binom n2(1-p)\bigl(1-(1-p^2)^{n-2}\bigr) =o\bigl(n^5p^{5/2}\bigr).
  \end{aligned}
\end{equation}
Therefore, by Lemma~\ref{lem:two-face-decomposition}, it remains to compute
$\E H_{n,p}$.
The following lemma gives the required asymptotic.
\begin{lemma}\label{lem:H-expectation}
Under the assumptions of Theorem~\ref{thm:expectation}, as $n\to\infty$, we have
\begin{align}
  \E H_{n,p}=\frac16\fall n6 p^3H_1(p)+\fall n5p^3H_2(p)+o\bigl(n^5p^{5/2}\bigr),
\end{align}
where
\begin{align}
  H_1(p)&=1-\frac32p^2-p^3+3p^4-\frac32p^5+\frac16p^6,\\
  H_2(p)&=1-p^2+\frac13p^3.
\end{align}
\end{lemma}
\begin{proof}
In order to compute the expectation of $H_{n,p}$, we distinguish the following
cases according to the overlap pattern of the underlying undirected arcs
$\{e,f,g\}$, illustrated in Figure~\ref{fig:overlap-types}:

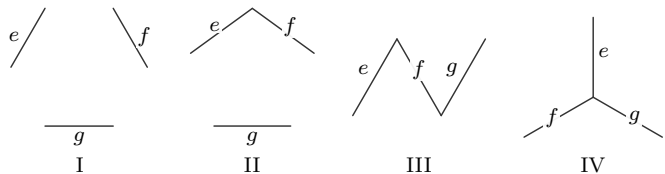
\begin{figure}[!b]
\centering
\setlength{\tabcolsep}{7pt}
\begin{tabular}{@{}cccc@{}}
\begin{tikzpicture}[scale=.82, baseline=(current bounding box.center), line cap=round, line join=round]
    \coordinate (A) at (-0.55,0.95);
    \coordinate (B) at (0.55,0.95);
    \coordinate (C) at (1.10,0.00);
    \coordinate (D) at (0.55,-0.95);
    \coordinate (E) at (-0.55,-0.95);
    \coordinate (F) at (-1.10,0.00);
    \draw[potB] (F)--(A);
    \putlab{(-1.05,0.50)}{$e$}
    \draw[potB] (B)--(C);
    \putlab{(1.05,0.50)}{$f$}
    \draw[potB] (E)--(D);
    \putlab{(0,-1.15)}{$g$}
\end{tikzpicture}
&
\begin{tikzpicture}[scale=.86, baseline=(current bounding box.center), line cap=round, line join=round]
    \coordinate (T)  at (0,1.00);
    \coordinate (L)  at (-0.95,0.31);
    \coordinate (R)  at (0.95,0.31);
    \coordinate (BR) at (0.59,-0.81);
    \coordinate (BL) at (-0.59,-0.81);
    \draw[potB] (L)--(T);
    \putlab{(-0.58,0.72)}{$e$}
    \draw[potB] (T)--(R);
    \putlab{(0.58,0.72)}{$f$}
    \draw[potB] (BL)--(BR);
    \putlab{(0,-1.00)}{$g$}
\end{tikzpicture}
&
\begin{tikzpicture}[scale=.78, baseline=(current bounding box.center), line cap=round, line join=round]
    \coordinate (A) at (0,0);
    \coordinate (B) at (0.75,1.299);
    \coordinate (C) at (1.50,0);
    \coordinate (D) at (2.25,1.299);
    \draw[potB] (A)--(B);
    \putlab{(0.18,0.78)}{$e$}
    \draw[potB] (B)--(C);
    \putlab{(1.12,0.78)}{$f$}
    \draw[potB] (C)--(D);
    \putlab{(1.68,0.78)}{$g$}
\end{tikzpicture}
&
\begin{tikzpicture}[scale=.88, baseline=(current bounding box.center), line cap=round, line join=round]
    \coordinate (O) at (0,0);
    \coordinate (A) at (90:1.20);
    \coordinate (B) at (210:1.20);
    \coordinate (C) at (330:1.20);
    \draw[potB] (O)--(A);
    \putlab{(0.16,0.68)}{$e$}
    \draw[potB] (O)--(B);
    \putlab{(-0.62,-0.30)}{$f$}
    \draw[potB] (O)--(C);
    \putlab{(0.62,-0.30)}{$g$}
\end{tikzpicture}
\\[1mm]
\scriptsize $\mathrm{I}$ & \scriptsize $\mathrm{II}$ &
\scriptsize $\mathrm{III}$ & \scriptsize $\mathrm{IV}$
\end{tabular}
\caption{The four overlap types of the underlying undirected arcs $e,f,g$.}
\label{fig:overlap-types}
\end{figure}

\begin{description}[leftmargin=2.8em, labelwidth=2.2em, style=nextline]

\item[$\mathrm{I}$:] The arcs $e,f,g$ are pairwise node-disjoint.
Equivalently, the union of the three arcs contains six distinct nodes.

\item[$\mathrm{II}$:] Exactly one pair among $e,f,g$ shares a common node,
while the remaining arc is node-disjoint from the other two.
Equivalently, the union of the three arcs contains five distinct nodes.

\item[$\mathrm{III}$:] The arcs $e,f,g$ form a path of length three.
Equivalently, the three arcs have exactly two distinct shared nodes, and their
union contains four distinct nodes.

\item[$\mathrm{IV}$:] The three arcs share a single common node.
Equivalently, they form a three-star whose union contains four distinct nodes.

\end{description}
We omit the remaining overlap pattern in which the three underlying arcs form a
triangle.
Indeed, for any choice of orientations, two of the three directed arcs lie on a
directed $3$-cycle.
Thus, such triples make no contribution to $\mathbb{E}H_{n,p}$.

Let $\mathcal{M}_i$, $i=1,2,3,4$, denote the sets of triples of directed arcs of
types $\mathrm{I}$, $\mathrm{II}$, $\mathrm{III}$, and $\mathrm{IV}$,
respectively.
Define $H_{n,p}^{(i)}=\sum_{\alpha\in\mathcal{M}_i}I_\alpha, i=1,2,3,4.$ Every
triple $\alpha \in \mathcal{M}$ that contributes to $H_{n,p}$ belongs to exactly
one of $\mathcal{M}_1,\ldots,\mathcal{M}_4$.
Hence,
\begin{equation}\label{eq:H-type-decomposition}
  H_{n,p}=\sum_{i=1}^4 H_{n,p}^{(i)}.
\end{equation}
We compute the expectations of the four terms in (\ref{eq:H-type-decomposition})
separately.
\par\medskip
\noindent\textbf{The term $\E H_{n,p}^{(1)}$.}
We first consider the triples in $\mathcal{M}_1$.
In this case, the arcs $e, f, g$ are pairwise node-disjoint.
Hence no two of these arcs can lie on a directed $3$-cycle, and the three arcs
cannot all lie on a directed $5$-cycle.
Therefore, it remains only to consider possible directed $4$-cycles and
$6$-cycles.

One way to obtain such a configuration is as follows.
Choose six nodes, partition them into three unordered pairs, and orient each of
the three resulting arcs in one of the two possible directions.
Hence the number of such configurations is
\[
  \binom n6\frac{\binom62\binom42}{3!}\,2^3=\frac16\fall n6.
\]
By symmetry, $\mathbb{E}I_\alpha$ is the same for all $\alpha\in\mathcal{M}_1$.
Hence we have
\[
  \E H_{n,p}^{(1)}=\frac16\fall n6\E I_\alpha.
\]
We compute $\E I_\alpha$ for $\alpha\in\mathcal{M}_1$.
For the subsequent analysis, let $\bara$ denote the set of arcs that may affect
the value of $I_\alpha$, namely, those that may form a directed $4$-cycle with
two arcs of $\alpha$ or a directed $6$-cycle with all three.
Without loss of generality, let
$\alpha=\{\overrightarrow e,\overrightarrow f,\overrightarrow g\}$, and label
the relevant surrounding arcs as in Figure~\ref{fig:caseI-potential-en}.

\begin{figure}[!b]
\centering
\FigCaseIPotential
\caption{$\alpha=\{\dir e,\dir f,\dir g\}, \bara=\{x_1,x_2,x_3,x_4,x_5,x_6\}$.}
\label{fig:caseI-potential-en}
\end{figure}

\noindent The possible local obstruction types involving $\alpha$ and
$\bar\alpha$ are shown in Figure~\ref{fig:caseI-allcases-en}.
Irrelevant arcs are omitted, and configurations that give rise to the same
obstruction are grouped together.

\begin{figure}[!t]
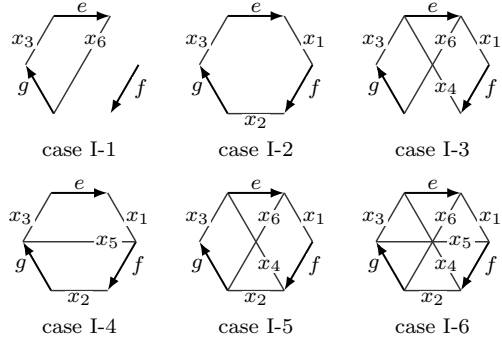

\centering
\begin{tabular}{ccc}
\FigCaseIone & \FigCaseItwo & \FigCaseIthree\\[1mm]
\scriptsize case I-1 & \scriptsize case I-2 & \scriptsize case I-3\\[1.5mm]
\FigCaseIfour & \FigCaseIfive & \FigCaseIsix\\[1mm]
\scriptsize case I-4 & \scriptsize case I-5 & \scriptsize case I-6
\end{tabular}
\caption{The six local obstruction types in Case I.}
\label{fig:caseI-allcases-en}
\end{figure}
\smallskip
\noindent\textbf{Case I-1.} This structure forms only one $4$-cycle. Put
$X_1=\min\{e,g,x_3,x_6\},$ and let $A_1=\{X_1\in\{e,g\}\}.$ Conditioned on the
existence of the arcs involved in this configuration, the probability that
$\alpha$ generates a $2$-face is $\Prob(A_1)=\frac12.$

\smallskip
\noindent\textbf{Case I-2.} This structure forms only one $6$-cycle. Put
$X_2=\min\{e,f,g,x_1,x_2,x_3\},$ and let $A_2=\{X_2\in\{e,f,g\}\}.$ Conditioned
on the existence of the arcs involved in this configuration, the probability
that $\alpha$ generates a $2$-face is $\Prob(A_2)=\frac12.$

\smallskip
\noindent\textbf{Case I-3.} This structure forms two $4$-cycles.  Put
$X_3=\min\{e,f,g,x_1,x_3,x_4,x_6\},$ and let
$A_3=\{\min\{e,f,x_1,x_4\}\in\{e,f\}\ \text{and}\ \min\{e,g,x_3,x_6\}\in\{e,g\}\}.$
Conditioned on the existence of the arcs involved in this configuration, the
probability that $\alpha$ generates a $2$-face is
\begin{align*}
  \Prob(A_3) &=\Prob(X_3=e)\Prob(A_3\mid X_3=e)+\Prob(X_3=f)\Prob(A_3\mid X_3=f)\\
  &\quad+\Prob(X_3=g)\Prob(A_3\mid X_3=g)\\
  &=\frac17\cdot1+\frac17\cdot\frac12+\frac17\cdot\frac12=\frac27.
\end{align*}
\smallskip
\noindent\textbf{Case I-4.} This structure forms one $6$-cycle and one $4$-cycle.  Put
$X_4=\min\{e,\allowbreak f,\allowbreak g,\allowbreak x_1,\allowbreak
x_2,\allowbreak x_3,\allowbreak x_5\},$ and let $A_4=\{\min\{f,\allowbreak
g,\allowbreak x_2,\allowbreak
x_5\}\in\nobreak\{f,g\}\ \text{and}\ \min\{e,\allowbreak f,\allowbreak
g,\allowbreak x_1,\allowbreak x_2,\allowbreak x_3\}\in\nobreak\{e,f,g\}\}.$
Conditioned on the existence of the arcs involved in this configuration, the
probability that $\alpha$ generates a $2$-face is
\begin{align*}
  \Prob(A_4) &=\Prob(X_4=e)\Prob(A_4\mid X_4=e)+\Prob(X_4=f)\Prob(A_4\mid X_4=f)\\
  &\quad+\Prob(X_4=g)\Prob(A_4\mid X_4=g)\\
  &=\frac17\cdot\frac12+\frac17\cdot1+\frac17\cdot1=\frac5{14}.
\end{align*}
\smallskip
\noindent\textbf{Case I-5.} This structure forms one $6$-cycle and two $4$-cycles. Put
$X_5=\min\{e,\allowbreak f,\allowbreak g,\allowbreak x_1,\allowbreak
x_2,\allowbreak x_3,\allowbreak x_4,\allowbreak x_6\},$ and let
$A_5=\{\min\{e,\allowbreak f,\allowbreak g,\allowbreak x_1,\allowbreak
x_2,\allowbreak x_3\}\in\nobreak\{e,f,g\},\allowbreak \ \min\{e,\allowbreak
f,\allowbreak x_1,\allowbreak x_4\}\in\nobreak\{e,f\},\allowbreak
\min\{e,\allowbreak g,\allowbreak x_3,\allowbreak x_6\}\in\nobreak\{e,g\}\}.$
Conditioned on the existence of the arcs involved in this configuration, the
probability that $\alpha$ generates a $2$-face is
\begin{align*}
  \Prob(A_5) &=\Prob(X_5=e)\Prob(A_5\mid X_5=e)+\Prob(X_5=f)\Prob(A_5\mid X_5=f)\\
  &\quad+\Prob(X_5=g)\Prob(A_5\mid X_5=g)\\
  &=\frac18\cdot1+\frac18\cdot\frac12+\frac18\cdot\frac12=\frac14.
\end{align*}
\smallskip
\noindent\textbf{Case I-6.} This structure forms two $6$-cycles and three $4$-cycles.  Put
$X_6=\min\{e,\allowbreak f,\allowbreak g,\allowbreak x_1,\allowbreak
x_2,\allowbreak x_3,\allowbreak x_4,\allowbreak x_5,\allowbreak x_6\}.$ and let
{\thinmuskip=3mu plus 2mu minus 1mu \thickmuskip=5mu plus 2mu minus 1mu
$A_6=\{\min\{e,\allowbreak f,\allowbreak g,\allowbreak x_1,\allowbreak
x_2,\allowbreak x_3\}\in\nobreak\{e,f,g\},\allowbreak \min\{e,\allowbreak
f,\allowbreak g,\allowbreak x_4,\allowbreak x_5,\allowbreak
x_6\}\in\nobreak\{e,f,g\},\allowbreak \min\{e,\allowbreak f,\allowbreak
x_1,\allowbreak x_4\}\in\nobreak\{e,f\},\allowbreak \min\{f,\allowbreak
g,\allowbreak x_2,\allowbreak x_5\}\in\nobreak\{f,g\},\allowbreak
\min\{e,\allowbreak g,\allowbreak x_3,\allowbreak x_6\}\in\nobreak\{e,g\}\}.$}
Conditioned on the existence of the arcs involved in this configuration, the
probability that $\alpha$ generates a $2$-face is
\begin{align*}
  \Prob(A_6) &=\Prob(X_6=e)\Prob(A_6\mid X_6=e)+\Prob(X_6=f)\Prob(A_6\mid X_6=f)\\
  &\quad+\Prob(X_6=g)\Prob(A_6\mid X_6=g)\\
  &=\frac19\cdot\frac12+\frac19\cdot\frac12+\frac19\cdot\frac12=\frac16.
\end{align*}
Let $K$ be the number of arcs in $\bar\alpha$ that are present in $G_{n,p}$, and
classify according to the value of $K$ as follows.
\begin{enumerate}[label=(\roman*)]
    \item If $K=0$ or $K=1$, there exists no structures listed above.

    \item If $K=2$, three choices correspond to Case~I-1, whereas the
    remaining twelve do not trigger an obstruction.

    \item If $K=3$, twelve choices correspond to Case~I-1, two to
    Case~I-2, and the remaining six do not trigger an obstruction.

    \item If $K=4$, six choices correspond to Case~I-1, three to
    Case~I-3, and the remaining six to Case~I-4.

    \item If $K=5$, each of the six choices corresponds to Case~I-5.

    \item If $K=6$, the unique choice corresponds to Case~I-6.
\end{enumerate}
The preceding classification gives
\begin{align*}
  \E(I_\alpha\mid K=0)&=p^3,\qquad \E(I_\alpha\mid K=1)=p^3,\\
  \E(I_\alpha\mid K=2)&=\frac{3\Prob(A_1)+12}{15}p^3=\frac{9}{10}p^3,\\
  \E(I_\alpha\mid K=3)&=\frac{12\Prob(A_1)+2\Prob(A_2)+6}{20}p^3=\frac{13}{20}p^3,\\
  \E(I_\alpha\mid
  K=4)&=\frac{6\Prob(A_1)+3\Prob(A_3)+6\Prob(A_4)}{15}p^3=\frac{2}{5}p^3,\\
  \E(I_\alpha\mid K=5)&=\Prob(A_5)p^3=\frac{1}{4}p^3,\\
  \E(I_\alpha\mid K=6)&=\Prob(A_6)p^3=\frac{1}{6}p^3.
\end{align*}
By the law of total expectation,
\begin{align*}
  \E I_\alpha &=\sum_{k=0}^6 \Prob(K=k)\E(I_\alpha\mid K=k)\\
  &=\sum_{k=0}^6\binom6k p^k(1-p)^{6-k}\E(I_\alpha\mid K=k)\\
  &=p^3\Bigl[(1-p)^6+6p(1-p)^5+\frac{27}{2}p^2(1-p)^4+13p^3(1-p)^3\\
  &\qquad +6p^4(1-p)^2+\frac32p^5(1-p)+\frac16p^6\Bigr]\\
  &=p^3\Bigl(1-\frac32p^2-p^3+3p^4-\frac32p^5+\frac16p^6\Bigr).
\end{align*}
Summing over $\alpha\in\mathcal M_1$ gives
\begin{equation}\label{eq:caseI-contribution-full-en}
  \E H_{n,p}^{(1)}=\sum_{\alpha\in\mathcal{M}_1}\E I_\alpha =\frac16\fall n6 p^3H_1(p).
\end{equation}
\medskip
\noindent\textbf{The term $\E H_{n,p}^{(2)}$.}
We then consider the triples in $\mathcal{M}_2$.
In this case, exactly two of the arcs $\{e,f,g\}$ share a node, while the
remaining arc is node-disjoint from them.
We denote the adjacent arcs by $e$ and $f$, and the disjoint arc by $g$.

There are $\binom n5$ ways to choose 5 nodes, $\binom 52$ ways to divide them
into two blocks of sizes $3,2$, and $3$ ways to choose the common node of the
$2$-path.
After that, we split $\mathcal M_2$ into two subcases according to whether the
arcs $e$ and $f$ have the same or opposite orientations at their common node.
Let $\mathcal M_{2}^{(1)}$ denote the same-orientation case, and let
$\mathcal M_{2}^{(2)}$ denote the opposite-orientation case.
Then
\[
  \mathcal M_2=\mathcal M_{2}^{(1)}\sqcup \mathcal M_{2}^{(2)}.
\]
The number of possible configurations in each of $\mathcal M_{2}^{(1)}$ and
$\mathcal M_{2}^{(2)}$ is
\[
  \binom n5\binom 52\cdot3\cdot4=\fall n5.
\]

\smallskip
\noindent\textbf{Arcs with the Same Orientation at the Common Node.}
We begin by considering triples in $\mathcal M_2^{(1)}$.
In this case, it suffices to consider directed $4$-cycles containing either
$\{e,g\}$ or $\{f,g\}$.
Without loss of generality, let
$\alpha=\{\overrightarrow e,\overrightarrow f,\overrightarrow g\}$, and label
the relevant surrounding arcs as in Figure~\ref{fig:IIa-base-en}.

\begin{figure}[!b]
\centering
\FigIIaBase
\caption{$\alpha=\{\dir e,\dir f,\dir g\}, \bara=\{x_1,x_2,x_3\}$.}
\label{fig:IIa-base-en}
\end{figure}

\noindent The possible local obstruction types involving $\alpha$ and
$\bar\alpha$ are shown in Figure~\ref{fig:IIa-allcases-en}.
Irrelevant arcs are omitted, and configurations that give rise to the same
obstruction are grouped together.

\begin{figure}[!b]
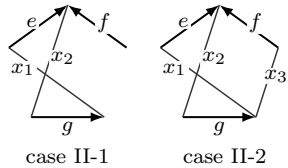

\centering
\begin{tabular}{cc}
\FigIIaOne & \FigIIaTwo\\[1mm]
\scriptsize case II-1 & \scriptsize case II-2
\end{tabular}
\caption{The two local obstruction types.}
\label{fig:IIa-allcases-en}
\end{figure}
\smallskip
\noindent\textbf{Case II-1.} This structure forms one $4$-cycle. Put
$Y_1=\min\{e,g,x_1,x_2\}$ and let $B_1=\{Y_1\in\{e,g\}\}.$ Conditioned on the
existence of the arcs involved in this configuration, the probability that
$\alpha$ generates a $2$-face is $\Prob(B_1)=\frac12.$

\smallskip
\noindent\textbf{Case II-2.} This structure forms two $4$-cycles.  Put
$Y_2=\min\{e,f,g,x_1,x_2,x_3\},$ and let
$B_2=\{\min\{e,g,x_1,x_2\}\in\{e,g\},\ \min\{f,g,x_2,x_3\}\in\{f,g\}\}.$
Conditioned on the existence of the arcs involved in this configuration, the
probability that $\alpha$ generates a $2$-face is
\begin{align*}
  \Prob(B_2) &=\Prob(Y_2=e)\Prob(B_2\mid Y_2=e)+\Prob(Y_2=f)\Prob(B_2\mid Y_2=f)\\
  &\quad+\Prob(Y_2=g)\Prob(B_2\mid Y_2=g)\\
  &=\frac16\cdot\frac12+\frac16\cdot\frac12+\frac16\cdot1=\frac13.
\end{align*}
\\
Let $K$ be the number of arcs in $\bar\alpha$ that appear in $G_{n,p}$, and
classify according to the value of $K$ as follows.
\begin{enumerate}[label=(\roman*)]
    \item If $K=0$ or $K=1$, there exists no structures listed above.

    \item If $K=2$, two choices correspond to Case~II-1, whereas the
    remaining choice does not trigger an obstruction.

    \item If $K=3$, the unique choice corresponds to Case~II-2.
\end{enumerate}
The preceding classification gives
\begin{align*}
  \E(I_\alpha\mid K=0)&=\E(I_\alpha\mid K=1)=p^3,\\
  \E(I_\alpha\mid K=2)&=\frac{2\Prob(B_1)+1}{3}p^3=\frac23p^3,\\
  \E(I_\alpha\mid K=3)&=\Prob(B_2)p^3=\frac13p^3.
\end{align*}
By the law of total expectation,
\begin{align*}
  \E I_\alpha &=\sum_{k=0}^3\Prob(K=k)\E(I_\alpha\mid K=k)\\
  &=p^3\Bigl[(1-p)^3+3p(1-p)^2+2p^2(1-p)+\frac13p^3\Bigr]\\
  &=p^3\Bigl(1-p^2+\frac13p^3\Bigr).
\end{align*}
Since the number of such configurations is $\fall n5$, the total contribution is
\begin{equation}\label{eq:caseIIa-contribution-full-en}
  \sum_{\alpha\in\mathcal M_2^{(1)}}\E I_\alpha=\fall n5p^3H_2(p).
\end{equation}
\smallskip
\noindent\textbf{Arcs with opposite orientations at the common node.}
We next consider triples in $\mathcal M_2^{(2)}$.
In this case, the possible obstructions may involve additional nodes beyond
those of $\alpha$; for instance, they may form a directed $4$-cycle together
with $\{\dir e,\dir f\}$, or a directed $6$-cycle together with
$\{\dir e,\dir f,\dir g\}$.
Since an exact enumeration of these configurations is not needed, we shall only
use an upper bound.
Without loss of generality, take $\alpha$ to be the configuration shown in
Figure~\ref{fig:IIb-base-en}.
\begin{figure}[!b]
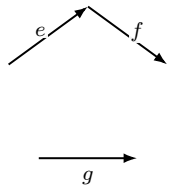

\centering
\FigIIbBase
\caption{The adjacent arcs have opposite orientations at the common endpoint.}
\label{fig:IIb-base-en}
\end{figure}

We bound $\E I_\alpha$, where $\alpha \in \mathcal{M}_2^{(2)}$.
Let $v$ be the common node of $e$ and $f$ and write $e=u_1v$ and $f=u_2v$.
Let $A$ be the event that the arc $u_1u_2$ is absent.
Let $B$ be the event that no directed 4-cycle containing $\dir e$ and $\dir f$
triggers the obstruction.
If $\alpha$ exists, then both $A$ and $B$ are necessary for
$\alpha=\{\dir e,\dir f,\dir g\}$ to generate a $2$-face.
Moreover, $A$ and $B$ are independent.
Hence we have
\[
  \E I_\alpha \leq p^3\mathbb{P}(A)\mathbb{P}(B) \leq p^3\mathbb{P}(B).
\]
It remains only to estimate $\mathbb{P}(B)$.
For each node $w_i\in V\setminus\{u_1,u_2,v\}$, there is a pair of arcs
$e_i=u_1w_i, f_i=u_2w_i$, such that $\{e_i,f_i\}$ forms a directed $4$-cycle
together with $\{\dir e,\dir f\}$.
Thus there are $n-3$ such pairs in total.
Equivalently, $B$ is the event that, for every $\{e_i,f_i,e,f\}$, which is an
existing directed 4-cycle containing $\{e,f\}$, we have
\begin{equation}\label{eq:4-cycle-rule}
  \min\{e_i,f_i,e,f\}\in\{e,f\}.
\end{equation}

We consider all possible underlying arcs on $[n]$,
regardless of their presence in $G_{n,p}$.
For $r=0,1,\ldots,n-3$, let $A_r$ be the event that
exactly $r$ of the $n-3$ potential directed $4$-cycles
$\{e_i,f_i,e,f\}$ violate (\ref{eq:4-cycle-rule}).
By symmetry of the uniform ordering on
$\binom{[n]}{2}$, we have
\[
  \mathbb{P}(A_r)=\frac{1}{n-2}.
\]
Conditioned on $A_r$, the event $B$ occurs if and only if
none of the $r$ corresponding pairs $\{e_i,f_i\}$
is simultaneously present in $G_{n,p}$.
Hence,
\[
  \mathbb{P}(B\mid A_r)=(1-p^2)^r.
\]
By the law of total probability,
\begin{align*}
  \mathbb{P}(B)
  &=\sum_{r=0}^{n-3}
    \mathbb{P}(B\mid A_r)\mathbb{P}(A_r)\\
  &=\frac1{n-2}\sum_{r=0}^{n-3}(1-p^2)^r\\
  &=\frac{1-(1-p^2)^{n-2}}{(n-2)p^2}\\
  &\le\min\left\{1,\frac2{np^2}\right\},
\end{align*}
where the last inequality holds for all sufficiently large $n$.

The number of such configurations is $\fall n5$.
Therefore,
\begin{equation}\label{eq:caseIIb-error-full-en}
  \sum_{\alpha\in\mathcal M_2^{(2)}} \E I_\alpha
  \le\fall n5 p^3 \min\left\{1,\frac{2}{np^2}\right\}
  =o\bigl(n^5p^{5/2}\bigr),\quad\text{as } n\to\infty.
\end{equation}
The last relation follows by splitting into the two regimes of $p$ as follows.
\\
If $p\leq \sqrt{2/n}$, then as $n\to\infty$,
\[
  \frac{n^5p^3}{n^5p^{5/2}}= \sqrt p\leq \left(\frac{2}{n}\right)^{1/4}\to0.
\]
If $p\geq \sqrt{2/n}$, then as $n\to\infty$,
\[
  \frac{n^5p^3\cdot 2/(np^2)}{n^5p^{5/2}}=\frac{2}{np^{3/2}}
  \leq \left(\frac{2}{n}\right)^{1/4}\to0.
\]
Combining (\ref{eq:caseIIa-contribution-full-en}) and
(\ref{eq:caseIIb-error-full-en}), we obtain
\begin{equation}\label{eq:caseII-contribution-full-en}
  \E H_{n,p}^{(2)} =\fall n5p^3H_2(p)+o(n^5p^{5/2}),\quad\text{as } n\to\infty.
\end{equation}

\medskip
\noindent\textbf{The terms $\E H_{n,p}^{(3)}$ and $\E H_{n,p}^{(4)}$.}
Each configuration of these cases involves $4$ nodes and $3$ arcs, hence the
number of such configurations is on the order of $n^4$.
Moreover, $\E I_\alpha\le p^3$ for every
$\alpha\in \mathcal{M}_3\sqcup \mathcal{M}_4$.
Hence
\[
  \E H_{n,p}^{(3)}+\E H_{n,p}^{(4)}=\sum_{\alpha\in\M_3\sqcup\M_4}\E I_\alpha\le
  Cn^4p^3,
\]
and
\[
  \frac{n^4p^3}{n^5p^{5/2}}\le \frac{\sqrt p}{n}\le \frac1n\to0,\quad\text{as }
  n\to\infty.
\]
It follows that, as $n\to\infty$,
\begin{equation}\label{eq:caseIII-IV-error-full-en}
  \E H_{n,p}^{(3)}=o(n^5p^{5/2}),\qquad \E H_{n,p}^{(4)}=o(n^5p^{5/2}).
\end{equation}
Combining (\ref{eq:caseI-contribution-full-en}),
(\ref{eq:caseII-contribution-full-en}), and (\ref{eq:caseIII-IV-error-full-en})
proves the lemma.
\end{proof}

\begin{proof}[Proof of Theorem~\ref{thm:expectation}]
By Lemma~\ref{lem:two-face-decomposition}, we have
$\E F_{n,p}=\E H_{n,p}+\E K_{n,p}-2\E E_{n,p}.$ Moreover,
Lemma~\ref{lem:H-expectation} gives,
\[
  \E H_{n,p}=\frac16\fall n6 p^3H_1(p)+\fall n5p^3H_2(p) +o\bigl(n^5p^{5/2}\bigr),
\]
while (\ref{eq:KminusE-full-en}) yields,
\[
  \E K_{n,p}-2\E E_{n,p}=o\bigl(n^5p^{5/2}\bigr).
\]
Combining the last two estimates, we obtain
\[
  \E F_{n,p}= \frac16\fall n6 p^3H_1(p)+\fall n5p^3H_2(p)+ o\bigl(n^5p^{5/2}\bigr).
\]
This proves Theorem~\ref{thm:expectation}.
\end{proof}
\subsection{Variance}\label{sec:variance}
We proceed to estimate the variance of $F_{n,p}$.
Recall from Section~\ref{sec:expectation} that the non-origin face count
$H_{n,p}$ is represented as a sum of indicators indexed by triples of directed
arcs.
The main task is therefore to control the covariance terms arising from pairs of
such configurations.
These terms are treated separately according to the possible intersections of
the two configurations.
The final estimate reads as follows.
\begin{theorem}\label{thm:variance}
If $p\gg n^{-2}$ and $\limsup\limits_{n\to\infty}p<1$, then as $n\to\infty$,
\begin{equation}\label{eq:Var-main-full-en}
  \Var(F_{n,p})=\frac18\fall n{10}p^5q^5H_3(p)^2+o(n^{10}p^5),
\end{equation}
where
\begin{equation}\label{eq:H3-full-en}
  H_3(p)=p^4-6p^3+p^2+4p+2>2, \quad\text{for } 0<p<1.
\end{equation}
Consequently, no degeneracy occurs.
\end{theorem}

To simplify the counting argument, we assign i.i.d.
uniform weights $w_e\sim U(0,1)$ to all arcs $e\in\binom{[n]}2$, and the order
is determined by $e\prec f \Longleftrightarrow  w_e<w_f.$ These weights are also
independent of $G_{n,p}$.
Let $w=\{w_e: e\in \binom{[n]}{2}\}$ denote the collection of weights assigned
to all arcs.
By comparing the values of the weights $w_e$, we obtain a random term order.
Hence $w$ naturally induces a total order on the arcs of $G_{n,p}$.
For simplicity, once a total order $w$ is fixed, we still write $F_{n,p}$ for
the number of $2$-faces of the triangulation $\Delta$ induced by $w$.
Since the value of $F_{n,p}$ does not depend on the total order on the arcs of
$G_{n,p}$, and since $w$ is independent of $G_{n,p}$, the conditional
expectation $\E(F_{n,p}\mid w)$ is constant in $w$.
Hence, by the law of total variance,
\begin{equation}\label{eq:total-variance}
  \Var(F_{n,p})=\E\bigl[\Var(F_{n,p}\mid w)\bigr].
\end{equation}
By Lemma~\ref{lem:two-face-decomposition} and (\ref{eq:H-type-decomposition}), we
have
\begin{equation}\label{eq:F-variance-decomposition}
  F_{n,p}=\sum_{i=1}^4H_{n,p}^{(i)}+K_{n,p}-2E_{n,p}.
\end{equation}
Expanding the conditional variance gives
\begin{equation}\label{eq:variance-decomposition}
  \begin{aligned}
    \Var(F_{n,p}) &=\E\bigl[\Var(F_{n,p}\mid w)\bigr] \\
    =&\EE \Biggl[ \sum_{i=1}^{4}\Var\left(H_{n,p}^{(i)}\mid w\right) + \Var
    (K_{n,p}\mid w) + 4\Var (E_{n,p}\mid w) \\
    & +2\sum_{1\le i<j\le 4} \Cov \left(H_{n,p}^{(i)},H_{n,p}^{(j)}\mid w\right) +
    2\sum_{i=1}^{4} \Cov \left(H_{n,p}^{(i)},K_{n,p}\mid w\right) \\
    & -4\sum_{i=1}^{4} \Cov \left(H_{n,p}^{(i)},E_{n,p}\mid w\right) - 4\Cov
    \left(K_{n,p},E_{n,p}\mid w\right) \Biggr].
  \end{aligned}
\end{equation}
We estimate each term in (\ref{eq:variance-decomposition}) separately.
It turns out that the leading contribution comes from the term
$\Var(H_{n,p}^{(1)}\mid w)$, which is evaluated in the following proposition.

\begin{proposition}\label{prop:variance-H1}
Under the assumptions of Theorem~\ref{thm:variance}, as $n\to\infty$, we have
\begin{equation}\label{eq:variance-H1-main}
  \EE [\Var(H_{n,p}^{(1)}\mid w)] =\frac18\fall n{10}p^5q^5H_3(p)^2+o(n^{10}p^5).
\end{equation}
\end{proposition}

\begin{proof}
Expanding the conditional variance of the indicator sum gives
\[
  \Var(H_{n,p}^{(1)}\mid w)
  = \sum_{\alpha\in\mathcal M_1} \sum_{\beta\in\mathcal M_1} \Cov(I_\alpha,I_\beta\mid
  w).
\]
Taking expectations easily gives
\[
  \E\bigl[\Var(H_{n,p}^{(1)}\mid w)\bigr]
  = \sum_{\alpha\in\mathcal M_1} \sum_{\beta\in\mathcal M_1}
  \E\bigl[\Cov(I_\alpha,I_\beta\mid w)\bigr].
\]
To estimate the covariance sum, for any set $\gamma$ of directed arcs, let
$V(\gamma)$ denote the set of nodes incident with the arcs of $\gamma$.
We rearrange the covariance sum according to the size of
$V(\alpha)\cap V(\beta)$ as follows:
\begin{equation}\label{eq:H1-overlap-decomposition}
  \sum_{\alpha,\beta\in\mathcal M_1} \E\bigl[\Cov(I_\alpha,I_\beta\mid w)\bigr]
  =\sum_{k=0}^{6} \sum_{\substack{\alpha,\beta\in\mathcal M_1\\
  |V(\alpha)\cap V(\beta)|=k}} \E\bigl[\Cov(I_\alpha,I_\beta\mid w)\bigr].
\end{equation}
We analyze the covariance sum in (\ref{eq:H1-overlap-decomposition}) by
distinguishing the following three overlap regimes:
\[
  \text{(i)}\quad |V(\alpha)\cap V(\beta)|\leq 1,
  \qquad \text{(ii)}\quad |V(\alpha)\cap V(\beta)|=2,
  \qquad \text{(iii)}\quad |V(\alpha)\cap V(\beta)|\geq 3.
\]

\noindent
\textbf{\boldmath$|V(\alpha)\cap V(\beta)|\leq 1$:} In this regime, $I_\alpha$
and $I_\beta$ are conditionally independent given $w$, and hence
$\operatorname{Cov}(I_\alpha,I_\beta\mid w)=0.$ Thus this regime makes no
contribution to the covariance sum.

\medskip
\noindent
\textbf{\boldmath$|V(\alpha)\cap V(\beta)|=2$:} We next consider the two-node
overlap, which gives the leading contribution to the covariance sum.
Write
\[
  \alpha=\{\overrightarrow e,\overrightarrow f,\overrightarrow g\},
  \qquad \beta=\{\overrightarrow k,\overrightarrow l,\overrightarrow m\}.
\]
Up to symmetry, there are three possible overlap configurations, shown in
Figure~\ref{fig:variance-threecases-en}.
\begin{figure}[!t]
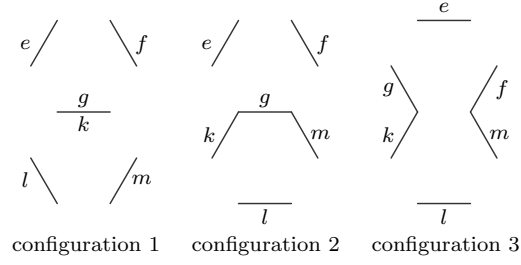

\centering
\begin{tabular}{ccc}
\VCaseOne & \VCaseTwo & \VCaseThree\\[1mm]
\scriptsize configuration 1 & \scriptsize configuration 2 & \scriptsize configuration 3
\end{tabular}
\caption{The three possible overlap configurations.}
\label{fig:variance-threecases-en}
\end{figure}

Let $u$ and $v$ denote the two common nodes of $\alpha$ and $\beta$, set $s=uv$,
and let $w_s$ be the random weight of $s$.
For all three configurations, the dependence between $I_\alpha$ and $I_\beta$ is
transmitted through the arc $s$.
We first reduce the covariance calculation by conditioning on $w_s$.
The required reduction follows from the conditional covariance identity below.

\begin{lemma}\label{lem:variance}
If $\mathcal{F} \subset \mathcal{G}$ are $\sigma$-fields, then
\begin{equation}\label{eq:cov-decomposition}
  \Cov(X,Y\mid\mathcal{F})=\EE\left[\Cov(X,Y\mid\mathcal{G})\mid\mathcal{F}\right]+
  \Cov\left(\EE[X\mid\mathcal{G}],\EE[Y\mid\mathcal{G}]\mid\mathcal{F}\right).
\end{equation}
\end{lemma}

\begin{proof}
According to the definition of conditional covariance,
\[
  \Cov(X,Y\mid\mathcal{G})=
  \EE[XY\mid\mathcal{G}]-\EE[X\mid\mathcal{G}]\EE[Y\mid\mathcal{G}].
\]
By the property of conditional expectation,
\[
  \EE\left[\EE[XY\mid\mathcal{G}]\mid\mathcal{F}\right]=\EE[XY\mid\mathcal{F}].
\]
Taking conditional expectations in the covariance identity gives
\[
  \EE\left[\Cov(X,Y\mid\mathcal{G})\mid\mathcal{F}\right]=\EE[XY\mid\mathcal{F}]
  -\EE\left[\EE[X\mid\mathcal{G}]\EE[Y\mid\mathcal{G}]\mid\mathcal{F}\right].
\]
For the covariance of the conditional expectations, we have
\[
  \begin{aligned}
    &\Cov\left( \EE[X\mid\mathcal{G}], \EE[Y\mid\mathcal{G}]\mid\mathcal{F} \right) \\
    =&\EE\left[\EE[X\mid\mathcal{G}]\EE[Y\mid\mathcal{G}]\mid\mathcal{F}\right]
    -\EE\left[\EE[X\mid\mathcal{G}]\mid\mathcal{F}\right]
    \EE\left[\EE[Y\mid\mathcal{G}]\mid\mathcal{F}\right] \\
    =&\EE\left[\EE[X\mid\mathcal{G}]\EE[Y\mid\mathcal{G}]\mid\mathcal{F}
    \right]-\EE[X\mid\mathcal{F}]\mathbb{E}[Y\mid\mathcal{F}].
  \end{aligned}
\]
The second equality is obtained in the same way from the properties of
conditional expectation.
Adding the two equations above yields
\[
  \begin{aligned}
    &\EE\left[\Cov(X,Y\mid\mathcal{G})\mid\mathcal{F}\right]
    +\Cov\left(\EE[X\mid\mathcal{G}],\EE[Y\mid\mathcal{G}]\mid\mathcal{F}\right) \\
    =&\EE[XY\mid\mathcal{F}]-\mathbb{E}[X\mid\mathcal{F}] \EE[Y\mid\mathcal{F}] \\
    =&\Cov(X,Y\mid\mathcal{F}).
  \end{aligned}
\]
This completes the proof.
\end{proof}
In the present setting, let $\mathcal F:=\sigma(w_s), \mathcal G:=\sigma(w).$
Then $\mathcal F\subset\mathcal G$.
Applying Lemma~\ref{lem:variance} to $X=I_\alpha$ and $Y=I_\beta$, we obtain
\[
  \Cov(I_\alpha,I_\beta\mid w_s)
  = \E\bigl[\Cov(I_\alpha,I_\beta\mid w)\mid w_s\bigr] + \Cov\bigl(\E(I_\alpha\mid w),
  \E(I_\beta\mid w)\mid w_s\bigr).
\]
Let $w_\alpha$ and $w_\beta$ denote the collections of all other weights
relevant to $I_\alpha$ and $I_\beta$, respectively.
We shall omit all weights which are irrelevant to both indicators.
Since the arc weights are independent and the two sets of relevant arcs are
disjoint apart from $s$, the sigma-fields $\sigma(w_\alpha)$ and
$\sigma(w_\beta)$ are conditionally independent given $\sigma(w_s)$.
Moreover, $\E(I_\alpha\mid w)$ is $\sigma(w_s,w_\alpha)$-measurable, while
$\E(I_\beta\mid w)$ is $\sigma(w_s,w_\beta)$-measurable.
Hence
\[
  \Cov\bigl(\E(I_\alpha\mid w),\E(I_\beta\mid w)\mid w_s\bigr)=0.
\]
Thus we have
\[
  \E[\Cov(I_\alpha,I_\beta\mid w)\mid w_s] =\Cov(I_\alpha,I_\beta\mid w_s).
\]
Taking expectations gives
\[
  \E[\Cov(I_\alpha,I_\beta\mid w)] = \E[\Cov(I_\alpha,I_\beta\mid w_s)].
\]
Since $w_s\sim U(0,1)$, we also have
\[
  \E[\Cov(I_\alpha,I_\beta\mid w_s)] = \int_0^1 \Cov(I_\alpha,I_\beta\mid w_s=t)\,dt.
\]
Let $I_s$ denote the indicator that the underlying arc $s$ is present in
$G_{n,p}$.
Conditioning on the presence or absence of $s$ gives
\begin{align*}
  \Cov(I_\alpha,I_\beta\mid w_s=t)
  &=\E(I_\alpha I_\beta\mid w_s=t) -\E(I_\alpha\mid w_s=t)\E(I_\beta\mid w_s=t)\\
  &=p\E(I_\alpha I_\beta\mid I_s=1,w_s=t)+q\E(I_\alpha I_\beta\mid I_s=0,w_s=t)\\
  &\quad-\bigl[p\E(I_\alpha\mid I_s=1,w_s=t)+q\E(I_\alpha\mid I_s=0,w_s=t)\bigr]\\
  &\qquad\times \bigl[p\E(I_\beta\mid I_s=1,w_s=t)+q\E(I_\beta\mid I_s=0,w_s=t)\bigr].
\end{align*}
Using the conditional independence of $I_\alpha$ and $I_\beta$ after
conditioning on $I_s$ and $w_s$, this becomes
\begin{align}
  \Cov(I_\alpha,I_\beta\mid w_s=t)
  &=pq\Bigl[\E(I_\alpha\mid I_s=1,w_s=t)-\E(I_\alpha\mid I_s=0,w_s=t)\Bigr]\notag\\
  &\quad\times \Bigl[\E(I_\beta\mid I_s=1,w_s=t)-\E(I_\beta\mid I_s=0,w_s=t)\Bigr].
  \label{eq:cov-factor-full-en}
\end{align}
This factorization shows that the contribution of a configuration is determined
by whether $s$ belongs to each triple.
Here and below, we write $s\in\alpha$ if $s$ is the underlying arc of a directed
arc in $\alpha$.
For $s\in\alpha$, define
\[
  f_a(t):=\E(I_\alpha\mid I_s=1,w_s=t)-\E(I_\alpha\mid I_s=0,w_s=t).
\]
For $s\notin\alpha$, define
\[
  f_b(t):=\E(I_\alpha\mid I_s=1,w_s=t)-\E(I_\alpha\mid I_s=0,w_s=t).
\]
Define the corresponding functions for $\beta$ in the same way.
Then by (\ref{eq:cov-factor-full-en}),
\[
  \EE[\Cov(I_\alpha,I_\beta\mid w_s)]=\int_0^1pq\,f_i(t)f_j(t)\dd t,
  \qquad i,j\in\{a,b\},
\]
where $i$ and $j$ are determined by whether $s$ belongs to $\alpha$ and $\beta$.

We now evaluate the three configurations in
Figure~\ref{fig:variance-threecases-en} in turn.
The first configuration consists of five pairwise node-disjoint arcs.
There are $\binom n{10}$ ways to choose $10$ nodes, and
$\frac{1}{5!}\binom{10}2\binom82\binom62\binom42\binom22$ ways to divide them
into five blocks of size $2$.
After that, there are $5$ choices for the common arc $s$, and $\binom42$ choices
for deciding which two of the remaining arcs are in the same triple as $s$.
In this configuration the two indicators are dependent.
The arcs have $2^6$ possible orientations.
Hence the total contribution of this type is
\begin{align*}
  &\binom n{10} \frac{\binom{10}2\binom82\binom62\binom42\binom22}{5!} \cdot
  5\cdot\binom42\cdot 2^6\, pq\int_0^1 f_a^2(t)\dd t \\
  &\quad =\frac12\fall n{10}pq\int_0^1 f_a^2(t)\dd t .
\end{align*}

For the second configuration, the union consists of one $3$-path on four nodes
and three additional arcs disjoint from this path.
There are $\binom n{10}$ ways to choose 10 nodes, and
$\frac{1}{3!}\binom{10}2\binom82\binom62$ ways to divide them into connected
blocks of sizes $4,2,2,2$.
After that, there are $\binom42$ choices for the middle arc and $2$ ways to
connect the remaining nodes.
There are $2$ choices for deciding whether the middle arc belongs to $\alpha$ or
to $\beta$, and another $3$ choices for deciding which two arcs other than $s$
are in the same triple as the middle arc, so there are $2 \times 3 = 6$ ways to
decide which pair of arcs belongs to $\alpha$.

The two indicators can be dependent only when this $3$-path is a subgraph of the
corresponding $4$-cycle or $6$-cycle.
The pair of arcs connected through this $3$-path has only $2$ possible
orientations, while the remaining four arcs have $2^4$ possible orientations.
Therefore the total contribution of this type is
\begin{align*}
  &\binom
  n{10}\frac{\binom{10}2\binom82\binom62}{3!}\binom42\cdot2\cdot2\cdot3\cdot2\cdot2^4\,pq\int_0^1f_a(t)f_b(t)\dd
  t\\
  &\quad=2\fall n{10}pq\int_0^1f_a(t)f_b(t)\dd t.
\end{align*}

Finally, for the third configuration, there are $\binom n{10}$ ways to choose 10
nodes, and $\frac{1}{2^2}\binom{10}2\binom82\binom62\binom32$ ways to group
them.
There are $3^2$ choices for the middle nodes of the two $2$-paths, and $2\cdot4$
choices for deciding which arcs belong to $\alpha$.
The two indicators can be dependent only when $s$ is a potential arc for both
$\alpha$ and $\beta$.
The two arcs adjacent to $s$ have $4$ possible orientations, and the remaining
arcs have another $2^2$ possible orientations.
Therefore the total contribution of this type is
\begin{align*}
  &\binom n{10} \frac{\binom{10}2\binom82\binom62\binom32}{2^2} \cdot 3^2\cdot
  2\cdot4\cdot4\cdot2^2\, pq\int_0^1 f_b^2(t)\dd t \\
  &\quad =2\fall n{10}pq\int_0^1 f_b^2(t)\dd t .
\end{align*}

Adding the three contributions gives
\begin{align}
  &\frac12\fall n{10}pq\int_0^1f_a^2(t)\dd t+2\fall n{10}pq\int_0^1f_a(t)f_b(t)\dd
  t+2\fall n{10}pq\int_0^1f_b^2(t)\dd t\notag\\
  =&\frac12\fall n{10}pq\int_0^1\bigl(f_a(t)+2f_b(t)\bigr)^2\dd t.
  \label{eq:two-node-overlap-contribution}
\end{align}
Thus the contribution in this overlap regime is reduced to evaluating
$f_a(t)+2f_b(t)$.

\begin{lemma}\label{lem:fa-fb-constant}
We have $f_a(t)+2f_b(t)\equiv f_a(0)+2f_b(0).$
\end{lemma}

\begin{proof}
Consider the common local configuration of the underlying undirected arcs shown
in Figure~\ref{fig:root-alpha-beta-en}.
\begin{figure}[!t]
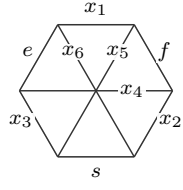

\centering
\ProofRootFig
\caption{The common underlying undirected configuration.}
\label{fig:root-alpha-beta-en}
\end{figure}
\noindent From this configuration we extract six directed triples, which we
denote by
\[
  \begin{aligned}
    \alpha_1&=\{\dir e,\dir f,\dir s\}, & \alpha_2&=\{\dir{x_1},\dir{x_4},\dir s\},\\
    \beta_1&=\{\dir{x_1},\dir{x_2},\dir{x_3}\},
    & \beta_2&=\{\dir{x_4},\dir{x_5},\dir{x_6}\},\\
    \beta_3&=\{\dir e,\dir{x_2},\dir{x_5}\}, & \beta_4&=\{\dir f,\dir{x_3},\dir{x_6}\}.
  \end{aligned}
\]
The orientations of these triples are shown in
Figure~\ref{fig:alpha-beta-six-en}.
\begin{figure}[!t]
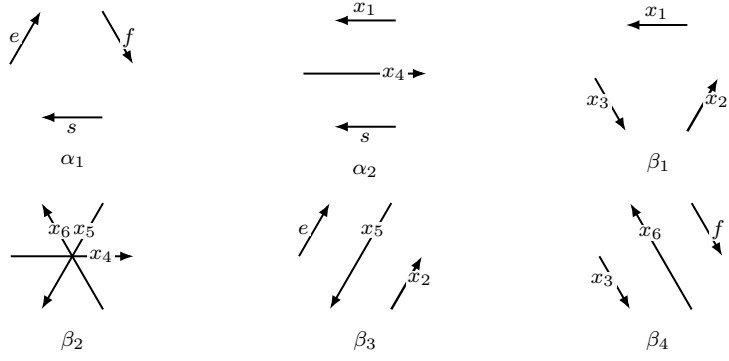

\centering
\ProofMini{\ProofAlphaOne}{$\alpha_1$} \ProofMini{\ProofAlphaTwo}{$\alpha_2$}
\ProofMini{\ProofBetaOne}{$\beta_1$}\\[3mm]
\ProofMini{\ProofBetaTwo}{$\beta_2$} \ProofMini{\ProofBetaThree}{$\beta_3$}
\ProofMini{\ProofBetaFour}{$\beta_4$}
\caption{
The six directed triples.
}
\label{fig:alpha-beta-six-en}
\end{figure}
Here $\alpha_1$ and $\alpha_2$ are the two triples containing the arc $s$, while
$\beta_1,\ldots,\beta_4$ are the four corresponding triples not containing $s$.

By symmetry of the local configuration,
\[
  \begin{aligned}
    2f_a(t)+4f_b(t) &= \E[I_{\alpha_1}+I_{\alpha_2}\mid I_s=1,w_s=t]  \\
    &\quad+ \E[I_{\beta_1}+I_{\beta_2}+I_{\beta_3}+I_{\beta_4}\mid I_s=1,w_s=t] \\
    &\quad- \E[I_{\beta_1}+I_{\beta_2}+I_{\beta_3}+I_{\beta_4}\mid I_s=0,w_s=t].
  \end{aligned}
\]
The last term is independent of $t$, since none of the triples
$\beta_1,\ldots,\beta_4$ contains the arc $s$.
Therefore it remains to show that
\[
  \E[I_{\alpha_1}+I_{\alpha_2}+I_{\beta_1}+I_{\beta_2}+I_{\beta_3}+I_{\beta_4} \mid
  I_s=1,w_s=t]
\]
is independent of $t$.

Let
$S:=I_{\alpha_1}+I_{\alpha_2} +I_{\beta_1}+I_{\beta_2}+I_{\beta_3}+I_{\beta_4}.$
We first show that $\E(S\mid w_s=t)$ is independent of $t$.
Indeed, after fixing all arcs and all weights except $w_s$, the six triples
above satisfy the following dichotomy.
If none of them exists, then $S=0$.
If at least one of them exists, then exactly one is counted as a $2$-face.
This is because any two distinct triples among
$\alpha_1,\alpha_2,\beta_1,\beta_2,\beta_3,\beta_4$ which appear simultaneously
are linked by a directed $4$-cycle or $6$-cycle.
The triangulation rule then allows at most one of them to survive.
On the other hand, as soon as at least one of these triples exists, the smallest
arc among the relevant local cycle determines one admissible triple.
Thus the value of $S$ is simply the indicator of the event that at least one of
the six triples exists.
In particular, this value does not depend on the position of $s$ in the order,
and hence does not depend on $t$.

Thus, $\E(S\mid w_s=t)$ is independent of $t$.
By conditioning further on whether the arc $s$ is present, we have
\[
  \begin{aligned}
    \E(S\mid w_s=t) &= p\,\E(S\mid I_s=1,w_s=t) + q\,\E(S\mid I_s=0,w_s=t).
  \end{aligned}
\]
When $I_s=0$, the two triples $\alpha_1$ and $\alpha_2$ cannot exist, so
\[
  \E(S\mid I_s=0,w_s=t)
  = \E[I_{\beta_1}+I_{\beta_2}+I_{\beta_3}+I_{\beta_4}\mid I_s=0,w_s=t],
\]
which is independent of $t$.
Since both $\E(S\mid w_s=t)$ and $\E(S\mid I_s=0,w_s=t)$ are independent of $t$,
it follows that $\E(S\mid I_s=1,w_s=t)$ is also independent of $t$.
Therefore $2f_a(t)+4f_b(t)$ is independent of $t$.
Dividing by $2$ gives $f_a(t)+2f_b(t)\equiv f_a(0)+2f_b(0).$ Here $t=0$ means
that the connecting arc $s$ is the smallest arc in the total order.
This proves the lemma.
\end{proof}

By Lemma~\ref{lem:fa-fb-constant}, the quantity $f_a(t)+2f_b(t)$ is independent of
$t$.
The lemma therefore reduces the remaining calculation to the endpoint $t=0$.
We evaluate the two terms separately, starting with
$f_a(0)=\E(I_\alpha\mid I_s=1,w_s=0), s\in\alpha$.
The possible local obstruction types involving $\alpha$ and $\bar\alpha$ are
shown in Figure~\ref{fig:fa-ten-cases-en}.
Irrelevant arcs are omitted, and configurations that give rise to the same
obstruction are grouped together.

\begin{figure}[!t]
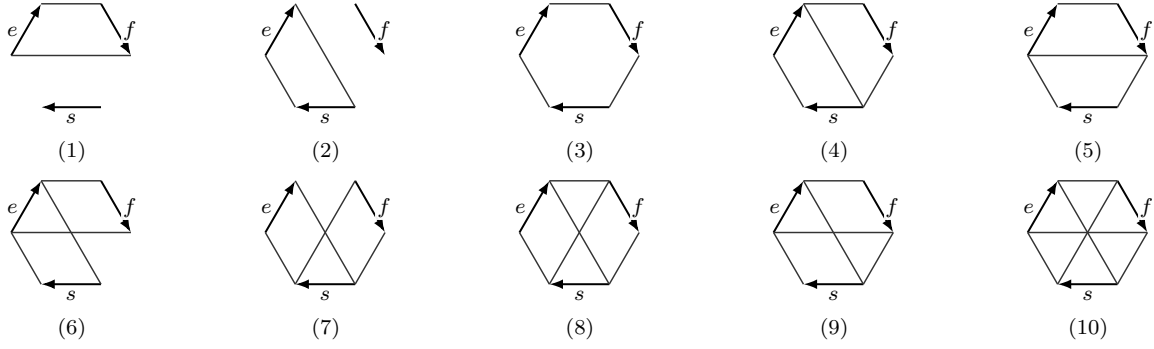

\centering
\Conf{\FAone}{1} \Conf{\FAtwo}{2} \Conf{\FAthree}{3} \Conf{\FAfour}{4}
\Conf{\FAfive}{5}\\[2mm]
\Conf{\FAsix}{6} \Conf{\FAseven}{7} \Conf{\FAeight}{8} \Conf{\FAnine}{9}
\Conf{\FAten}{10}
\caption{The ten local obstruction types for $f_a(0)$.}
\label{fig:fa-ten-cases-en}
\end{figure}

When $s$ is the smallest arc, cases~$(2)$, $(3)$, $(4)$, $(7)$, and $(8)$ do not
trigger an obstruction, whereas cases~$(1)$, $(5)$, $(6)$, $(9)$, and $(10)$ do
so with probability $1/2$.
Let $K$ be the number of arcs in $\bar{\alpha}\setminus\{s\}$ that are present
in $G_{n,p}$, and classify according to the value of $K$ as follows.

\begin{enumerate}[label=(\roman*)]
    \item If $K=0$ or $K=1$, there exists no structures listed above.

    \item If $K=2$, one choice corresponds to Case~$(1)$, two to Case~$(2)$,
    and the remaining twelve do not trigger an obstruction.

    \item If $K=3$, four choices correspond to Case~$(1)$, eight to
    Case~$(2)$, two to Case~$(3)$, and the remaining six do not trigger an
    obstruction.

    \item If $K=4$, two choices correspond to Case~$(1)$, four to Case~$(2)$,
    four to Case~$(4)$, two to Case~$(5)$, two to Case~$(6)$, and the remaining
    choice to Case~$(7)$.

    \item If $K=5$, two choices correspond to Case~$(8)$ and the remaining
    four to Case~$(9)$.

    \item If $K=6$, the unique choice corresponds to Case~$(10)$.
\end{enumerate}

\noindent Summing over the cases above gives
\begin{align}
  f_a(0)
  &=p^2\biggl(q^6+6pq^5+\Bigl(\frac12+2+12\Bigr)p^2q^4
  +\Bigl(4\cdot\frac12+8+2+6\Bigr)p^3q^3\notag\\
  &\quad+\Bigl(2\cdot\frac12+4+4+2\cdot\frac12+2\cdot\frac12+1\Bigr)p^4q^2
  +\Bigl(2+4\cdot\frac12\Bigr)p^5q+\frac12 p^6\biggr)\notag\\
  &=p^2q^6+6p^3q^5+\frac{29}{2}p^4q^4+18p^5q^3+12p^6q^2+4p^7q+\frac12 p^8.
  \label{eq:fa1}
\end{align}

\noindent It remains to compute
\[
  f_b(0)=\E(I_\alpha\mid I_s=1,w_s=0) -\E(I_\alpha\mid I_s=0,w_s=0),
  \qquad s\notin\alpha.
\]
We first evaluate the conditional expectation with $I_s=1$.
The possible local obstruction types involving $\alpha$ and $\bar\alpha$ are
shown in Figure~\ref{fig:fb-one-cases-en}.
Irrelevant arcs are omitted, and configurations that give rise to the same
obstruction are grouped together.

\begin{figure}[!t]
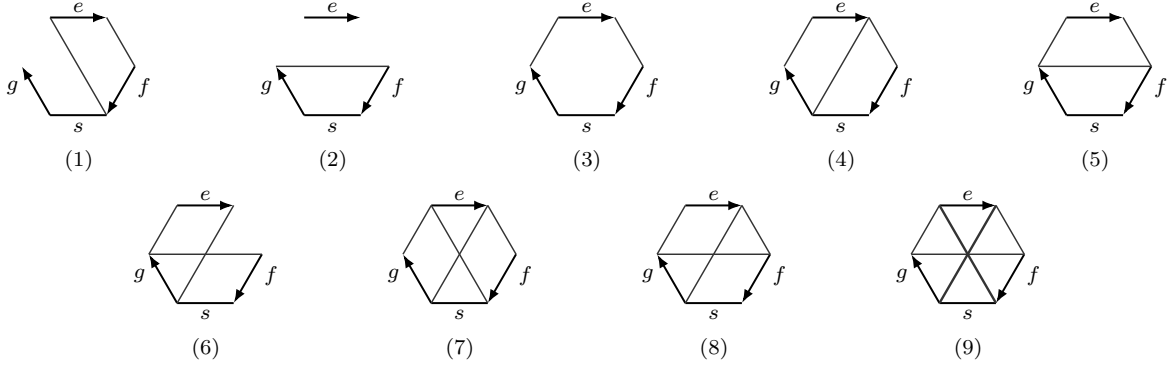

\centering
\Conf{\FBone}{1} \Conf{\FBtwo}{2} \Conf{\FBthree}{3} \Conf{\FBfour}{4}
\Conf{\FBfive}{5}\\[2mm]
\Conf{\FBsix}{6} \Conf{\FBseven}{7} \Conf{\FBeight}{8} \Conf{\FBnine}{9}
\caption{The nine local obstruction types for
$\E(I_\alpha\mid I_s=1,w_s=0)$.}
\label{fig:fb-one-cases-en}
\end{figure}

When $w_s=0$, cases~$(2)$, $(3)$, $(4)$, $(5)$, $(6)$, $(7)$, $(8)$, and $(9)$
always trigger an obstruction, whereas case~$(1)$ does so with probability
$1/2$.
Let $K$ be the number of arcs in $\bar{\alpha}\setminus\{s\}$ that are present
in $G_{n,p}$, and classify according to the value of $K$ as follows.

\begin{enumerate}[label=(\roman*)]
    \item If $K=0$, there exists no structures listed above.

    \item If $K=1$, one choice corresponds to Case~$(2)$, whereas the
    remaining four do not trigger an obstruction.

    \item If $K=2$, two choices correspond to Case~$(1)$, four to Case~$(2)$,
    one to Case~$(3)$, and the remaining three do not trigger an obstruction.

    \item If $K=3$, two choices correspond to each of Cases~$(1)$, $(2)$,
    $(4)$, $(5)$, and $(6)$.

    \item If $K=4$, one choice corresponds to Case~$(7)$, whereas the
    remaining four correspond to Case~$(8)$.

    \item If $K=5$, the unique choice corresponds to Case~$(9)$.
\end{enumerate}

\noindent Summing over the cases above gives
\begin{align}
  \E(I_\alpha\mid I_s=1,w_s=0)
  &=p^3\biggl(q^5+4pq^4+\Bigl(2\cdot \frac12+4\cdot 0+0+3\Bigr)p^2q^3\notag\\
  &\quad+\Bigl(2\cdot \frac12+2\cdot 0+2\cdot 0+2\cdot 0+2\cdot 0\Bigr)p^3q^2\notag\\
  &\quad+(0+4\cdot 0)p^4q+0\cdot p^5\biggr)\notag\\
  &=p^3\bigl(q^5+4pq^4+4p^2q^3+p^3q^2\bigr).
  \label{eq:fb1}
\end{align}

\noindent We then evaluate the conditional expectation with $I_s=0$.
The possible local obstruction types involving $\alpha$ and $\bar\alpha$ are
shown in Figure~\ref{fig:fb-zero-cases-en}.
Irrelevant arcs are omitted, and configurations that give rise to the same
obstruction are grouped together.

\begin{figure}[!t]
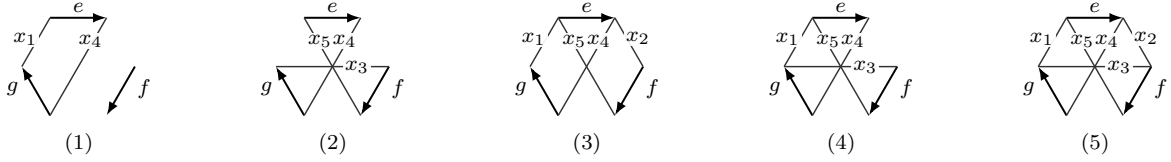

\centering
\Conf{\FBZone}{1} \Conf{\FBZtwo}{2} \Conf{\FBZthree}{3} \Conf{\FBZfour}{4}
\Conf{\FBZfive}{5}
\caption{The five local obstruction types for
$\E(I_\alpha\mid I_s=0,w_s=0)$.}
\label{fig:fb-zero-cases-en}
\end{figure}

\noindent \textbf{Case (1).} This structure forms only one $4$-cycle.
Put $Z_1=\min\{e,g,x_1,x_4\},$ and let $C_1=\{Z_1\in\{e,g\}\}.$ Conditioned on
the existence of the arcs involved in this configuration, the probability that
$\alpha$ generates a $2$-face is $\Prob(C_1)=\frac12.$

\noindent \textbf{Case (2).} This structure forms only one $6$-cycle.
Put $Z_2=\min\{e,f,g,x_3,x_4,x_5\},$ and let $C_2=\{Z_2\in\{e,f,g\}\}.$
Conditioned on the existence of the arcs involved in this configuration, the
probability that $\alpha$ generates a $2$-face is $\Prob(C_2)=\frac12.$

\noindent \textbf{Case (3).} This structure forms two $4$-cycles.
Put $Z_3=\min\{e,f,g,x_1,x_2,x_4,x_5\},$ and let
$C_3=\{\min\{e,g,x_1,x_4\}\in\{e,g\},\min\{e,f,x_2,x_5\}\in\{e,f\}\}.$
Conditioned on the existence of the arcs involved in this configuration, the
probability that $\alpha$ generates a $2$-face is
\begin{align*}
  \Prob(C_3) &=\Prob(Z_3=e)\Prob(C_3\mid Z_3=e) +\Prob(Z_3=f)\Prob(C_3\mid Z_3=f)\\
  &\quad+\Prob(Z_3=g)\Prob(C_3\mid Z_3=g)\\
  &=\frac17\cdot1+\frac17\cdot\frac12+\frac17\cdot\frac12 =\frac27.
\end{align*}

\noindent \textbf{Case (4).} This structure forms one $4$-cycle and one
$6$-cycle.
Put $Z_4=\min\{e,\allowbreak f,\allowbreak g,\allowbreak x_1,\allowbreak
x_3,\allowbreak x_4,\allowbreak x_5\},$ and let $C_4=\{\min\{e,\allowbreak
g,\allowbreak x_1,\allowbreak x_4\}\in\nobreak\{e,g\},\allowbreak
\min\{e,\allowbreak f,\allowbreak g,\allowbreak x_3,\allowbreak x_4,\allowbreak
x_5\}\in\nobreak\{e,f,g\}\}.$ Conditioned on the existence of the arcs involved
in this configuration, the probability that $\alpha$ generates a $2$-face is
\begin{align*}
  \Prob(C_4) &=\Prob(Z_4=e)\Prob(C_4\mid Z_4=e) +\Prob(Z_4=f)\Prob(C_4\mid Z_4=f)\\
  &\quad+\Prob(Z_4=g)\Prob(C_4\mid Z_4=g)\\
  &=\frac17\cdot1+\frac17\cdot\frac12+\frac17\cdot1 =\frac5{14}.
\end{align*}

\begin{samepage}
\noindent \textbf{Case (5).} This structure forms two $4$-cycles and one
$6$-cycle.
Put $Z_5=\min\{e,\allowbreak f,\allowbreak g,\allowbreak x_1,\allowbreak
x_2,\allowbreak x_3,\allowbreak x_4,\allowbreak x_5\},$ and let
$C_5=\{\min\{e,\allowbreak g,\allowbreak x_1,\allowbreak
x_4\}\in\nobreak\{e,g\},\allowbreak \min\{e,\allowbreak f,\allowbreak
x_2,\allowbreak x_5\}\in\nobreak\{e,f\},\allowbreak \min\{e,\allowbreak
f,\allowbreak g,\allowbreak x_3,\allowbreak x_4,\allowbreak
x_5\}\in\nobreak\{e,f,g\}\}.$ Conditioned on the existence of the arcs involved
in this configuration, the probability that $\alpha$ generates a $2$-face is
\begin{align*}
  \Prob(C_5) &=\Prob(Z_5=e)\Prob(C_5\mid Z_5=e) +\Prob(Z_5=f)\Prob(C_5\mid Z_5=f)\\
  &\quad+\Prob(Z_5=g)\Prob(C_5\mid Z_5=g)\\
  &=\frac18\cdot1+\frac18\cdot\frac12+\frac18\cdot\frac12 =\frac14.
\end{align*}
\end{samepage}

Let $K$ be the number of arcs in $\bar\alpha$ that are present in $G_{n,p}$, and
classify according to the value of $K$ as follows.
\begin{enumerate}[label=(\roman*)]
    \item If $K=0$ or $K=1$, there exists no structures listed above.

    \item If $K=2$, two choices correspond to Case~$(1)$, whereas the
    remaining eight do not trigger an obstruction.

    \item If $K=3$, six choices correspond to Case~$(1)$, one to Case~$(2)$,
    and the remaining three do not trigger an obstruction.

    \item If $K=4$, two choices correspond to Case~$(1)$, one to Case~$(3)$,
    and the remaining two to Case~$(4)$.

    \item If $K=5$, the unique choice corresponds to Case~$(5)$.
\end{enumerate}

\noindent Summing over these cases gives
\begin{align}
  \E(I_\alpha\mid I_s=0,w_s=0)
  &=p^3\Bigl(q^5+5pq^4+(2\Prob(C_1)+8)p^2q^3+(6\Prob(C_1)+\Prob(C_2)+3)p^3q^2\notag\\
  &\qquad +(2\Prob(C_1)+\Prob(C_3)+2\Prob(C_4))p^4q+\Prob(C_5)p^5\Bigr)\notag\\
  &=p^3\left(q^5+5pq^4+9p^2q^3+\frac{13}{2}p^3q^2+2p^4q+\frac14p^5\right)\notag\\
  &=p^3q^5+5p^4q^4+9p^5q^3+\frac{13}{2}p^6q^2+2p^7q+\frac14p^8.
  \label{eq:fb2}
\end{align}

\noindent Combining (\ref{eq:fa1}), (\ref{eq:fb1}), and (\ref{eq:fb2}), we
obtain
\begin{align*}
  f_a(0)+2f_b(0)
  &=p^2q^6+6p^3q^5+\frac{29}{2}p^4q^4+18p^5q^3+12p^6q^2+4p^7q+\frac12p^8\\
  &\quad+2\bigl[p^3q^5+4p^4q^4+4p^5q^3+p^6q^2\bigr]\\
  &\quad-2\bigl[p^3q^5+5p^4q^4+9p^5q^3+\frac{13}{2}p^6q^2+2p^7q+\frac14p^8\bigr]\\
  &=p^2q^6+6p^3q^5+\frac{25}{2}p^4q^4+8p^5q^3+p^6q^2\\
  &=p^2q^2\left(q^4+6pq^3+\frac{25}{2}p^2q^2+8p^3q+p^4\right)\\
  &=\frac12p^2q^2\bigl(p^4-6p^3+p^2+4p+2\bigr).
\end{align*}
Substituting this into (\ref{eq:two-node-overlap-contribution}), we obtain
\[
  \frac12 (n)_{10}pq\int_0^1 \bigl(f_a(t)+2f_b(t)\bigr)^2\,dt=\frac18 (n)_{10}p^5q^5
  \bigl(p^4-6p^3+p^2+4p+2\bigr)^2 .
\]

\medskip
\noindent
\textbf{\boldmath$|V(\alpha)\cap V(\beta)|\geq 3$:} It remains to show that the
pairs in the final overlap regime are negligible.
Set $k=|V(\alpha)\cap V(\beta)|$, so that $3\le k\le 6$.
The union of the two configurations involves $12-k$ nodes.
Among the six underlying arcs of the two triples, at most $\lfloor k/2\rfloor$
arcs can be identified through the overlap.
Hence the two triples use at least $6-\lfloor k/2\rfloor$ distinct underlying
arcs, and the contribution from this overlap pattern is at most of order
$n^{12-k}p^{\,6-\lfloor k/2\rfloor}$.
For $k=3,4,5,6$, the ratios of these bounds to $n^{10}p^5$ are
\[
  \frac1n,\quad\frac1{n^2p},\quad\frac1{n^3p}=\frac1n\cdot\frac1{n^2p},
  \quad\frac1{n^4p^2}=\left(\frac1{n^2p}\right)^2.
\]
Since $p\gg n^{-2}$, all these ratios tend to zero as $n\to\infty$.
Consequently, pairs sharing at least three nodes contribute only $o(n^{10}p^5)$.

Having treated all three overlap regimes, we now combine the estimates.
The contributions with $|V(\alpha)\cap V(\beta)|\le 1$ vanish, while the
contribution with $|V(\alpha)\cap V(\beta)|=2$ is
\[
  \frac18 (n)_{10}p^5q^5 \bigl(p^4-6p^3+p^2+4p+2\bigr)^2,
\]
and the total contribution from the pairs with $|V(\alpha)\cap V(\beta)|\ge 3$
is $o(n^{10}p^5)$ as $n\to\infty$.
Therefore,
\[
  \E\bigl[\Var(H_{n,p}^{(1)}\mid w)\bigr]=\frac18 (n)_{10}p^5q^5
  \bigl(p^4-6p^3+p^2+4p+2\bigr)^2+o(n^{10}p^5),\quad\text{as }n\to\infty.
\]
This completes the proof.
\end{proof}
\begin{remark}
We would like to highlight two key ingredients in the proof of
Proposition~\ref{prop:variance-H1}.
First, among all possible overlap patterns of two triples in $\mathcal M_1$,
only the case $|V(\alpha)\cap V(\beta)|=2$ contributes on the scale $n^{10}p^5$.
The cases $|V(\alpha)\cap V(\beta)|\le 1$ give zero covariance, whereas the
cases $|V(\alpha)\cap V(\beta)|\ge 3$ are of smaller order.
Second, in the remaining case, the dependence between $I_\alpha$ and $I_\beta$
is transmitted through the connecting arc $s$.
Conditioning on the weight of this arc reduces the covariance calculation to the
integral $\int_0^1 \bigl(f_a(t)+2f_b(t)\bigr)^2\,dt$.
By Lemma~\ref{lem:fa-fb-constant}, the combination $f_a(t)+2f_b(t)$ is independent of
$t$, and hence it is enough to compute $f_a(0)+2f_b(0)$.
This avoids an enumeration over all possible positions of the connecting arc in
the order, and leaves only a finite local computation.
\end{remark}
The preceding proposition gives the contribution of the leading term
$H_{n,p}^{(1)}$ on the scale $n^{10}p^5$.
It remains to show that all remaining terms in the decomposition are negligible
on this scale.
\begin{proposition}\label{prop:variance-remainder}
If $p\gg n^{-2}$ and $\limsup\limits_{n\to\infty}p<1$, then, as $n\to\infty$,
for $i=2,3,4$,
\[
  \E\bigl[\Var(H_{n,p}^{(i)}\mid w)\bigr]=o(n^{10}p^5),
\]
while the edge and arc counts satisfy
\[
  \E\bigl[\Var(K_{n,p}\mid w)\bigr]=o(n^{10}p^5),
\]
and
\[
  \E\bigl[\Var(E_{n,p}\mid w)\bigr]=o(n^{10}p^5).
\]
\end{proposition}
\begin{proof}
We estimate each term separately.
Throughout the proof, all asymptotic notation is with respect to $n\to\infty$,
and all implicit estimates hold for sufficiently large $n$.
Conditioning on the order $w$, we expand the conditional variances into
covariance sums and group the summands according to the number of common nodes
between the two triples.
Each resulting contribution is then bounded by counting the nodes in the joint
configuration and the distinct underlying arcs required for the corresponding
indicators.

\medskip
\noindent\textbf{The term $H_{n,p}^{(2)}$.}
Expanding the conditional variance gives
\[
  \E\bigl[\Var(H_{n,p}^{(2)}\mid w)\bigr]
  =\sum_{\alpha,\beta\in\mathcal M_2}\E\bigl[\Cov(I_\alpha,I_\beta\mid w)\bigr].
\]
We classify the contributions according to the number of common nodes
$|V(\alpha)\cap V(\beta)|$.
\begin{enumerate}[label=(\roman*)]
\item If $|V(\alpha)\cap V(\beta)|=0$, by the definitions introduced in the proof of Lemma~\ref{lem:H-expectation},
we decompose $\mathcal M_2$ as
$\mathcal M_2=\mathcal M_2^{(1)}\sqcup \mathcal M_2^{(2)}, $ where
$\mathcal M_2^{(1)}$ corresponds to the same-orientation case at the common
node, and $\mathcal M_2^{(2)}$ corresponds to the opposite-orientation case.
If at least one of $\alpha$ and $\beta$ belongs to $\mathcal M_2^{(1)}$, then
the corresponding indicators are independent after conditioning on $w$.
Hence $\E\bigl[\Cov(I_\alpha,I_\beta\mid w)\bigr]=0$ unless
$\alpha,\beta\in\mathcal M_2^{(2)}$.

It remains to consider the case $\alpha,\beta\in\mathcal M_2^{(2)}$.
For such a pair, we have
\[
  \begin{aligned}
    \left|\E\bigl[\Cov(I_\alpha,I_\beta\mid w)\bigr]\right|
    &=\left|\E\Bigl[\E(I_\alpha I_\beta\mid w)-\E(I_\alpha\mid w)\E(I_\beta\mid
    w)\Bigr]\right|\\
    &\le\max\left\{\E(I_\alpha I_\beta),\E\Bigl[\E(I_\alpha\mid w)\E(I_\beta\mid
    w)\Bigr]\right\}.
  \end{aligned}
\]
Let $B_\alpha$ denote the event that no directed $4$-cycle containing the two
intersecting directed arcs in $\alpha$ triggers the obstruction.
By the estimate used in the proof of Lemma~\ref{lem:H-expectation}, we have
\[
  \mathbb P(B_\alpha)\le\min\left(1,\frac{2}{np^2}\right).
\]
Since $|V(\alpha)\cap V(\beta)|=0$, the six underlying arcs appearing in
$\alpha$ and $\beta$ are distinct.
Hence
\[
  \E(I_\alpha I_\beta)\le p^6\mathbb P(B_\alpha).
\]
The corresponding conditional expectations satisfy
\[
  \E(I_\alpha\mid w)\le p^3\mathbb P(B_\alpha\mid w),
  \qquad\E(I_\beta\mid w)\le p^3.
\]
Thus, taking expectations yields
\[
  \E\Bigl[\E(I_\alpha\mid w)\E(I_\beta\mid w)\Bigr]
  \le p^6\E\bigl[\mathbb P(B_\alpha\mid w)\bigr]= p^6\mathbb P(B_\alpha).
\]
Combining these estimates gives
\begin{equation}\label{eq:covariance-bound}
  \left|\E\bigl[\Cov(I_\alpha,I_\beta\mid w)\bigr]\right|
  \le p^6\mathbb P(B_\alpha)\le p^6\min\left(1,\frac{2}{np^2}\right).
\end{equation}
We use the elementary bound
\begin{equation}\label{eq:elementary-bound}
  p\min\left(1,\frac{2}{np^2}\right) \le \sqrt{\frac{2}{n}} .
\end{equation}
Indeed, if $p<\sqrt{2/n}$, then
\[
  p\min\left(1,\frac{2}{np^2}\right) \le p<\sqrt{\frac{2}{n}},
\]
whereas if $p\ge \sqrt{2/n}$, then
\[
  p\min\left(1,\frac{2}{np^2}\right) \le\frac{2}{np} \le\sqrt{\frac{2}{n}}.
\]
Combining (\ref{eq:covariance-bound}) and (\ref{eq:elementary-bound}), we obtain
\[
  \left|\E\bigl[\Cov(I_\alpha,I_\beta\mid w)\bigr]\right|
  \le p^6\min\left(1,\frac{2}{np^2}\right)\le\sqrt{\frac{2}{n}}\,p^5.
\]
Since $|V(\alpha)\cap V(\beta)|=0$, the joint configuration uses $10$ nodes.
Hence the contribution from this case is at most on the order of $n^{19/2}p^5$.

\item[(ii)] If $|V(\alpha)\cap V(\beta)|=1$, the joint configuration uses
$9$ nodes and at least $6$ arcs.
The contribution is at most of order $n^9p^6$.

\item[(iii)] If $|V(\alpha)\cap V(\beta)|=2$, the joint configuration uses
$8$ nodes and at least $5$ arcs.
The contribution is at most of order $n^8p^5$.

\item[(iv)] If $|V(\alpha)\cap V(\beta)|=3$, the joint configuration uses
$7$ nodes and at least $4$ arcs.
The contribution is at most of order $n^7p^4$.

\item[(v)] If $|V(\alpha)\cap V(\beta)|=4$, the joint configuration uses
$6$ nodes and at least $4$ arcs.
The contribution is at most of order $n^6p^4$.

\item[(vi)] If $|V(\alpha)\cap V(\beta)|=5$, the joint configuration uses
$5$ nodes and at least $3$ arcs.
The contribution is at most of order $n^5p^3$.
\end{enumerate}

Dividing these six bounds by $n^{10}p^5$ gives the ratios
\[
  \frac{1}{\sqrt{n}},\qquad \frac{p}{n},\qquad \frac1{n^2},\qquad \frac1{n^3p}
  = \frac1n\cdot\frac1{n^2p},\qquad \frac1{n^4p}
  = \frac1{n^2}\cdot\frac1{n^2p},\qquad \frac1{n^5p^2}
  = \frac1n\left(\frac1{n^2p}\right)^2 .
\]
Since $p\in(0,1)$ and $p\gg n^{-2}$, all these ratios tend to zero.
Thus all overlap contributions in $\mathcal M_2$ are dominated by $n^{10}p^5$.
It follows that
\[
  \E\bigl[\Var(H_{n,p}^{(2)}\mid w)\bigr]=o(n^{10}p^5).
\]
\medskip
\noindent\textbf{The terms $H_{n,p}^{(3)}$ and $H_{n,p}^{(4)}$.}
We first estimate $H_{n,p}^{(3)}$.
Expanding the conditional variance gives
\[
  \E\bigl[\Var(H_{n,p}^{(3)}\mid w)\bigr]
  =\sum_{\alpha,\beta\in\mathcal M_3}\E\bigl[\Cov(I_\alpha,I_\beta\mid w)\bigr].
\]
We classify the contributions according to the number of common nodes
$|V(\alpha)\cap V(\beta)|$.
\begin{enumerate}[label=(\roman*)]
\item If $|V(\alpha)\cap V(\beta)|=0$, the joint configuration uses $8$ nodes
and at least $6$ arcs.
The contribution is at most of order $n^8p^6$.
\item If $|V(\alpha)\cap V(\beta)|=1$, the joint configuration uses $7$ nodes
and at least $6$ arcs.
The contribution is at most of order $n^7p^6$.
\item If $|V(\alpha)\cap V(\beta)|=2$, the joint configuration uses $6$ nodes
and at least $5$ arcs.
The contribution is at most of order $n^6p^5$.
\item If $|V(\alpha)\cap V(\beta)|=3$, the joint configuration uses $5$ nodes
and at least $4$ arcs.
The contribution is at most of order $n^5p^4$.
\item If $|V(\alpha)\cap V(\beta)|=4$, the joint configuration uses $4$ nodes
and at least $3$ arcs.
The contribution is at most of order $n^4p^3$.
\end{enumerate}

Dividing these five bounds by $n^{10}p^5$ gives the ratios
\[
  \frac{p}{n^2},\qquad\frac{p}{n^3},\qquad\frac1{n^4},\qquad\frac1{n^5p}
  =\frac1{n^3}\cdot\frac1{n^2p},\qquad\frac1{n^6p^2}=\frac1{n^2}\left(\frac1{n^2p}\right)^2.
\]
Since $p\in(0,1)$ and $p\gg n^{-2}$, all these ratios tend to zero.
Hence
\[
  \E\bigl[\Var(H_{n,p}^{(3)}\mid w)\bigr]=o(n^{10}p^5).
\]
The same argument applies to $H_{n,p}^{(4)}$, since triples in $\mathcal M_4$
are also supported on four nodes.
Therefore
\[
  \E\bigl[\Var(H_{n,p}^{(4)}\mid w)\bigr]=o(n^{10}p^5).
\]

\medskip
\noindent\textbf{The terms $K_{n,p}$ and $E_{n,p}$.}
Since $K_{n,p}$ does not depend on the random order $w$, it follows from
\cite[Theorem~1.1(a)]{DonzelmannJuhnkeRednossThaele2026} that
\[
  \E\bigl[\Var(K_{n,p}\mid w)\bigr]=\Var(K_{n,p}) \asymp n^6p^3.
\]
Moreover, since $p\gg n^{-2}$, we have
\[
  \frac{n^6p^3}{n^{10}p^5}=\left(\frac1{n^2p}\right)^2 \longrightarrow 0.
\]
The variance contribution of $K_{n,p}$ is therefore negligible:
\[
  \E\bigl[\Var(K_{n,p}\mid w)\bigr] =o(n^{10}p^5).
\]
Similarly, $E_{n,p}\sim \operatorname{Bin}\left(\binom n2,p\right)$, and $E_{n,p}$ is
independent of $w$.
Thus
\[
  \E\bigl[\Var(E_{n,p}\mid w)\bigr]=\Var(E_{n,p})=\binom n2 p(1-p)\asymp n^2p.
\]
Again, since $p\gg n^{-2}$, we have
\[
  \frac{n^2p}{n^{10}p^5}=\left(\frac1{n^2p}\right)^4 \longrightarrow 0.
\]
It follows that
\[
  \E\bigl[\Var(E_{n,p}\mid w)\bigr]=o(n^{10}p^5).
\]

The preceding estimates show that each of the remaining terms is negligible on
the scale $n^{10}p^5$.
This establishes all the assertions of the proposition.
\end{proof}
\begin{proof}[Proof of Theorem~\ref{thm:variance}]
Recall (\ref{eq:total-variance}) and (\ref{eq:F-variance-decomposition}).
We compare the terms in (\ref{eq:variance-decomposition}).
By Proposition~\ref{prop:variance-H1},
\[
  \E\bigl[\Var(H_{n,p}^{(1)}\mid w)\bigr]
  = \frac18 (n)_{10}p^5q^5 \bigl(p^4-6p^3+p^2+4p+2\bigr)^2 +o(n^{10}p^5).
\]
Also, Proposition~\ref{prop:variance-remainder} gives
\[
  \E\bigl[\Var(H_{n,p}^{(i)}\mid w)\bigr] = o(n^{10}p^5), \qquad i=2,3,4,
\]
and
\[
  \E\bigl[\Var(K_{n,p}\mid w)\bigr]
  = o(n^{10}p^5),
  \qquad \E\bigl[\Var(E_{n,p}\mid w)\bigr]
  = o(n^{10}p^5).
\]

It remains only to check that the covariance terms in
(\ref{eq:variance-decomposition}) do not contribute to the leading order.
Note that a simple algebraic argument gives, for any $X$, $Y$, and $w$,
\[
  \left| \E\bigl[\Cov(X,Y\mid w)\bigr] \right|
  \le \Bigl(\E\bigl[\Var(X\mid w)\bigr]\Bigr)^{1/2} \Bigl(\E\bigl[\Var(Y\mid
  w)\bigr]\Bigr)^{1/2}.
\]
It is therefore easy to see that each covariance term involving at least one of
\[
  H_{n,p}^{(2)},\quad H_{n,p}^{(3)},\quad H_{n,p}^{(4)}, \quad K_{n,p},\quad E_{n,p}
\]
is $o(n^{10}p^5)$.
Since there are only finitely many such terms, their total contribution is also
$o(n^{10}p^5)$.
Therefore, we have
\begin{align*}
  \Var(F_{n,p}) &=\E\bigl[\Var(H_{n,p}^{(1)}\mid w)\bigr] +o(n^{10}p^5)\\
  &=\frac18 (n)_{10}p^5q^5\bigl(p^4-6p^3+p^2+4p+2\bigr)^2 +o(n^{10}p^5).
\end{align*}
This is the asserted variance asymptotic.
Finally, note that
\[
  p^4-6p^3+p^2+4p+2=2+p(1-p)(4+5p-p^2)>2, \qquad 0<p<1.
\]
This completes the proof.
\end{proof}
\begin{remark}
The variance scale in Theorem~\ref{thm:variance} should be compared with the edge-count
result for unimodular triangulations of random symmetric edge polytopes.
For the edge count $K_{n,p}$, the variance is of order $n^6p^3$; see
\cite[Theorem~1.1(a)]{DonzelmannJuhnkeRednossThaele2026}.
In contrast, the $2$-face statistic $F_{n,p}$ has variance of order $n^{10}p^5$.
This is consistent with the variance scale $n^{4k+2}p^{2k+1}$ for a fixed $k$.
The leading covariance contribution comes from pairs of configurations whose
arcs are pairwise node-disjoint within each configuration and whose node sets
intersect in exactly two nodes.
Dependence is carried by the arc joining these common nodes, which need not
belong to either configuration.
Section~\ref{sec:higher-dimensional} uses this mechanism to establish the
variance order in Theorem~\ref{thm:k-face-clt}.
\end{remark}

\subsection{Discrete gradients}\label{sec:discrete-gradients}
After computing the expectation and variance of $F_{n,p}$, the next step is to
control the first- and second-order discrete gradients of $F_{n,p}$ and estimate
their fourth moments respectively.
According to the preceding definitions and assumptions, $G_{n,p}$ is an
Erd\H{o}s--R\'{e}nyi random graph, $E_{n,p}$ denotes the number of arcs of
$G_{n,p}$, and $K_{n,p}$ denotes the number of edges in any unimodular
triangulation of $P_{n,p}$, the random symmetric edge polytope generated by
$G_{n,p}$.
We introduce some additional notation below.

Fix an arc $e=\{s,t\}\in\binom{[n]}{2}.$ Let $G_{n,p}+e$ (respectively,
$G_{n,p}-e$) denote the graph obtained from $G_{n,p}$ by adding (respectively,
deleting) the arc $e$.
Let $(E_{n,p})^-_e$ denote the number of arcs of $G_{n,p}-e$.

If $s=v_0,v_1,v_2,\ldots,v_i=t$ is a sequence of distinct nodes of $G_{n,p}$,
and if for every $l=0,1,\ldots,i-1$ the pair $\{v_l,v_{l+1}\}$ is an arc of
$G_{n,p}$, then the arcs in this sequence are said to form an $s$-$t$ path of
length $i$ in $G_{n,p}$.
The length $i$ refers to the number of arcs traversed.
The number of $s$-$t$ paths of length $i$ in $G_{n,p}$ is denoted by
$W_{n,p}(i)$.

\begin{lemma}\label{lem:first-order-gradient-F}
Fix an arc $e=\{s,t\}\in\binom{[n]}{2}$.
Then the following inequality holds almost surely:
\begin{align}
  |D_eF_{n,p}|
  &\le \sqrt{pq}\Bigl( 4+2W_{n,p}(2)+14W_{n,p}(3)+8W_{n,p}(4)+20W_{n,p}(5) \notag\\
  &\hspace{3.2cm} +4E_{n,p}+4E_{n,p}^2+4W_{n,p}(2)E_{n,p}+12W_{n,p}(3)E_{n,p} \Bigr).
  \label{eq:first-order-gradient-bound}
\end{align}
\end{lemma}

\begin{proof}
By Lemma~\ref{lem:two-face-decomposition}, we have
\[
  D_eF_{n,p}=D_eH_{n,p}+D_eK_{n,p}-2D_eE_{n,p},
\]
and therefore
\[
  |D_eF_{n,p}| \le |D_eH_{n,p}|+|D_eK_{n,p}|+2|D_eE_{n,p}|.
\]
Flipping the state of $e$ changes $E_{n,p}$ by exactly one, so the gradient
definition gives
\begin{equation}\label{eq:first-gradient-Enp}
  |D_eE_{n,p}|=\sqrt{pq}.
\end{equation}
Also, by \cite[Lemma~4.4]{DonzelmannJuhnkeRednossThaele2026}, we have
\begin{equation}\label{eq:first-gradient-Knp}
  |D_eK_{n,p}| \le \sqrt{pq}\bigl(2+4E_{n,p}+2W_{n,p}(2)+6W_{n,p}(3)\bigr).
\end{equation}
Thus it remains only to estimate $D_eH_{n,p}$.
Note that
\[
  D_eH_{n,p}=\sqrt{pq}\bigl((H_{n,p})_e^+-(H_{n,p})_e^-\bigr),
\]
and that the absolute difference $|(H_{n,p})_e^+-(H_{n,p})_e^-|$ is bounded by
the number of non-origin $2$-faces created by adding $e$, plus the number of
non-origin $2$-faces destroyed by adding $e$.
We treat these two numbers separately.

\textbf{Created $2$-faces.} A newly created $2$-face must involve $e$ and must
also involve two arcs of $G_{n,p}-e$.
After these arcs are chosen, each of the three arcs may be oriented in two
possible directions.
Therefore, the number of newly created $2$-faces is at most
\[
  2^3\binom{(E_{n,p})_e^-}{2} \le 8\frac{E_{n,p}(E_{n,p}-1)}{2} \le 4E_{n,p}^2.
\]

\noindent\textbf{Destroyed $2$-faces.}
Suppose that adding the arc $e=\{s,t\}$ destroys a $2$-face.
By Proposition~\ref{prop:two-face-criterion}, if this $2$-face is generated by the directed
arcs $\{\dir k,\dir l,\dir m\}$, then at least two of these arcs lie on an
obstruction cycle containing $e$.
Thus, for an upper bound, it suffices to count the triples
$\{\dir k,\dir l,\dir m\}$ for which some of their arcs lie on an $s$--$t$ path
that closes with $e$ to a directed cycle.
More precisely, either two arcs in the triple lie on an $s$--$t$ path of length
$2$ or $3$, or all three lie on an $s$--$t$ path of length $3$, $4$, or $5$.
We now count these possibilities according to the path length.

If the path has length $2$, then the two arcs forming the path must have
opposite directions with respect to the middle node.
Thus there are two possible orientations.
The third directed arc in the triple may be chosen arbitrarily, with either
orientation.
Consequently, the number of $2$-faces destroyed by all length-$2$ paths is at
most $4W_{n,p}(2)E_{n,p}$.

If the path has length $3$, and if $\{\dir k,\dir l,\dir m\}$ itself forms a
$3$-path, then no matter what the individual directions are, together with $e$
it forms a directed $4$-cycle.
Hence there are $8$ possible orientations.

If two arcs among $\{\dir k,\dir l,\dir m\}$ lie on a $3$-path, then there are
$\binom{3}{2}$ choices for the two arcs.
For these two arcs to lie on an obstruction cycle of length $4$ containing $e$,
there are two possible orientations.
The remaining directed arc in the triple may be chosen arbitrarily, with either
orientation.
Thus the number of $2$-faces destroyed by all length-$3$ paths is at most
$8W_{n,p}(3)+12W_{n,p}(3)E_{n,p}$.

If the path has length $4$, then the three arcs $\{\dir k,\dir l,\dir m\}$ all
lie on this path.
There are $\binom{4}{3}$ choices for these arcs.
They lie on an obstruction cycle of length $5$ containing $e$, and there are two
possible orientations.
Therefore, the number of $2$-faces destroyed by all length-$4$ paths is at most
$8W_{n,p}(4)$.

If the path has length $5$, then the three arcs $\{\dir k,\dir l,\dir m\}$ all
lie on this path.
There are $\binom{5}{3}$ choices for these arcs.
They lie on an obstruction cycle of length $6$ containing $e$, and there are two
possible orientations.
Therefore, the number of $2$-faces destroyed by all length-$5$ paths is at most
$20W_{n,p}(5)$.

Combining these estimates together easily yields
\begin{align}
  |(H_{n,p})_e^+-(H_{n,p})_e^-|
  &\le 8W_{n,p}(3)+8W_{n,p}(4)+20W_{n,p}(5)+4E_{n,p}^2 \notag\\
  &\quad +4W_{n,p}(2)E_{n,p}+12W_{n,p}(3)E_{n,p}.
\end{align}
It follows that
\begin{align}
  |D_eH_{n,p}|
  &\le \sqrt{pq}\Bigl( 8W_{n,p}(3)+8W_{n,p}(4)+20W_{n,p}(5)+4E_{n,p}^2 \notag\\
  &\hspace{1cm} +4W_{n,p}(2)E_{n,p}+12W_{n,p}(3)E_{n,p} \Bigr).
  \label{eq:first-gradient-Hnp}
\end{align}
Combining (\ref{eq:first-gradient-Enp}), (\ref{eq:first-gradient-Knp}), and
(\ref{eq:first-gradient-Hnp}) proves (\ref{eq:first-order-gradient-bound}).
\end{proof}

\begin{lemma}\label{lem:fourth-moment-DeFnp}
Let $e=\{s,t\}\in \binom{[n]}{2}$.
If $p\gg n^{-2}$, then
\[
  \mathbb{E}[(D_eF_{n,p})^4]\lesssim n^{16}p^{10}q^2 .
\]
\end{lemma}

\begin{proof}
Recall that $E_{n,p}$ denotes the number of arcs of $G_{n,p}$.
The binomial moment estimates give
\[
  E_{n,p}\sim \operatorname{Bin}\left(\binom{n}{2},p\right), \quad
  \mathbb{E}[E_{n,p}^4]\asymp n^8p^4, \quad \mathbb{E}[E_{n,p}^8]\asymp n^{16}p^8 .
\]

We estimate the moments of $W_{n,p}(i)$ for $i=2,3,4,5$.
Let $S_i$ denote the set of all simple $s$--$t$ paths of length $i$ in the
complete graph on $[n]$.
Then it follows that
\[
  W_{n,p}(i)=\sum_{\gamma\in S_i}\one_{\{\gamma\in G_{n,p}\}},
\]
where $\gamma\in G_{n,p}$ means that every arc of the path $\gamma$ exists in
$G_{n,p}$.

First consider $W_{n,p}(2)$.
To form an $s$-$t$ path of length $2$, one needs to choose one node $u$ from the
node set of $G_{n,p}$ outside $s$ and $t$, and the arcs $su$ and $tu$ must both
be present in $G_{n,p}$.
Therefore $S_2$ contains $n-2$ elements.
Moreover, the random variables $\one_{\{\gamma_1\in G_{n,p}\}}$ and
$\one_{\{\gamma_2\in G_{n,p}\}}$ are independent when $\gamma_1\ne \gamma_2$,
since the involved arcs are distinct.
Therefore
\[
  W_{n,p}(2)\sim \operatorname{Bin}(n-2,p^2), \quad \mathbb{E}[W_{n,p}(2)^4]\asymp
  n^4p^8+np^2, \quad \mathbb{E}[W_{n,p}(2)^8]\asymp n^8p^{16}+np^2 .
\]
As for $W_{n,p}(3)$, we have
\[
  \begin{aligned}
    \mathbb{E}[W_{n,p}(3)^4]
    &= \sum_{\gamma_1\in S_3} \sum_{\gamma_2\in S_3} \sum_{\gamma_3\in S_3}
    \sum_{\gamma_4\in S_3} \mathbb{E}[ \one_{\{\gamma_1\in G_{n,p}\}}
    \one_{\{\gamma_2\in G_{n,p}\}} \one_{\{\gamma_3\in G_{n,p}\}} \one_{\{\gamma_4\in
    G_{n,p}\}} ] .
  \end{aligned}
\]
To bound this sum, we classify the contributions according to the overlap
patterns of $\gamma_1,\gamma_2,\gamma_3,\gamma_4$.
If exactly $k$ nodes other than $s$ and $t$ are involved, then there must be at
least $k+\lceil k/2\rceil$ arcs.
Therefore the total contribution of this type is at most of order
$n^k p^{k+\lceil k/2\rceil}$.
Thus we have the upper bound
\begin{align*}
  \mathbb{E}\bigl[W_{n,p}(3)^4\bigr] &\lesssim n^8p^{12}+n^7p^{11}+n^6p^9+n^5p^8 \\
  &\quad +n^4p^6+n^3p^5+n^2p^3 \\
  &\lesssim n^8p^{12}+n^2p^3.
\end{align*}

On the other hand, it is easy to see that
\[
  \mathbb{E}[W_{n,p}(3)^4]\gtrsim n^8p^{12}+n^2p^3 .
\]
Combining the upper and lower bounds gives
\[
  \mathbb{E}[W_{n,p}(3)^4]\asymp n^8p^{12}+n^2p^3 .
\]
Similarly, we have
\[
  \mathbb{E}[W_{n,p}(3)^8]\asymp n^{16}p^{24}+n^2p^3,
\]
and
\[
  \mathbb{E}[W_{n,p}(4)^4]\asymp n^{12}p^{16}+n^3p^4,
  \qquad \mathbb{E}[W_{n,p}(5)^4]\asymp n^{16}p^{20}+n^4p^5 .
\]
Finally, applying the inequality
\[
  \left(\sum_{i=1}^9 x_i\right)^4 \leq 9^3\sum_{i=1}^9 x_i^4
\]
and the Cauchy--Schwarz inequality, we obtain
\begin{align*}
  \E[(D_eF_{n,p})^4]
  &\le Cp^2q^2\Bigl\{\E[E_{n,p}^8] +\E[W_{n,p}(4)^4]+\E[W_{n,p}(5)^4]\\
  &\qquad\qquad+\E[W_{n,p}(2)^4E_{n,p}^4] +\E[W_{n,p}(3)^4E_{n,p}^4]\Bigr\}\\
  &\le Cp^2q^2\Bigl\{\E[E_{n,p}^8] +\E[W_{n,p}(4)^4]+\E[W_{n,p}(5)^4]\\
  &\qquad\qquad+\E[W_{n,p}(2)^8]^{1/2}\E[E_{n,p}^8]^{1/2}\\
  &\qquad\qquad+\E[W_{n,p}(3)^8]^{1/2}\E[E_{n,p}^8]^{1/2}\Bigr\}\\
  &\lesssim n^{16}p^{10}q^2,
\end{align*}
which completes the proof.
\end{proof}

\begin{lemma}[Second-order gradient bound for $F_{n,p}$]
\label{lem:second-order-gradient-F}
Let $e,f\in\binom{[n]}2$ be distinct arcs.
Under the assumptions of Theorem~\ref{thm:variance}, we have
\begin{equation}\label{eq:second-gradient-F-fourth-moment}
  \E[(D_eD_fF_{n,p})^4] \lesssim
  \begin{cases}
    n^8p^8q^4,& |e\cap f|=0,\\
    n^8p^8q^4+n^{12}p^{16}q^4,& |e\cap f|=1.
  \end{cases}
\end{equation}
\end{lemma}

\begin{proof}
Recall from Lemma~\ref{lem:two-face-decomposition} that
$F_{n,p}=H_{n,p}+K_{n,p}-2E_{n,p}.$ Since $E_{n,p}$ is linear in the arc
indicators, $D_eD_fE_{n,p}=0$ for $e\ne f$.
Also, by \cite[Lemma~4.5]{DonzelmannJuhnkeRednossThaele2026},
\begin{equation}\label{eq:l313-K-fourth-moment}
  \E[(D_eD_fK_{n,p})^4]\lesssim
  \begin{cases}
    p^4q^4,&|e\cap f|=0,\\
    p^4q^4+n^4p^{12}q^4,&|e\cap f|=1.
  \end{cases}
\end{equation}
Thus it is enough to estimate $D_eD_fH_{n,p}$ below.

Put $\bar G = G_{n,p}-e-f.$ We compare the non-origin $2$-faces of the
order-induced unimodular triangulations of
$P_{\bar G}, P_{\bar G+e}, P_{\bar G+f}, P_{\bar G+e+f}, $ and classify their
changes below.

For a graph $G$, let $\mathcal F(G)$ denote the set of non-origin $2$-faces in
the order-induced unimodular triangulation of the symmetric edge polytope $P_G$.
For $\{i,j\}=\{e,f\}$, let $S_i$ count the $2$-faces in $\mathcal F(\bar G+i)$
that are not in $\mathcal F(\bar G)$, and let $R_i^{(j)}$ count the $2$-faces in
$\mathcal F(\bar G)$ that disappear after adding $i$ but remain after adding
$j$.
Let $R_{e,f}$ count the 2-faces in $\mathcal F(\bar G)$ that disappear whenever
either $e$ or $f$ is added.
Further, let $S_{e,f}$ count the $2$-faces that belong to
$\mathcal F(\bar G+e+f)$ but to none of $\mathcal F(\bar G)$,
$\mathcal F(\bar G+e)$, and $\mathcal F(\bar G+f)$.
Conversely, let $T_{e,f}$ count the $2$-faces that belong to
$\mathcal F(\bar G)$, $\mathcal F(\bar G+e)$, and $\mathcal F(\bar G+f)$, but
not to $\mathcal F(\bar G+e+f)$.
Finally, let $T_i^{(j)}$ count the $2$-faces that belong to
$\mathcal F(\bar G+i)$ but to neither $\mathcal F(\bar G)$ nor
$\mathcal F(\bar G+i+j)$.
It easily follows that
\begin{align*}
  (H_{n,p})_{e,f}^{+,+}
  &=(H_{n,p})_{e,f}^{-,-}+S_e+S_f+S_{e,f} -\bigl(R_e^{(f)}+R_f^{(e)}+R_{e,f}\bigr)\\
  &\quad-\bigl(T_e^{(f)}+T_f^{(e)}+T_{e,f}\bigr),\\
  (H_{n,p})_{e,f}^{+,-} &=(H_{n,p})_{e,f}^{-,-}+S_e-\bigl(R_e^{(f)}+R_{e,f}\bigr),\\
  (H_{n,p})_{e,f}^{-,+} &=(H_{n,p})_{e,f}^{-,-}+S_f-\bigl(R_f^{(e)}+R_{e,f}\bigr).
\end{align*}
Substituting these identities into the mixed difference gives
\begin{align*}
  D_eD_fH_{n,p}
  =&pq\bigl((H_{n,p})_{e,f}^{+,+}-(H_{n,p})_{e,f}^{+,-}-(H_{n,p})_{e,f}^{-,+}+(H_{n,p})_{e,f}^{-,-}\bigr)\\
  =&pq\bigl(S_{e,f}+R_{e,f}-T_e^{(f)}-T_f^{(e)}-T_{e,f}\bigr).
\end{align*}
Taking absolute values gives
\begin{equation}\label{eq:l313-absolute-bound}
  |D_eD_fH_{n,p}| \le pq\bigl(S_{e,f}+R_{e,f}+T_e^{(f)}+T_f^{(e)}+T_{e,f}\bigr).
\end{equation}
We estimate these five terms separately.
For upper bounds, we do not need to consider the order conditions and count
every 3-, 4-, 5-, or 6-cycle under consideration as an obstruction.
We can also ignore the positions and orientations of the arcs, since there are
only finitely many possible choices.
Throughout the proof, $C>0$ denotes a constant that may change from line to line
and from case to case, but is independent of $n$, $p$, $e$, and $f$.

Recall that $E_{n,p}$ denotes the number of arcs of $G_{n,p}$.
Fix $e,f$, and put
\[
  W(k):=\max_{\{u,w\}\in\binom{e\cup f}{2}} W_{uw}(k),
\]
where $W_{uw}(k)$ denotes the number of simple $u$--$w$ paths of length $k$ in
$\bar G$.
The maximum ranges over six endpoint pairs when $e,f$ are disjoint, and three
endpoint pairs when they are adjacent.

\medskip
\noindent\textit{Bounding $S_{e,f}$.}
The three directed arcs generating a $2$-face counted by $S_{e,f}$ can be
represented as $\{\vec e,\vec f,\vec x\},$ where $x$ is an arc of $\bar G$.
Hence,
\begin{equation}\label{eq:l313-S-bound}
  S_{e,f}\lesssim E_{n,p}.
\end{equation}

\medskip
\noindent\textit{Bounding $T_f^{(e)}$ and $T_e^{(f)}$.}
The three directed arcs generating a $2$-face counted by $T_f^{(e)}$ can be
represented as $\{\vec f,\vec x_1,\vec x_2\},$ where $x_1,x_2$ are arcs of
$\bar G$.
It disappears after adding $e$, so there is an obstruction cycle containing $e$.
By Proposition~\ref{prop:two-face-criterion}, there are three possible types of
configurations:
\begin{enumerate}
\renewcommand{\labelenumi}{(\roman{enumi})}
\item All three arcs $f,x_1,x_2$ lie on an obstruction cycle of length
$4$, $5$, or $6$ containing $e$.

\item The arcs $f,x_1$ lie on an obstruction cycle of length $3$ or $4$ containing
$e$, while $x_2$ is not an arc of this cycle.

\item The arcs $x_1,x_2$ lie on an obstruction cycle of length $3$ or $4$
containing $e$, while $f$ is not an arc of this cycle.
\end{enumerate}

First suppose that $|e\cap f|=1$, and write $e=\{s,v\}$ and $f=\{v,t\}$.
In case (i), deleting $e,f$ leaves an $s$--$t$ path of length $2$, $3$, or $4$
containing $x_1,x_2$.
In case (ii), deleting $e,f$ leaves an $s$--$t$ path of length $1$ or $2$
containing $x_1$, while $x_2$ can be chosen arbitrarily from the arcs of
$\bar G$.
In case (iii), deleting $e$ leaves a path of length $2$ or $3$ containing
$x_1,x_2$.
Thus, the numbers of configurations in cases (i), (ii), and (iii) are bounded by
$C\bigl(W(2)+W(3)+W(4)\bigr), CE_{n,p}\bigl(1+W(2)\bigr), C\bigl(W(2)+W(3)\bigr),$
respectively.

Next suppose that $|e\cap f|=0$, and write $e=\{a,b\}$ and $f=\{c,d\}$.
In case (i), deleting $e,f$ breaks the relevant cycle into two simple paths in
$\bar G$.
Since the cycle has length $4$, $5$, or $6$, the possible pairs of path lengths
are $(1,1),(1,2),(1,3),(2,2)$, up to interchange.
In case (ii), the relevant cycle can only be a $4$-cycle.
Deleting $e,f$ leaves two length-$1$ paths with fixed endpoints, while $x_2$ can
be chosen arbitrarily from the arcs of $\bar G$.
In case (iii), deleting $e$ leaves a path of length $2$ or $3$ containing
$x_1,x_2$.
Thus, the numbers of configurations in cases (i), (ii), and (iii) are bounded by
$C\bigl(1+W(2)+W(3)+W(2)^2\bigr), CE_{n,p}, C\bigl(W(2)+W(3)\bigr),$
respectively.

Combining the preceding bounds, we obtain
\begin{equation}\label{eq:l313-T-single-bound}
  T_f^{(e)} \lesssim
  \begin{cases}
    1+E_{n,p}+W(2)+W(3)+W(2)^2, & |e\cap f|=0,\\
    E_{n,p}+E_{n,p}W(2)+W(2)+W(3)+W(4), & |e\cap f|=1.
  \end{cases}
\end{equation}
By symmetry, the same bound holds for $T_e^{(f)}$.

\medskip
\noindent\textit{Bounding $T_{e,f}$.}
The three directed arcs generating a $2$-face counted by $T_{e,f}$ can be
represented as $\{\vec x_1,\vec x_2,\vec x_3\}$, where $x_1,x_2,x_3$ are arcs of
$\bar G$.
It disappears only after both $e$ and $f$ are added, so there is an obstruction
cycle containing both $e$ and $f$.
By Proposition~\ref{prop:two-face-criterion}, there are two possible types:
\begin{enumerate}
\renewcommand{\labelenumi}{(\roman{enumi})}
\item All three arcs $x_1,x_2,x_3$ lie on an obstruction cycle of length
$5$ or $6$ containing $e$ and $f$.

\item The arcs $x_1,x_2$ lie on an obstruction cycle of length $4$ containing $e$ and $f$,
while $x_3$ is not an arc of this cycle.
\end{enumerate}

First suppose that $|e\cap f|=1$, and write $e=\{s,v\}$ and $f=\{v,t\}$.
In case (i), deleting $e,f$ leaves an $s$--$t$ path of length $3$ or $4$
containing all three arcs $x_1,x_2,x_3$.
In case (ii), deleting $e,f$ leaves an $s$--$t$ path of length $2$ containing
$x_1,x_2$, while $x_3$ can be chosen arbitrarily from the arcs of $\bar G$.
Thus, the numbers of configurations in cases (i) and (ii) are bounded by
$C\bigl(W(3)+W(4)\bigr), CE_{n,p}W(2),$ respectively.

Next suppose that $|e\cap f|=0$, and write $e=\{a,b\}$ and $f=\{c,d\}$.
In case (i), deleting $e,f$ breaks the relevant cycle into two simple paths in
$\bar G$.
The possible pairs of path lengths are $(1,2),(1,3),(2,2)$, up to interchange.
In case (ii), the relevant cycle is a $4$-cycle.
Deleting $e,f$ leaves two length-$1$ paths with fixed endpoints, while $x_3$ can
be chosen arbitrarily from the arcs of $\bar G$.
Thus, the numbers of configurations in cases (i) and (ii) are bounded by
$C\bigl(W(2)+W(3)+W(2)^2\bigr), CE_{n,p},$ respectively.

Combining the preceding bounds, we obtain
\begin{equation}
  T_{e,f} \lesssim
  \begin{cases}
    E_{n,p}+W(2)+W(3)+W(2)^2, & |e\cap f|=0,\\
    E_{n,p}W(2)+W(3)+W(4), & |e\cap f|=1.
  \end{cases}
\end{equation}

\medskip
\noindent\textit{Bounding $R_{e,f}$.}
The three directed arcs generating a $2$-face counted by $R_{e,f}$ can be
represented as $\{\vec{x}_1,\vec{x}_2,\vec{x}_3\}$, where $x_1,x_2,x_3$ are arcs
of $\bar G$.
Let $\Gamma_e$ and $\Gamma_f$ denote the obstruction cycles for $e$ and $f$,
respectively, where $\Gamma_e$ contains $e$ but not $f$, and $\Gamma_f$ contains
$f$ but not $e$.
By Proposition~\ref{prop:two-face-criterion}, each of $\Gamma_e$ and $\Gamma_f$ is either a
$5$- or $6$-cycle containing all three arcs $x_1,x_2,x_3$, or a $3$- or
$4$-cycle containing two of them.
The pair $(\Gamma_e,\Gamma_f)$ is therefore of one of the following four types:
\begin{enumerate}
\renewcommand{\labelenumi}{(\roman{enumi})}
\item[(i)] Both cycles contain all three arcs $x_1,x_2,x_3$.
\item[(ii)] One cycle contains all three arcs, while the other contains $x_1,x_2$.
\item[(iii)] One cycle contains $x_1,x_2$, while the other contains $x_2,x_3$.
\item[(iv)] Both cycles contain $x_1,x_2$. In this case, the
arcs $x_1$ and $x_2$ must share a node; otherwise $x_1$ and $x_2$ would lie on a
unique directed $4$-cycle, so $\Gamma_e$ and $\Gamma_f$ could not be distinct.
\end{enumerate}

We shall use the following observation.
Among the directed 2-paths with the same pair of endpoints, at most one can
generate an edge in the triangulation of $P_{\bar G}$, and among the directed
3-paths with the same pair of endpoints, at most one can generate a 2-face.
Indeed, two distinct such paths of the same length would give rise to an
obstruction cycle.
It therefore remains to count the configurations arising from $\Gamma_e$ and
$\Gamma_f$.

The representative configurations that account for the different terms in the
bound are shown in Figure~\ref{fig:R-ef-configurations}.
\begin{figure}[!b]
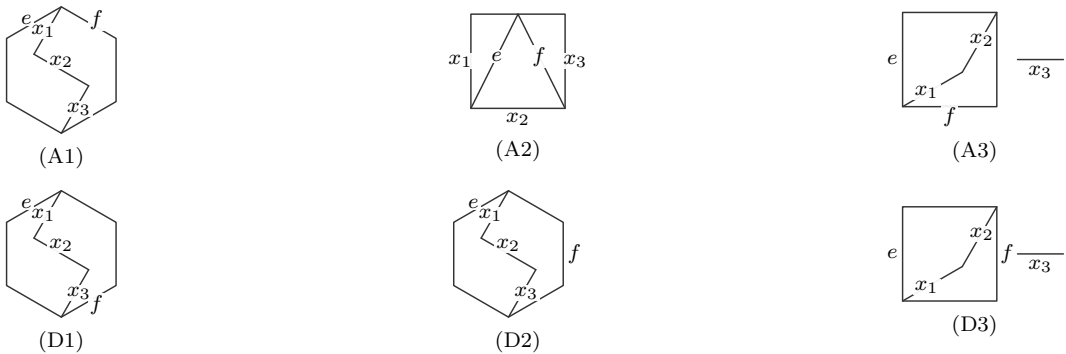

\centering

\RefConf{\RefAdjOne}{\textnormal{(A1)}} \hfill
\RefConf{\RefAdjTwo}{\textnormal{(A2)}} \hfill
\RefConf{\RefAdjThree}{\textnormal{(A3)}}

\vspace{2.5mm}

\RefConf{\RefDisOne}{\textnormal{(D1)}} \hfill
\RefConf{\RefDisTwo}{\textnormal{(D2)}} \hfill
\RefConf{\RefDisThree}{\textnormal{(D3)}}

\caption{Representative configurations used to bound $R_{e,f}$.
Configurations (A1)--(A3) correspond to the case $|e\cap f|=1$, while (D1)--(D3)
correspond to the case $|e\cap f|=0$.}
\label{fig:R-ef-configurations}
\end{figure}

First suppose that $|e\cap f|=1$.
In configuration (A1), the remaining arcs form a path of length four between
fixed endpoints.
Hence the number of such configurations is bounded by $CW(4)$.
Including the corresponding cases with paths of smaller length, the total number
is bounded by $C\bigl(1+W(2)+W(3)+W(4)\bigr).$ Configuration (A2) leaves two
length-two paths between fixed endpoint pairs.
Hence the number of such configurations is bounded by $CW(2)^2$.
In configuration (A3), $x_3$ can be chosen arbitrarily, while $x_1,x_2$ form a
length-two path.
Hence the number of such configurations is bounded by $CE_{n,p}W(2)$.
Including the corresponding cases with paths of smaller length, the total number
is bounded by $CE_{n,p}\bigl(1+W(2)\bigr).$

Next suppose that $|e\cap f|=0$.
Configuration (D1) leaves two length-two paths between fixed endpoint pairs.
Hence the number of such configurations is bounded by $CW(2)^2$.
In configuration (D2), the remaining arcs form a length-three path between fixed
endpoints.
Hence the number of such configurations is bounded by $CW(3)$.
Including the corresponding cases with paths of smaller length, the total number
is bounded by $C\bigl(1+W(2)+W(3)\bigr).$ Finally, in configuration (D3), $x_3$
can be chosen arbitrarily, while the remaining choices are determined by their
ordered endpoints.
Hence the number of such configurations is bounded by $CE_{n,p}$.

The remaining configurations are handled in the same way as above.
In particular, the symmetric difference of $\Gamma_e$ and $\Gamma_f$
may consist of two triangles, whose contributions are bounded by
$CW(2)^2$ or $C(1+E_{n,p})$, depending on whether $e$ and $f$
belong to different triangles or the same triangle.
In all cases, the path lengths and the number of free arc choices
are no larger than those in the representative configurations
in Figure~\ref{fig:R-ef-configurations}. Hence, the numbers of these configurations are
bounded by the terms already obtained above.

Combining the preceding bounds, we obtain
\begin{equation}
  \label{eq:l313-R-bound}
  R_{e,f} \lesssim
  \begin{cases}
    1+E_{n,p}+W(2)+W(3)+W(2)^2, & |e\cap f|=0,\\[2mm]
    1+E_{n,p}+E_{n,p}W(2) +W(2)+W(3)+W(2)^2+W(4), & |e\cap f|=1.
  \end{cases}
\end{equation}

\smallskip
\noindent\emph{Fourth moments.}
Substituting (\ref{eq:l313-S-bound})--(\ref{eq:l313-R-bound}) into
(\ref{eq:l313-absolute-bound}), we obtain
\begin{equation}\label{eq:l313-combined-bound}
  |D_eD_fH_{n,p}|\lesssim
  \begin{cases}
    pq+pqE_{n,p}+pqW(2)+pqW(3)+pqW(2)^2, & |e\cap f|=0,\\[1mm]
    \begin{aligned}[t]
      &pq+pqE_{n,p}+pqE_{n,p}W(2)+pqW(2)\\
      &\quad+pqW(3)+pqW(2)^2+pqW(4),
    \end{aligned}
    & |e\cap f|=1.
  \end{cases}
\end{equation}
By the definition of $W(k)$, for every fixed $r\ge1$,
\[
  \mathbb E[W(k)^r] \le C\,\mathbb E[W_{n,p}(k)^r].
\]
Applying the moment estimates from the proof of Lemma~\ref{lem:fourth-moment-DeFnp}
to (\ref{eq:l313-combined-bound}), we obtain
\begin{equation}\label{eq:l313-H-fourth-moment}
  \mathbb E[(D_eD_fH_{n,p})^4] \lesssim
  \begin{cases}
    n^8p^8q^4, & |e\cap f|=0,\\[1mm]
    n^8p^8q^4+n^{12}p^{16}q^4, & |e\cap f|=1.
  \end{cases}
\end{equation}
Finally, combining (\ref{eq:l313-K-fourth-moment}) and
(\ref{eq:l313-H-fourth-moment}) with $D_eD_fF_{n,p}=D_eD_fH_{n,p}+D_eD_fK_{n,p}$
proves (\ref{eq:second-gradient-F-fourth-moment}), as desired.
\end{proof}

\section{Proof of \texorpdfstring{Theorem~\ref{thm:main}}{Theorem \ref{thm:main}}}\label{sec:clt}
\begin{proof}
Part~(i) follows directly from Theorems~\ref{thm:expectation} and~\ref{thm:variance}.
It remains to prove part~(ii).
For $1\le i\le5$, define $B_i(F_{n,p})$ by the following formulas:
\begin{align*}
  B_1(F_{n,p}) &=\sum_{e,f,g} \E[(D_eF_{n,p})^2(D_fF_{n,p})^2]^{1/2}\\
  &\qquad\times \E[(D_gD_eF_{n,p})^2(D_gD_fF_{n,p})^2]^{1/2},\\
  B_2(F_{n,p}) &=\frac{1}{pq}\sum_{e,f,g} \E[(D_gD_eF_{n,p})^2(D_gD_fF_{n,p})^2],\\
  B_3(F_{n,p}) &=\frac{1}{pq}\sum_e\E[(D_eF_{n,p})^4],\\
  B_4(F_{n,p})
  &=\frac{1}{pq}\sum_{e,g} \E[(D_eF_{n,p})^4]^{1/2} \E[(D_gD_eF_{n,p})^4]^{1/2},\\
  B_5(F_{n,p}) &=\frac{1}{(pq)^2}\sum_{e,g}\E[(D_gD_eF_{n,p})^4].
\end{align*}
Here every index ranges over $\binom{[n]}2$.
These are the five quantities arising from Proposition~\ref{prop:MS} when it is applied to
the standardized face count.
Applying Proposition~\ref{prop:MS} to
\[
  \frac{F_{n,p}-\E F_{n,p}}{\sqrt{\Var(F_{n,p})}},
\]
and using the linearity of the first- and second-order discrete gradients, we
obtain
\begin{equation}\label{eq:MS-unnormalized-Fnp}
  d_K\left( \frac{F_{n,p}-\E F_{n,p}}{\sqrt{\Var(F_{n,p})}},N \right)
  \lesssim \frac{1}{\Var(F_{n,p})} \sum_{i=1}^5\sqrt{B_i(F_{n,p})}.
\end{equation}

It remains to bound the five terms on the right-hand side of
(\ref{eq:MS-unnormalized-Fnp}).
To this end, recall from Lemma~\ref{lem:fourth-moment-DeFnp} that
$\mathbb{E}\left[(D_e F_{n,p})^4\right] \lesssim n^{16}p^{10}q^2$, and hence, by
the Cauchy--Schwarz inequality, for any $e,f \in \binom{[n]}{2}$,
\[
  \mathbb{E}\left[ (D_eF_{n,p})^2(D_fF_{n,p})^2 \right]
  \le \sqrt{ \mathbb{E}\left[(D_eF_{n,p})^4\right] \mathbb{E}\left[(D_fF_{n,p})^4\right]
  }
  \lesssim n^{16}p^{10}q^2.
\]
Also, Lemma~\ref{lem:second-order-gradient-F} gives
\[
  \mathbb{E}\left[(D_eD_fF_{n,p})^4\right] \lesssim
  \begin{cases}
    n^8p^8q^4, & |e\cap f|=0,\\[2mm]
    n^8p^8q^4+n^{12}p^{16}q^4, & |e\cap f|=1.
  \end{cases}
\]
Another application of the Cauchy--Schwarz inequality gives, for any
$e,f,g\in\binom{[n]}{2}$,
\[
  \mathbb{E}\left[ (D_gD_eF_{n,p})^2(D_gD_fF_{n,p})^2 \right] \lesssim
  \begin{cases}
    n^8p^8q^4, & |g\cap e|=|g\cap f|=0,\\[2mm]
    n^8p^8q^4+n^{10}p^{12}q^4, & |g\cap e|+|g\cap f|=1,\\[2mm]
    n^8p^8q^4+n^{12}p^{16}q^4, & |g\cap e|=|g\cap f|=1.
  \end{cases}
\]

We start with the estimate for $B_1(F_{n,p})$.
Since $D_gD_eF_{n,p}=0$ whenever $g=e$, only terms with $g\neq e,f$ contribute.
According to the three cases above, the numbers of triples $(e,f,g)$ satisfying
$|g\cap e|=|g\cap f|=0, |g\cap e|+|g\cap f|=1, |g\cap e|=|g\cap f|=1$ are
$O(n^6)$, $O(n^5)$, and $O(n^4)$, respectively.
Combining these counts with the moment bounds gives
\begin{align*}
  B_1(F_{n,p}) &\lesssim n^6(n^8p^5q)(n^4p^4q^2) +n^5(n^8p^5q)(n^4p^4q^2+n^5p^6q^2)\\
  &\qquad +n^4(n^8p^5q)(n^4p^4q^2+n^6p^8q^2) \lesssim n^{18}p^9q^3.
\end{align*}
Similarly, for $B_2$,
\begin{align*}
  B_2(F_{n,p})
  &\lesssim \frac{ n^6n^8p^8q^4 +n^5(n^8p^8q^4+n^{10}p^{12}q^4)
  +n^4(n^8p^8q^4+n^{12}p^{16}q^4) }{pq}\\
  &\lesssim n^{14}p^7q^3+n^{16}p^{15}q^3.
\end{align*}
The term $n^{15}p^{11}q^3$ is bounded by the geometric mean of the two displayed
terms and may therefore be absorbed.

For $B_3$, the first-order fourth-moment estimate gives directly
\[
  B_3(F_{n,p}) \lesssim \frac{n^2}{pq}\,n^{16}p^{10}q^2 \lesssim n^{18}p^9q.
\]
For $B_4$ and $B_5$, the numbers of ordered pairs $(e,g)$ with $|e\cap g|=0$ and
$|e\cap g|=1$ are $O(n^4)$ and $O(n^3)$, respectively.
Hence
\begin{align*}
  B_4(F_{n,p})
  &\lesssim \frac{ n^4(n^8p^5q)(n^4p^4q^2) +n^3(n^8p^5q)(n^4p^4q^2+n^6p^8q^2) }{pq}\\
  &\lesssim n^{16}p^8q^2+n^{17}p^{12}q^2,\\
  B_5(F_{n,p}) &\lesssim \frac{ n^4n^8p^8q^4 +n^3(n^8p^8q^4+n^{12}p^{16}q^4) }{(pq)^2}\\
  &\lesssim n^{12}p^6q^2+n^{15}p^{14}q^2.
\end{align*}
Combining the preceding estimates, we obtain
\begin{align}
  B_1(F_{n,p})&\lesssim n^{18}p^9q^3,
  &\qquad B_2(F_{n,p})&\lesssim n^{14}p^7q^3+n^{16}p^{15}q^3, \notag\\
  B_3(F_{n,p})&\lesssim n^{18}p^9q,
  &\qquad B_4(F_{n,p})&\lesssim n^{16}p^8q^2+n^{17}p^{12}q^2, \notag\\
  B_5(F_{n,p})&\lesssim n^{12}p^6q^2+n^{15}p^{14}q^2.
  \label{eq:five-B-Fnp}
\end{align}
Taking square roots in (\ref{eq:five-B-Fnp}) and using
$\sqrt{a+b}\le \sqrt a+\sqrt b$, we obtain
\[
  \begin{aligned}
    \sqrt{B_1(F_{n,p})}
    &\lesssim n^9p^{9/2}q^{3/2},
    &\qquad \sqrt{B_2(F_{n,p})}
    &\lesssim n^7p^{7/2}q^{3/2} +n^8p^{15/2}q^{3/2},\\
    \sqrt{B_3(F_{n,p})}
    &\lesssim n^9p^{9/2}q^{1/2},
    &\qquad \sqrt{B_4(F_{n,p})}
    &\lesssim n^8p^4q+n^{17/2}p^6q,\\
    \sqrt{B_5(F_{n,p})} &\lesssim n^6p^3q+n^{15/2}p^7q.
  \end{aligned}
\]
Since $p\gg n^{-2}$ and $q$ is bounded away from zero, all the preceding terms
are bounded by a constant multiple of $n^9p^{9/2}q^{1/2}$.
Summing over $i$ gives
\[
  \sum_{i=1}^5\sqrt{B_i(F_{n,p})} \lesssim n^9p^{9/2}q^{1/2}.
\]
By Theorem~\ref{thm:variance}, we have
\[
  \Var(F_{n,p}) = \frac18 (n)_{10}p^5q^5H_3(p)^2+o(n^{10}p^5),
\]
and $H_3(p)>2$.
Therefore
\[
  \Var(F_{n,p})\asymp n^{10}p^5q^5.
\]
Substituting the last two estimates into (\ref{eq:MS-unnormalized-Fnp}) yields
\[
  d_K\left( \frac{F_{n,p}-\E F_{n,p}}{\sqrt{\Var(F_{n,p})}},N \right)
  \lesssim \frac{1}{n\sqrt p\,q^{9/2}}
  \lesssim\frac{1}{n\sqrt p}.
\]
This proves the Berry--Esseen bound.
Since $p\gg n^{-2}$, the central limit theorem follows as well.
\end{proof}

\section{Higher-dimensional face counts}
\label{sec:higher-dimensional}
In this section, we prove Theorem~\ref{thm:k-face-clt}, retaining the triangulation
$\Delta_{n,p}$ and the notation $F^{[k]}_{n,p}=f_k(\Delta_{n,p})$ introduced in
Section~\ref{sec:intro}.
For $r\ge0$, let $H^{[r]}_{n,p}$ denote the number of non-origin $r$-faces.
Adjoining the origin gives a bijection between non-origin $(k-1)$-faces and
$k$-faces containing the origin.
Hence
\begin{equation}\label{eq:k-cone-decomposition}
  F^{[k]}_{n,p}=H^{[k]}_{n,p}+H^{[k-1]}_{n,p},\qquad k\ge1.
\end{equation}
This notation agrees with that of the preceding sections:
\begin{equation}\label{eq:k-H01}
  \begin{gathered}
    F^{[1]}_{n,p}=K_{n,p},\qquad F^{[2]}_{n,p}=F_{n,p},\\
    H^{[0]}_{n,p}=2E_{n,p},\qquad H^{[1]}_{n,p}=K_{n,p}-2E_{n,p},\qquad
    H^{[2]}_{n,p}=H_{n,p}.
  \end{gathered}
\end{equation}
In particular, (\ref{eq:k-cone-decomposition}) reduces to
Lemma~\ref{lem:two-face-decomposition} when $k=2$.

The proof begins with a local characterization of non-origin faces in terms of
obstructions arising from directed cycles.
The Gr\"obner-basis description in
\cite[Proposition~3.8]{HigashitaniJochemkoMichalek2019} yields the following
generalization of Proposition~\ref{prop:two-face-criterion} to arbitrary dimension.

\begin{proposition}\label{prop:k-face-criterion}
Let $G$ be a graph equipped with a total order $\prec$ on its arcs, and let
$\mathcal T_{\prec}$ be the corresponding order-induced regular unimodular
triangulation of $P_G$, with the origin adjoined as a vertex.
Let $\alpha$ be a set of $k+1$ directed arcs of $G$.
The corresponding vertices span a non-origin $k$-face of $\mathcal T_{\prec}$ if
and only if all of the following conditions hold:
\begin{enumerate}[label=\textup{(\roman*)}]
\item No two arcs of $\alpha$ have the same underlying undirected arc;
\item For every $2\le j\le k+1$, no $j$ arcs in $\alpha$
      lie on a directed cycle of length $2j-1$;
\item For every $2\le j\le k+1$ and every $j$-element subset
      $\beta\subseteq\alpha$ whose arcs lie on a directed cycle $C$ of length
      $2j$, the $\prec$-minimal arc of $C$ belongs to $\beta$.
\end{enumerate}
\end{proposition}
For the remainder of this section, we decompose $H^{[k]}_{n,p}$ according to
whether the $k+1$ arcs generating each face are pairwise node-disjoint:
\begin{equation}\label{eq:k-HDR-decomposition}
  H^{[k]}_{n,p}=D^{[k]}_{n,p}+R^{[k]}_{n,p},
\end{equation}
where $D^{[k]}_{n,p}$ counts the non-origin $k$-faces generated by $k+1$
pairwise node-disjoint arcs, and $R^{[k]}_{n,p}$ counts the remaining non-origin
$k$-faces.

Let $\mathcal M^{[k]}_1$ be the collection of sets of $k+1$ directed arcs on
$[n]$ whose underlying arcs are pairwise node-disjoint.
For $\alpha\in\mathcal M^{[k]}_1$, let $I_\alpha$ be the indicator that $\alpha$
generates a non-origin $k$-face of $\Delta_{n,p}$.
Then $D^{[k]}_{n,p} =\sum_{\alpha\in\mathcal M^{[k]}_1} I_\alpha.$

\subsection{Expectation}

We first determine the order of the expectation.
\begin{proposition}\label{prop:k-expectation}
Under the assumptions of Theorem~\ref{thm:k-face-clt},
\begin{equation}\label{eq:k-expectation-order}
  \mathbb E F_{n,p}^{[k]} \asymp_k n^{2k+2}p^{k+1}.
\end{equation}
\end{proposition}

\begin{proof}Fix
$\alpha=\{a_i\to b_i:0\le i\le k\}\in\mathcal M_1^{[k]},$ where $a_i\to b_i$
denotes the directed arc from $a_i$ to $b_i$.
If $\alpha$ generates a non-origin $k$-face, then all $k+1$ arcs of $\alpha$
must be present.
On the other hand, since the arcs of $\alpha$ are pairwise node-disjoint,
Proposition~\ref{prop:k-face-criterion} implies that the arcs outside $\alpha$ that may
affect the value of $I_\alpha$ are those joining $a_i$ and $b_j$ with $i\ne j$.
The number of such arcs is $k(k+1)$.
Thus, if all arcs of $\alpha$ are present and all these additional arcs are
absent, then $\alpha$ generates a non-origin $k$-face.
Hence we have
\[
  p^{k+1}q^{k(k+1)} \le \mathbb E I_\alpha \le p^{k+1}.
\]
Since $q$ is bounded away from zero and the number of configurations in
$\mathcal M_1^{[k]}$ is $(n)_{2k+2}/(k+1)!\asymp_k n^{2k+2}$, we obtain
\begin{equation}\label{eq:k-D-expectation}
  \mathbb E D_{n,p}^{[k]} \asymp_k n^{2k+2}p^{k+1}.
\end{equation}
Every configuration counted by $R_{n,p}^{[k]}$ uses $k+1$ distinct underlying
arcs and at most $2k+1$ nodes.
Hence
\[
  \mathbb E R_{n,p}^{[k]} \lesssim_k n^{2k+1}p^{k+1} =o\!\left(n^{2k+2}p^{k+1}\right).
\]
Together with (\ref{eq:k-D-expectation}), this gives
\[
  \mathbb E H_{n,p}^{[k]} \asymp_k n^{2k+2}p^{k+1}.
\]
Similarly,
\[
  \mathbb E H_{n,p}^{[k-1]} \asymp_k n^{2k}p^k =o\!\left(n^{2k+2}p^{k+1}\right),
\]
Combining the preceding estimates with (\ref{eq:k-cone-decomposition}) completes
the proof.
\end {proof}

\subsection{Variance estimates}
\label{subsec:k-variance-lower}
As in the proof of Proposition~\ref{prop:variance-H1}, assign independent weights
$w_e\sim U(0,1)$ to all possible underlying arcs, independently of $G_{n,p}$,
and order the arcs by their weights.
Let $w$ denote the collection of these weights.
For simplicity, we write $\Delta_{n,p}$ for the triangulation induced by the
order determined by $w$.
Since $H^{[0]}_{n,p}=2E_{n,p}$ by (\ref{eq:k-H01}), it follows inductively from
(\ref{eq:k-cone-decomposition}) that $H^{[r]}_{n,p}$ does not depend on $w$ for
any $r\ge0$.

Recall the decomposition in (\ref{eq:k-HDR-decomposition}) and the indicator
representation $D^{[k]}_{n,p}=\sum_{\alpha\in\mathcal M^{[k]}_1} I_\alpha$.
As in (\ref{eq:H1-overlap-decomposition}), we expand the conditional variance
and group the covariance terms according to the size of the node overlap:
\begin{equation}\label{eq:k-overlap-decomposition}
  \E\bigl[\Var(D^{[k]}_{n,p}\mid w)\bigr]
  =\sum_{\ell=0}^{2k+2} \sum_{\substack{\alpha,\beta\in\mathcal M^{[k]}_1\\
  |V(\alpha)\cap V(\beta)|=\ell}} \E\bigl[\Cov(I_\alpha,I_\beta\mid w)\bigr].
\end{equation}

For $\alpha\in\mathcal M^{[k]}_1$, no $j$ arcs in $\alpha$ can lie on a directed
cycle of length $2j-1$.
If a directed cycle of length $2j$ contains $j$ arcs of $\alpha$, then the
remaining $j$ arcs join the heads of these arcs to the tails of the others.
Thus, once $w$ is fixed, whether $\alpha$ generates a face depends on the
presence or absence of these additional arcs.

If $|V(\alpha)\cap V(\beta)|\le1$, the corresponding families of underlying arcs
are disjoint.
Hence $I_\alpha$ and $I_\beta$ are conditionally independent given $w$, and the
corresponding terms in (\ref{eq:k-overlap-decomposition}) vanish.

If $|V(\alpha)\cap V(\beta)|=2$, write
$V(\alpha)\cap V(\beta)=\{u,v\}, s=\{u,v\}.$ The conditional covariance vanishes
unless $s$ is relevant to both configurations.
Conditional on $I_s$ and $w_s$, the remaining arc indicators that can affect the
value of $I_\alpha$ are then independent of those that can affect the value of
$I_\beta$.
Hence, the conditioning argument in the proof of Proposition~\ref{prop:variance-H1}, based
on Lemma~\ref{lem:variance}, gives
\[
  \E\bigl[\Cov(I_\alpha,I_\beta\mid w)\bigr]
  = \E\bigl[\Cov(I_\alpha,I_\beta\mid w_s)\bigr].
\]

Let $I_s$ be the indicator that $s$ is present.
As in Section~\ref{sec:variance}, define
\begin{equation}\label{eq:k-fa-fb}
  \begin{aligned}
    f_a(t)
    &:=\E[I_\alpha\mid I_s=1,w_s=t] -\E[I_\alpha\mid I_s=0,w_s=t],
    && s\in\alpha,\\
    f_b(t)
    &:=\E[I_\alpha\mid I_s=1,w_s=t] -\E[I_\alpha\mid I_s=0,w_s=t],
    && s\notin\alpha.
  \end{aligned}
\end{equation}
By symmetry, the definitions of $f_a$ and $f_b$ do not depend on the particular
choice of $\alpha$ and $s$, so both are well defined.
The conditional independence above, together with (\ref{eq:cov-factor-full-en}),
yields
\begin{equation}\label{eq:k-covariance-factorization}
  \E[\Cov(I_\alpha,I_\beta\mid w)]
  = pq\int_0^1 f_i(t)f_j(t)\,\dd t,
  \qquad i,j\in\{a,b\},
\end{equation}
where $i=a$ when $s\in\alpha$ and $i=b$ otherwise, with $j$ defined similarly
for $\beta$.
This yields three cases according to the membership of $s$ in $\alpha$ and
$\beta$.

\smallskip\noindent
\textit{Case 1: $s\in\alpha\cap\beta$.} The union of the underlying arcs
consists of $2k+1$ pairwise node-disjoint arcs.
Choose $4k+2$ nodes from $[n]$ and partition them into $2k+1$ unordered pairs.
Choose the common arc $s$ in $2k+1$ ways, and choose the $k$ remaining arcs of
$\alpha$ from the other $2k$ arcs in $\binom{2k}{k}$ ways; the rest belong to
$\beta$.

The $2k$ non-common arcs may be oriented independently, and the common
underlying arc may be oriented independently in $\alpha$ and $\beta$.
Thus there are $2^{2k}\cdot 2^2=2^{2k+2}$ possible orientations in total.
The contribution of this case is therefore
\begin{align*}
  &\binom{n}{4k+2} \frac{(4k+2)!}{2^{2k+1}(2k+1)!} (2k+1)\binom{2k}{k}2^{2k+2}
  pq\int_0^1 f_a(t)^2\,\dd t\\
  &\qquad= \frac{2\fall n{4k+2}}{(k!)^2}\, pq\int_0^1 f_a(t)^2\,\dd t.
\end{align*}

\smallskip\noindent
\textit{Case 2: $s$ belongs to exactly one of $\alpha$ and $\beta$.} The union
of the underlying arcs consists of a path of length three, with middle arc $s$,
together with $2k-1$ isolated arcs.
Choose $4k+2$ nodes and partition them into one four-node block and $2k-1$
unordered pairs.
Within the four-node block, choose $s$ in $\binom{4}{2}$ ways and join each of
the two remaining nodes to a different endpoint of $s$ in $2$ ways.
There are $2$ choices for which of $\alpha$ and $\beta$ contains $s$.
The configuration containing $s$ contains $k$ of the $2k-1$ isolated arcs,
giving $\binom{2k-1}{k}$ choices; the other contains the two end arcs of the
path and the remaining isolated arcs.

The two indicators can be dependent only when the two end arcs of the 3-path are
oriented toward one of $u,v$ and away from the other.
There are $2$ such joint orientations.
The middle arc $s$ and the $2k-1$ isolated arcs may be oriented freely, giving
$2\cdot2^{2k}=2^{2k+1}$ possible orientations in total.
Hence the contribution of this case is
\begin{align*}
  &\binom{n}{4k+2} \frac{(4k+2)!}{4!\,2^{2k-1}(2k-1)!} \binom{4}{2}\cdot2\cdot2
  \binom{2k-1}{k}2^{2k+1}\\
  &\qquad\qquad{}\times pq\int_0^1 f_a(t)f_b(t)\,\dd t\\
  &\qquad= \frac{4k\fall n{4k+2}}{(k!)^2}\, pq\int_0^1 f_a(t)f_b(t)\,\dd t.
\end{align*}

\smallskip\noindent
\textit{Case 3: $s\notin\alpha\cup\beta$.} The union of the underlying arcs
consists of two node-disjoint paths of length two, centered at $u$ and $v$,
together with $2k-2$ isolated arcs.
Choose the nodes and partition them into two unordered three-node blocks and
$2k-2$ unordered pairs.
Choose the middle node of each three-node block in $3^2$ ways.
On each path, assign one arc to $\alpha$ and the other to $\beta$, giving $2^2$
assignments.
Each configuration then contains $k-1$ of the isolated arcs, and there are
$\binom{2k-2}{k-1}$ ways to assign the isolated arcs to the two configurations.

The two indicators can be dependent only when the two path arcs in each
configuration are oriented toward one of $u$ and $v$ and away from the other.
The four path arcs then have only $2^2$ possible orientations, while the $2k-2$
isolated arcs have $2^{2k-2}$ possible orientations.
Thus there are $2^{2k}$ possible orientations in total.
The contribution is therefore
\begin{align*}
  &\binom{n}{4k+2} \frac{(4k+2)!} {(3!)^2\,2!\,2^{2k-2}(2k-2)!} 3^2\cdot2^2
  \binom{2k-2}{k-1}2^{2k}\\
  &\qquad\qquad{}\times pq\int_0^1 f_b(t)^2\,\dd t\\
  &\qquad= \frac{2k^2\fall n{4k+2}}{(k!)^2}\, pq\int_0^1 f_b(t)^2\,\dd t.
\end{align*}

Summing the three cases gives
\begin{equation}\label{eq:k-two-node-overlap}
  \sum_{\substack{\alpha,\beta\in\mathcal M^{[k]}_1\\
  |V(\alpha)\cap V(\beta)|=2}} \E[\Cov(I_\alpha,I_\beta\mid w)] = \frac{2\fall
  n{4k+2}}{(k!)^2}\,pq \int_0^1\bigl(f_a(t)+kf_b(t)\bigr)^2\,\dd t.
\end{equation}
Thus it remains to evaluate $f_a(t)+kf_b(t)$.

\begin{lemma}\label{lem:k-matching-kernel}
The sum $f_a(t)+k f_b(t)$ does not depend on $t$ and is a polynomial in $p$
satisfying
\begin{equation}\label{eq:k-kernel-bounds}
  \frac{p^k q^k}{k!} \leq f_a(t)+k f_b(t) \leq p^k, \qquad 0<p<1.
\end{equation}
Moreover, $D_{n,p}^{[k]}$ and $R_{n,p}^{[k]}$ do not depend on $w$.
\end{lemma}

\begin{proof}
Recall that a perfect matching between two disjoint node sets $A$ and $B$ of
equal cardinality is a set of pairwise node-disjoint arcs, each joining a node
of $A$ to a node of $B$, that covers every node of $A\cup B$.
For two matchings $M$ and $M'$, their symmetric difference is
$M\triangle M'=(M\setminus M')\cup(M'\setminus M).$ The symmetric difference of
two perfect matchings between $A$ and $B$ is a disjoint union of even cycles.

Fix disjoint node sets $A=\{a_0,\ldots,a_k\}, B=\{b_0,\ldots,b_k\},$ and let
$G_{n,p}[A,B]$ be the bipartite graph consisting of the arcs of $G_{n,p}$
joining $A$ to $B$, all oriented from $A$ to $B$.
Let $\mathcal P(A,B)$ denote the set of all $(k+1)!$ possible perfect matchings
between $A$ and $B$, with the same orientation.
We show that whenever $G_{n,p}[A,B]$ contains a perfect matching, exactly one
such matching generates a non-origin $k$-face.

Order its perfect matchings lexicographically: for distinct matchings $M$ and
$M'$, declare $M$ to precede $M'$ if the $\prec$-minimal underlying arc in
$M\triangle M'$ belongs to $M$.
If $G_{n,p}[A,B]$ contains a perfect matching, let $M^*$ be the first of its
perfect matchings in this order.
If $M^*$ did not generate a non-origin $k$-face, then by
Proposition~\ref{prop:k-face-criterion}, some $j$ arcs of $M^*$ would lie on a directed
cycle $\Gamma$ of length $2j$ whose $\prec$-minimal underlying arc does not
belong to $M^*$.
The other $j$ underlying arcs of $\Gamma$ which are not in $M^*$, together with
the arcs of $M^*$ outside $\Gamma$, form another perfect matching $M'$ between
$A$ and $B$, oriented from $A$ to $B$.
Since the $\prec$-minimal underlying arc of $\Gamma$ belongs to $M'$, this
matching precedes $M^*$, a contradiction.
Thus $M^*$ generates a non-origin $k$-face.

On the other hand, for any other perfect matching $M$ contained in
$G_{n,p}[A,B]$, the $\prec$-minimal underlying arc in $M\triangle M^*$ belongs
to $M^*$.
Hence $M$ can not generate a face by Proposition~\ref{prop:k-face-criterion}(iii).
Thus we have
\begin{equation}\label{eq:k-perfect-matching-indicator}
  \sum_{\alpha\in\mathcal P(A,B)} I_\alpha
  = \mathbf{1}\{G_{n,p}[A,B]\text{ contains a perfect matching}\}.
\end{equation}

Since each set of pairwise node-disjoint directed arcs uniquely determines its
tail and head sets, summing (\ref{eq:k-perfect-matching-indicator}) over all
ordered pairs $(A,B)$ shows that $D_{n,p}^{[k]}$ does not depend on $w$.
The same holds for $R_{n,p}^{[k]}$ by (\ref{eq:k-HDR-decomposition}).

Now set $s=\{a_0,b_0\}$.
We say that $s$ is pivotal if $G_{n,p}[A,B]+s$ contains a perfect matching,
whereas $G_{n,p}[A,B]-s$ does not.
Among the $(k+1)!$ matchings in $\mathcal P(A,B)$, exactly $k!$ contain $s$,
while $k\,k!$ do not.
Taking the difference between the conditional expectations of
(\ref{eq:k-perfect-matching-indicator}) given $I_s=1$, $w_s=t$ and $I_s=0$,
$w_s=t$, and using the definitions of $f_a$ and $f_b$, we obtain
\begin{equation}\label{eq:k-pivotal-probability}
  k!\bigl(f_a(t)+k f_b(t)\bigr) = \mathbb{P}\{s\text{ is pivotal}\}.
\end{equation}
Since this event depends only on which of the other arcs are present, its
probability is a polynomial in $p$ and does not depend on $t$.

For the lower bound, consider the event that, for every $1\leq i\leq k$, the arc
$\{a_i,b_i\}$ is present and $\{a_0,b_i\}$ is absent.
On this event, $a_0$ is isolated in $G_{n,p}[A,B]-s$, whereas $G_{n,p}[A,B]+s$
contains the diagonal perfect matching.
Thus $s$ is pivotal, and this event has probability $p^kq^k$.
For the upper bound, if $s$ is pivotal, then $G_{n,p}[A,B]+s$ contains a perfect
matching using $s$.
There are $k!$ possible such matchings, each requiring $k$ other arcs to be
present.
Hence, by a union bound, $\mathbb{P}\{s\text{ is pivotal}\}\leq k!p^k.$ Together
with (\ref{eq:k-pivotal-probability}), this proves (\ref{eq:k-kernel-bounds}).
\end{proof}

\begin{lemma}\label{lem:k-variance-lower}
Under the assumptions of Theorem~\ref{thm:k-face-clt},
\begin{equation}\label{eq:k-D-variance-asymptotic}
  \Var(D^{[k]}_{n,p})
  =\frac{2\fall n{4k+2}}{(k!)^2}\,pq \bigl(f_a(0)+kf_b(0)\bigr)^2
  +o\bigl(n^{4k+2}p^{2k+1}\bigr),
\end{equation}
and
\begin{equation}\label{eq:k-R-variance-negligible}
  \Var(R^{[k]}_{n,p})=o\bigl(n^{4k+2}p^{2k+1}\bigr).
\end{equation}
Consequently,
\begin{equation}\label{eq:k-F-variance-order}
  \Var(F^{[k]}_{n,p}) \asymp_k n^{4k+2}p^{2k+1}.
\end{equation}
\end{lemma}
\begin{proof}
We first complete the estimate of $\operatorname{Var}(D_{n,p}^{[k]})$.
By (\ref{eq:k-overlap-decomposition}) and (\ref{eq:k-two-node-overlap}), it
remains to bound the terms with $|V(\alpha)\cap V(\beta)|\ge3$.
To this end, we count the nodes and the distinct underlying arcs in the joint
configuration.
If $|V(\alpha)\cap V(\beta)|=\ell$, the joint configuration uses $4k+4-\ell$
nodes and at least $2k+2- \lfloor\ell/2\rfloor$ arcs.
Thus we have
\[
  \sum_{\substack{\alpha,\beta\in\mathcal M^{[k]}_1\\
  |V(\alpha)\cap V(\beta)|=\ell}} \bigl|\E[\Cov(I_\alpha,I_\beta\mid w)]\bigr|
  \lesssim_k n^{4k+4-\ell}p^{2k+2-\lfloor\ell/2\rfloor}.
\]
The ratio of the preceding bound to $n^{4k+2}p^{2k+1}$ is $(n^2p)^{1-m}$ if
$\ell=2m\ge4$, and $n^{-1}(n^2p)^{1-m}$ if $\ell=2m+1\ge3$.
Since $n^2p\to\infty$, both ratios tend to zero.
Hence all terms with $|V(\alpha)\cap V(\beta)|\ge3$ are $o(n^{4k+2}p^{2k+1})$.
Moreover, $D_{n,p}^{[k]}$ is independent of $w$, and therefore
$\Var(D_{n,p}^{[k]})= \E\!\left[\Var(D_{n,p}^{[k]}\mid w)\right].$ Combining
this identity with (\ref{eq:k-overlap-decomposition}),
(\ref{eq:k-two-node-overlap}), and Lemma~\ref{lem:k-matching-kernel} yields
(\ref{eq:k-D-variance-asymptotic}).

We next estimate $\operatorname{Var}(R_{n,p}^{[k]})$.
Every configuration counted by $R_{n,p}^{[k]}$ contains a pair of intersecting
arcs.
Fix two configurations $\alpha$ and $\beta$ contributing to $R_{n,p}^{[k]}$, and
suppose that they have no common underlying arc.
We first consider the case $|V(\alpha)\cap V(\beta)|=0$.
If neither $\alpha$ nor $\beta$ contains a directed $2$-path, then $I_\alpha$
and $I_\beta$ are conditionally independent.
Hence, by symmetry, we may assume that $\alpha$ contains a directed $2$-path.
Fix such a directed $2$-path in $\alpha$, and let $B_\alpha$ be the event that
no directed $4$-cycle containing this directed $2$-path and a node outside
$V(\alpha)\cup V(\beta)$ triggers an obstruction.
Writing $r=n-|V(\alpha)\cup V(\beta)|$, by the estimate used in the proof of
Lemma~\ref{lem:H-expectation},
\[
  \mathbb P(B_\alpha)
  = \frac{1}{r+1}\sum_{j=0}^{r}(1-p^2)^j
  = \frac{1-(1-p^2)^{r+1}}{(r+1)p^2}
  \le \min\left\{1,\frac{C_k}{np^2}\right\}
\]
Since the joint configuration uses $2k+2$ distinct underlying arcs, we have
\[
  \begin{aligned}
    \left|\mathbb E[\operatorname{Cov}(I_\alpha,I_\beta\mid w)]\right|
    &\le \max\left\{ \mathbb E(I_\alpha I_\beta), \mathbb E\!\left[
    \mathbb E(I_\alpha\mid w)\mathbb E(I_\beta\mid w) \right] \right\} \\
    &\le p^{2k+2}\mathbb P(B_\alpha)\le p^{2k+2}\min\left\{1,\frac{C_k}{np^2}\right\}\\
    &\lesssim p^{2k+1}n^{-1/2}.
  \end{aligned}
\]
Also the joint configuration uses at most $4k+2$ nodes.
Hence the contribution is at most of order $n^{4k+3/2}p^{2k+1}$.
We next consider the case $|V(\alpha)\cap V(\beta)|\ge 1$.
The joint configuration uses at most $4k+1$ nodes and $2k+2$ distinct underlying
arcs.
Hence the contribution is at most of order $n^{4k+1}p^{2k+2}$.

It remains to consider the case where $\alpha$ and $\beta$ have $m\ge1$ common
underlying arcs.
Their joint configuration uses $2k+2-m$ distinct underlying arcs.
If the common arcs are not pairwise node-disjoint, they use at most $2m-1$
nodes, while the non-common arcs use at most $4(k+1-m)$ additional nodes.
Hence the number of nodes in the joint configuration is at most
\[
  (2m-1)+4(k+1-m)=4k+3-2m,
\]
and the contribution is at most of order $n^{4k+3-2m}p^{2k+2-m}$.
If the common arcs are pairwise node-disjoint, they use $2m$ nodes.
Since both $\alpha$ and $\beta$ are counted by $R_{n,p}^{[k]}$, the non-common
arcs in each configuration use at most $2(k+1-m)-1$ additional nodes.
Hence the number of nodes in the joint configuration is at most
\[
  2m+2\bigl(2(k+1-m)-1\bigr)=4k+2-2m,
\]
and the contribution is at most of order $n^{4k+2-2m}p^{2k+2-m}$.
Combining the preceding estimates, we obtain
\[
  \begin{aligned}
    \Var(R_{n,p}^{[k]}) &=\E\!\left[\Var(R_{n,p}^{[k]}\mid w)\right] \\
    &\lesssim_k n^{4k+3/2}p^{2k+1}+n^{4k+1}p^{2k+2} +n^{4k}p^{2k+1} +\sum_{m=2}^{k+1}
    n^{4k+3-2m}p^{2k+2-m}\\
    &=o\bigl(n^{4k+2}p^{2k+1}\bigr),
  \end{aligned}
\]
which proves (\ref{eq:k-R-variance-negligible}).
Combining (\ref{eq:k-D-variance-asymptotic}),
(\ref{eq:k-R-variance-negligible}), and the Cauchy--Schwarz inequality gives
\[
  \Var(H_{n,p}^{[k]})
  = \Var(D_{n,p}^{[k]}) +o\!\left(n^{4k+2}p^{2k+1}\right)
  \asymp_k n^{4k+2}p^{2k+1}.
\]
The same argument with $k-1$ in place of $k$ yields
\[
  \Var(H_{n,p}^{[k-1]}) \asymp_k n^{4k-2}p^{2k-1} =o\!\left(n^{4k+2}p^{2k+1}\right).
\]
Finally, (\ref{eq:k-cone-decomposition}) and Cauchy--Schwarz give
\[
  \Var(F_{n,p}^{[k]})
  = \Var(H_{n,p}^{[k]}) +o\!\left(n^{4k+2}p^{2k+1}\right)
  \asymp_k n^{4k+2}p^{2k+1}.
\]
This completes the proof.
\end{proof}

\subsection{Discrete gradients and normal approximation}
\label{subsec:k-normal-approximation}
As in Section~\ref{sec:discrete-gradients}, when estimating the first- and
second-order discrete gradients, we ignore the order conditions and count every
cycle under consideration as an obstruction for upper bounds.
We also suppress the choices of orientations and positions of the arcs.

\begin{lemma}\label{lem:k-first-order-gradients}
Fix $k\ge3$ and let $e=\{s,t\}\in\binom{[n]}2$.
Under the assumptions of Theorem~\ref{thm:k-face-clt}, we have
\begin{equation}\label{eq:k-first-gradient-moment}
  \mathbb E\!\left[(D_eH_{n,p}^{[k]})^4\right] \lesssim_k n^{8k}p^{4k+2}q^2.
\end{equation}
The same estimate holds with $F_{n,p}^{[k]}$ in place of $H_{n,p}^{[k]}$.
\end{lemma}
\begin{proof}
Fix $e=\{s,t\}$, and let $W_{n,p}(\ell)$ denote the number of $s$--$t$ paths of
length $\ell$.
We distinguish the non-origin $k$-faces created and destroyed by adding $e$.
For each such face, let $\alpha$ be the set of $k+1$ directed arcs generating
it.
If the face is newly created, then $\alpha$ contains an orientation of $e$, and
the remaining $k$ arcs in $\alpha$ can be chosen freely.
Hence the number of such faces is at most $C_kE_{n,p}^{\,k}$.

Now suppose that the face disappears after adding $e$.
By Proposition~\ref{prop:k-face-criterion}, for some $2\le j\le k+1$, $j$ arcs in $\alpha$
must lie on an obstruction cycle containing $e$, whose length is between $j+1$
and $2j$.
Removing $e$ leaves an $s$--$t$ path of length between $j$ and $2j-1$, while the
remaining $k+1-j$ arcs in $\alpha$ may be chosen freely.
Therefore,
\begin{equation}\label{eq:k-first-gradient-count}
  \frac{|D_eH_{n,p}^{[k]}|}{\sqrt{pq}}
  \lesssim_k E_{n,p}^k + \sum_{j=2}^{k+1} \left( \sum_{\ell=j}^{2j-1}W_{n,p}(\ell)
  \right) E_{n,p}^{k+1-j}.
\end{equation}

As in the proof of Lemma~\ref{lem:fourth-moment-DeFnp}, for every fixed $\ell\ge2$
and $m\ge1$,
\[
  \mathbb E[W_{n,p}(\ell)^m]
  \lesssim_{\ell,m} n^{\ell-1}p^\ell + \bigl(n^{\ell-1}p^\ell\bigr)^m.
\]
The binomial moment bounds and the Cauchy--Schwarz inequality therefore yield,
for $j\le\ell\le2j-1$,
\[
  \mathbb E\!\left[ W_{n,p}(\ell)^4E_{n,p}^{4(k+1-j)} \right]
  \lesssim_k (n^2p)^{4k},
  \qquad \mathbb E[E_{n,p}^{4k}]
  \lesssim_k (n^2p)^{4k}.
\]
Taking fourth moments in (\ref{eq:k-first-gradient-count}) easily gives
\[
  \mathbb E\!\left[(D_eH_{n,p}^{[k]})^4\right] \lesssim_k n^{8k}p^{4k+2}q^2,
\]
which proves (\ref{eq:k-first-gradient-moment}).
Finally, by (\ref{eq:k-cone-decomposition}), applying the same estimate in
dimensions $k$ and $k-1$ proves the corresponding bound for $F_{n,p}^{[k]}$.
\end{proof}
\begin{lemma}
\label{lem:second-order-gradient-Fk}
\label{lem:k-gradient-moments}
Fix $k\ge3$ and let $e,f\in\binom{[n]}2$ be distinct arcs.
Under the assumptions of Theorem~\ref{thm:variance}, we have
\begin{equation}\label{eq:second-gradient-Hk-fourth-moment}
  \E[(D_eD_fH_{n,p}^{[k]})^4] \lesssim_k
  \begin{cases}
    n^{8k-8}p^{4k}q^4, & |e\cap f|=0,\\[1mm]
    n^{8k-8}p^{4k}q^4+n^{8k-4}p^{4k+8}q^4, & |e\cap f|=1.
  \end{cases}
\end{equation}
The same estimate holds with $F_{n,p}^{[k]}$ in place of $H_{n,p}^{[k]}$.
\end{lemma}

\begin{proof}
As in the proof of Lemma~\ref{lem:second-order-gradient-F}, the decomposition in
(\ref{eq:l313-absolute-bound}) yields
\[
  |D_eD_fH_{n,p}^{[k]}| \le pq\bigl(S_{e,f}+T_e^{(f)}+T_f^{(e)}+T_{e,f}+R_{e,f}\bigr).
\]
We estimate these five terms separately.

\medskip
\noindent\emph{Bounding $S_{e,f}$.}
Let $\alpha$ be the set of $k+1$ directed arcs generating a $k$-face counted by
$S_{e,f}$.
Then $\alpha$ contains an orientation of each of $e$ and $f$, and the remaining
$k-1$ arcs in $\alpha$ can be chosen freely.
Hence
\begin{equation}\label{eq:l56-S-bound}
  S_{e,f}\lesssim_k(E_{n,p})^{k-1}.
\end{equation}
For the remaining terms, we count the free nodes, namely the additional nodes
that must be chosen to specify a contributing configuration.

\medskip
\noindent\emph{Bounding $T_f^{(e)}$ and $T_e^{(f)}$.}
Let $\alpha$ be the set of $k+1$ directed arcs generating a $k$-face counted by
$T_f^{(e)}$.
Then $\alpha$ contains an orientation of $f$.
By Proposition~\ref{prop:k-face-criterion}, there is an obstruction cycle containing $e$.

If the cycle also contains $f$, let $j$ be the number of other arcs in $\alpha$
on it.
Its length is between $j+2$ and $2j+2$, and the remaining $k-j$ arcs in $\alpha$
can be chosen freely.
Deleting $e$ and $f$ leaves one path if $|e\cap f|=1$ and two paths if
$|e\cap f|=0$.
Since their endpoints are fixed and their total length is at most $2j$, the
number of free nodes is at most $2(k-j)+2j-2+|e\cap f| =2k-2+|e\cap f|.$

If the cycle does not contain $f$, let $j$ be the number of arcs in $\alpha$ on
it.
Its length is between $j+1$ and $2j$, and the remaining $k-j$ arcs in $\alpha$
other than the orientation of $f$ can be chosen freely.
Deleting $e$ leaves a path with fixed endpoints and length at most $2j-1$.
Hence the number of free nodes is at most $2(k-j)+2j-2=2k-2.$ The same estimates
hold for $T_e^{(f)}$ by symmetry.

\medskip
\noindent\emph{Bounding $T_{e,f}$.}
Let $\alpha$ be the set of $k+1$ directed arcs generating a $k$-face counted by
$T_{e,f}$.
By Proposition~\ref{prop:k-face-criterion}, $e$, $f$, and $j$ arcs in $\alpha$ lie on an
obstruction cycle of length between $j+2$ and $2j$.
The remaining $k+1-j$ arcs in $\alpha$ can be chosen freely.
Deleting $e$ and $f$ from the cycle leaves one path if $|e\cap f|=1$ and two
paths if $|e\cap f|=0$.
Since their endpoints are fixed by $e$ and $f$ and their total length is at most
$2j-2$, the number of free nodes is at most
$2(k+1-j)+2j-4+|e\cap f| =2k-2+|e\cap f|.$

\medskip
\noindent\emph{Bounding $R_{e,f}$.}
Let $\alpha$ be the set of $k+1$ directed arcs generating a $k$-face counted by
$R_{e,f}$.
By Proposition~\ref{prop:k-face-criterion}, there are two obstruction cycles $\Gamma_e$ and
$\Gamma_f$, where $\Gamma_e$ contains $e$ but not $f$, and $\Gamma_f$ contains
$f$ but not $e$.

As in the proof of Lemma~\ref{lem:second-order-gradient-F}, among directed paths with
the same pair of endpoints and length, at most one can be formed from arcs
generating a $k$-face.
Thus, its internal nodes do not contribute to the count of free nodes.

Let $j$ be the number of arcs in $\alpha$ outside $\Gamma_e\cup\Gamma_f$.
Suppose that $\Gamma_e$ and $\Gamma_f$ share $m_1$ directed arcs in $\alpha$ and
$n_1$ other arcs.
Let $a$ be the number of directed paths formed by these $m_1$ arcs, and let $b$
be the number of path components formed by all $m_1+n_1$ common arcs.
Let $m_2$ and $n_2$ be the numbers of remaining arcs of $\Gamma_e$ in and
outside $\alpha$, respectively, and define $m_3$ and $n_3$ similarly for
$\Gamma_f$.
Then we have
\begin{align*}
  \begin{cases}
    &m_1+m_2+m_3=k+1-j,\\
    &n_1+n_2\le m_1+m_2,\\
    &n_1+n_3\le m_1+m_3,\\
    &a\le n_1+b.
  \end{cases}
\end{align*}

The last inequality follows since deleting the $n_1$ common arcs creates at most
$n_1$ additional path components.
Deleting the common arcs from both $\Gamma_e$ and $\Gamma_f$ leaves a graph that
can be decomposed into arc-disjoint cycles of total length $m_2+m_3+n_2+n_3$.
Since the nodes of $e\cup f$ are fixed, the number of free nodes is at most
$m_2+m_3+n_2+n_3-4+|e\cap f|.$ Since the endpoints of the common paths are
already fixed, contracting the $a$ directed subpaths in $\alpha$ leaves $b$
paths of total length $n_1+a$, contributing at most $n_1+a-b$ additional free
nodes.
Combining all the preceding estimates, it follows that the total number of free
nodes is at most
\begin{align*}
  & m_2+n_2+m_3+n_3-4+|e\cap f| +n_1+a-b+2j\\
  \le &m_2+n_2+m_3+n_3+2n_1 -4+2j+|e\cap f|\\
  \le &2(m_1+m_2+m_3)-4+2j+|e\cap f|\\
  =&2k-2+|e\cap f|.
\end{align*}

\noindent\emph{Fourth moments.}
The preceding estimates show that each of the five terms
can be bounded by finite sums involving $E_{n,p}$,
$W(\ell)$, and possibly simple cycles, with at most
$2k-2+|e\cap f|$ free nodes.
Using the corresponding moment estimates and H\"older's
inequality, it suffices to count the free nodes and arcs.
For simple cycles, the same moment estimate applies since
their unions have at least as many arcs as nodes.

A freely chosen arc contributes two free nodes and one arc,
a path of length $\ell$ with fixed endpoints contributes
$\ell-1$ free nodes and $\ell$ arcs, and a simple cycle
of length $\ell$ contributes $\ell$ free nodes and $\ell$ arcs.
Therefore, a term with $2m$ free nodes requires at least
$m$ arcs, while a term with $2m+1$ free nodes requires at
least $m+2$ arcs, since the latter necessarily involves
an even-length path or an odd-length cycle.

If $|e\cap f|=0$, there are at most $2k-2$ free nodes.
A configuration with $2k-2$ free nodes requires at least
$k-1$ arcs, while configurations with fewer free nodes
contribute no more. Thus we have
\[
\mathbb E\left[(D_eD_fH_{n,p}^{[k]})^4\right]
\lesssim_k n^{8k-8}p^{4k}q^4.
\]

If $|e\cap f|=1$, there are at most $2k-1$ free nodes.
Terms with at most $2k-2$ free nodes are controlled as
above. Those with $2k-1$ free nodes require at least
$k+1$ arcs. Thus we have
\[
\mathbb E\left[(D_eD_fH_{n,p}^{[k]})^4\right]
\lesssim_k n^{8k-8}p^{4k}q^4
+n^{8k-4}p^{4k+8}q^4.
\]
This proves (\ref{eq:second-gradient-Hk-fourth-moment}).
Finally, (\ref{eq:k-cone-decomposition}) and
$(a+b)^4\le8(a^4+b^4)$ give the assertion for
$F_{n,p}^{[k]}$ by applying the corresponding bounds
in dimensions $k$ and $k-1$, with
(\ref{eq:l313-H-fourth-moment}) covering dimension
two when $k=3$.
\end{proof}

\begin{proof}[Proof of Theorem~\ref{thm:k-face-clt}]
Part~(i) follows from the expectation estimate (\ref{eq:k-expectation-order}) in
Proposition~\ref{prop:k-expectation} and the variance order in Lemma~\ref{lem:k-variance-lower}.
For part~(ii), it suffices to bound the five quantities $B_i(F_{n,p}^{[k]})$,
$1\le i\le5$, as in Section~\ref{sec:clt}.
Combining the estimates in
Lemmas~\ref{lem:k-first-order-gradients} and~\ref{lem:k-gradient-moments} with the
Cauchy--Schwarz inequality and then summing as in Section~\ref{sec:clt}, we
obtain
\begin{equation}\label{eq:k-MS-terms}
  \begin{aligned}
    B_1(F^{[k]}_{n,p})&\lesssim_k n^{8k+2}p^{4k+1}q^3,\\
    B_2(F^{[k]}_{n,p})&\lesssim_k n^{8k-2}p^{4k-1}q^3+n^{8k}p^{4k+7}q^3,\\
    B_3(F^{[k]}_{n,p})&\lesssim_k n^{8k+2}p^{4k+1}q,\\
    B_4(F^{[k]}_{n,p})&\lesssim_k n^{8k}p^{4k}q^2+n^{8k+1}p^{4k+4}q^2,\\
    B_5(F^{[k]}_{n,p})&\lesssim_k n^{8k-4}p^{4k-2}q^2+n^{8k-1}p^{4k+6}q^2.
  \end{aligned}
\end{equation}
Since $p\gg n^{-2}$ and $q$ is bounded away from zero, all five terms in
(\ref{eq:k-MS-terms}) are bounded by a constant multiple of $n^{8k+2}p^{4k+1}q$.
Taking square roots and summing over $i$ gives
\begin{equation}\label{eq:k-MS-sum}
  \sum_{i=1}^5\sqrt{B_i(F^{[k]}_{n,p})} \lesssim_k n^{4k+1}p^{2k+1/2}q^{1/2}.
\end{equation}
Applying Proposition~\ref{prop:MS} with the normalization and Lemma~\ref{lem:k-variance-lower},
we obtain
\[
  d_K\left( \frac{F^{[k]}_{n,p}-\E F^{[k]}_{n,p}} {\sqrt{\Var(F^{[k]}_{n,p})}},N\right)
  \lesssim_k \frac{n^{4k+1}p^{2k+1/2}q^{1/2}} {n^{4k+2}p^{2k+1}}
  \lesssim_k \frac1{n\sqrt p}.
\]
This proves the Berry--Esseen bound in part~(ii).
Since $n\sqrt p\to\infty$, the central limit theorem follows.
\end{proof}

\section*{Acknowledgements}
Su Z. and Zheng H. were partially supported by the National Natural Science
Foundation of China (Grant Nos.~12271475 and U23A2064).

\end{document}